\documentclass[12pt]{amsart}
\usepackage[margin=1in]{geometry} 
\usepackage{enumerate}
\usepackage{enumitem}
\usepackage[mathscr]{eucal}
\usepackage{color} 
\usepackage[colorlinks=true, linkcolor=red, citecolor=blue, menucolor=black]{hyperref}

\usepackage{amsmath,amsthm, amsfonts, amssymb}
\numberwithin{equation}{subsection}
\usepackage{graphicx}
\usepackage{hyperref}
\usepackage{listings}
\usepackage{enumerate}
\usepackage{centernot}
\usepackage{tikz, tikz-cd}
\usepackage{comment}

\theoremstyle{definition}
\newtheorem{introthm}{Theorem}

\newtheorem{theorem}{Theorem}[section]
\newtheorem{lemma}[theorem]{Lemma}
\newtheorem{corollary}[theorem]{Corollary}
\newtheorem{proposition}[theorem]{Proposition}

\newtheorem{observation}[theorem]{Observation}

\theoremstyle{definition}
\newtheorem{definition}[theorem]{Definition}
\newtheorem{example}[theorem]{Example}
\newtheorem{examples}[theorem]{Examples}

\newtheorem{question}[theorem]{Question}
\newtheorem{remark}[theorem]{Remark}

\theoremstyle{remark}

\newcommand{\supp}{\text{Supp}}
\newcommand{\diam}{\operatorname{diam}}

\newcommand{\norm}[1]{\left\Vert#1\right\Vert}
\newcommand{\id}{\operatorname{id}}
\newcommand{\aut}{\operatorname{Aut}}
\newcommand{\stab}{\operatorname{stab}}
\newcommand{\estab}{\operatorname{estab}}

\newcommand{\sym}{\operatorname{Sym}}

\DeclareMathOperator{\dom}{Dom}
\DeclareMathOperator{\image}{Im}
\newcommand{\SL}{\operatorname{SL}}
\newcommand{\psl}{\operatorname{PSL}}

\DeclareMathOperator{\fin}{Fin}

\newcommand{\floor}[1]{\left\lfloor #1 \right\rfloor}

\DeclareMathOperator{\ind}{Ind}

\title[vNAs from Shift-similar \& Houghton-like Groups]{Property Gamma and McDuff Factors from Shift-similar and Houghton-like Groups}

\author[E. Bashwinger]{Eli Bashwinger}
\address{Department of Mathematics, University of Utah, 155 S 1400 E, Salt Lake City, Utah, 84112, USA}
\email{elibashwinger@outlook.com}

\author[P. DeBonis]{Patrick DeBonis}
\address{TU Delft, EWI/DIAM, P.O.Box 5031, 2600 GA Delft, The Netherlands}
\email{phdebonis@gmail.com}

\begin{document}

\begin{abstract}
\looseness=-1 We initiate the study of von Neumann algebras of shift-similar and Houghton-like groups, introduced by Mallery and Zaremsky as analogs of self-similar and R\"over-Nekrashevych groups, respectively. Given a shift-similar group $G$ of permutations of $\mathbb{N}$, one can form the Houghton-like group $H_n(G)$, a ``mashup" of the classical Houghton group $H_n$ with $G$, which acts on $X_n$, the disjoint union of $n$ copies of $\mathbb{N}$. One main result is that when $G$ and the canonical shift map $\sigma$ admit a joint F\o lner sequence for their action on $\mathbb{N}$, then both group von Neumann algebras $L(G)$ and $L(H_n(G))$ are II$_1$ factors with property Gamma and hence $G$ and $H_n(G)$ are inner amenable. More generally, let $\Gamma$ be an ICC group of permutations preserving the ends of $X_n$. We prove that if the natural near action on $\mathbb{N}$ of the germs of almost permutations arising from $\Gamma$ along some ray is amenable, then $L(\Gamma)$ has property Gamma; if this group of germs is itself amenable, then $L(\Gamma)$ is McDuff. Then we construct a shift-similar group $G$, using property (T) groups with certain actions, such that neither $G$ nor $H_n(G)$ is inner amenable. Finally, we produce extremely broad and diverse sources of shift-similar groups satisfying the joint amenability condition. These examples come from residually finite groups and from wobbling groups of a variety of metric spaces. In particular, we completely characterize all metrics $d$ on $\mathbb{N}$ for which the wobbling group $W(\mathbb{N},d)$ is shift-similar in terms of a coarse geometric property---namely, that $\sigma^{-1}$ is bornologous with respect to $d$. As a consequence, this gives us necessary and sufficient conditions for when the wobbling group of any countable metric space admits a shift-similar embedding.
\end{abstract}

\maketitle

\tableofcontents

\section{Introduction}

\addtocontents{toc}{\protect\setcounter{tocdepth}{0}}
\subsubsection*{\textbf{Shift-Similar and Houghton-like Groups}.}
\addtocontents{toc}{\protect\setcounter{tocdepth}{2}}
\looseness = -1 In \cite{mallery25}, shift-similar and Houghton-like groups were introduced by Mallery and Zaremsky as analogs of self-similar and R\"over-Nekrashevych groups, respectively. Shift-similar groups are certain groups of permutations of $\mathbb{N}$ closed under a family of deletion-shift maps on $\sym (\mathbb{N})$, while self-similar groups are certain groups of automorphisms of the infinite regular $d$-ary tree which are closed under the restriction to the complete subtree rooted at any vertex (after relabeling the subtree). Shift-similar groups are to Houghton groups what self-similar groups are to the Higman-Thompson groups: just as a self-similar group and a Higman–Thompson group can be ``mashed up'' to form a R\"over–Nekrashevych group, a shift-similar group and a Houghton group can be mashed up to form a Houghton-like group.

While Mallery and Zaremsky were inspired by self-similar groups and the R\"over-  Nekrashevych groups, they showed that shift-similar and Houghton-like groups are decisively different.  For example, every finitely generated group embeds into a finitely generated (strongly) shift-similar group. Consequently, there are uncountably many isomorphism classes of finitely generated shift-similar groups (see \cite[Theorem 3.28]{mallery25}), in stark contrast to self-similar groups, of which there are only countably many isomorphism classes of finitely generated self-similar groups (see \cite[Section 1.5.3]{Nekrashevych2005}). Moreover, self-similar groups are necessarily residually finite, while shift-similar groups have no such restriction. In contrast, given that any countable group embeds into $\sym (\mathbb{N})$ via Cayley's theorem, after identifying the underlying set $G$ with $\mathbb{N}$, one can take the shift-similar closure of $G$ to be the smallest shift-similar group in $\sym (\mathbb N)$ (see Section \ref{sec:ShiftSimandHoughtonLike} for a more precise definition). Residual finiteness, on the other hand, prevents any analogous notion of a self-similar closure of an arbitrary countable group. As we explore in this paper, the shift-similar closure of a group is highly dependent on the embedding into $\sym (\mathbb N)$ (compare for example Theorem \ref{mainthm:NonInnerAmean} and Theorem \ref{mainthm:ResidFiniteClosureGamma} below).

Shift-similar groups also have their own structural restrictions, particularly due to the presence of $S_{\infty}$, the group of finitely supported permutations on $\mathbb{N}$, as a subgroup of any infinite shift-similar group. What this implies, for example, is that infinite shift-similar groups cannot have any of the following properties (even virtually): torsion-free, solvable, residually finite, being a p-group, or being hyperbolic (see \cite[Corollary 3.13]{mallery25}). But the presence of $S_{\infty}$ does not prevent us from taking shift-similar closures. 

In this paper, we initiate the study of the analytic properties of shift-similar and Houghton-like groups and how they relate to their resulting group von Neumann algebras. Given a countable discrete group $G$ with infinite conjugacy class (ICC), one can associate a II$_1$ factor $L(G)$ by taking the weak operator topology closure of the group ring $\mathbb C[G]$ inside $\mathcal B (\ell^2(G))$, the bounded operators on the Hilbert space $\ell^2(G)$. II$_1$ factors are a fundamental object of study in the theory of von Neumann algebras and many times inherit key properties of the underlying group (e.g. amenability, property (T), Haagerup Property). 

\subsubsection*{\textbf{Property Gamma and the McDuff Property}}

A II$_1$ factor $M$ has \emph{property Gamma} if it admits a non-trivial central sequence (i.e. $M' \cap M^\omega \neq \mathbb C$ for some ultrapower $\omega$). Property Gamma has played an important role 
in the theory and classification of II$_1$ factors ever since it was introduced by Murray and von Neumann to prove that the free group factor $L(\mathbb F_2)$ is not isomorphic to $L(S_\infty)$ \cite{murray1943rings}. Later, Effros proved in \cite{Effros1975} that if a group factor has property Gamma, then the underlying group must be inner amenable. The converse is not true, as was first proved in \cite{VaesNonGamma2012}. More recently, activity has increased around property Gamma and the stronger McDuff property \cite{mcduff1970central}---which requires two non-commuting, non-trivial central sequences (i.e., $M' \cap M^\omega$ is non-abelian)---as well as non-amenable inner amenable groups and rigidity results related to McDuff group factors. See for example, \cite{tucker2020invariant} \cite{bashwinger2021neumann}, \cite{KerrTucker-Drob2023}, \cite{kerr25}, \cite{AMCOS2025}. 

Our first main result is to prove that there is a wealth of shift-similar and Houghton-like groups which yield property Gamma II$_1$ factors. More precisely, if $\sigma$ is the canonical shift map defined by $\sigma (j) = j+1$, we prove the following.

\begin{introthm}[Theorem \ref{thm:PropertyGamma}]\label{mainthm:Gamma}
     Let $G \leq \sym(\mathbb N)$ be an infinite shift-similar group. If $G$ acts amenably on $\mathbb{N}$, then $L(G)$ and $L(H_1(G))$ have property Gamma. If we further assume $G$ and $\sigma$ act jointly amenably on $\mathbb N$, then $L(H_n(G))$ has property Gamma for $n \ge 1$. 
\end{introthm}

In fact, we prove something much more general about countable subgroup of $\sym_{ends}(X_n)$, the permutations of $X_n$ preserving its ends (see Section \ref{sec:ShiftSimandHoughtonLike} for the definition). For each $k=1, \dots , n$, $\sym_{ends}(X_n)$ admits a homomorphism $\phi_k$ to $\sym^*(\mathbb{N})$, the group of germs of almost permutations of $\mathbb{N}$. The map $\phi_k$ encodes the eventual behavior of a permutation in $\sym_{ends}(X_n)$ along the $k$-th ray. With these ideas, we prove the following.

\begin{introthm}[Theorem \ref{thm:PropertyGamma} and Theorem \ref{thm:McDuff}]\label{mainthm:McDuff}
Let $\Gamma \le \sym_{ends}(X_n)$ be a countable group which is ICC (e.g., it contains $S_{\infty}(X_n))$. If the natural near action of $\phi_k(\Gamma)$ on $\mathbb{N}$ is amenable for some $k=1, \dots , n$, then $L(\Gamma)$ has property Gamma. If $\phi_k(\Gamma)$ is amenable, then $L(\Gamma)$ is McDuff.  
\end{introthm}

Theorem \ref{mainthm:Gamma} and Theorem \ref{mainthm:McDuff} should be compared with the property Gamma (and later McDuff) factors from \cite{KerrTucker-Drob2023} and \cite{kerr25}, given that the work of Kerr and Tucker-Drob and Kerr and Petrakos proves property Gamma and the McDuff property for certain groups acting amenably on the Cantor space. In contrast, Theorem \ref{mainthm:Gamma} can be interpreted as a result for shift-similar groups acting amenably on a finite product of copies of the natural numbers (see Theorem \ref{thm:amenablenearaction}). More specifically, \cite[Theorem C]{kerr25} proves that if a group $\Gamma$ acts amenably on a set $Y$ with no finite orbits and $H$ is any non-trivial ICC group, then the generalized wreath product $H \wr_{Y} \Gamma$ yields a type II$_1$ McDuff factor.  It is therefore natural to ask whether their results subsume or overlap with ours. However, in Theorem \ref{thm:NonWreathProduct}, we prove neither shift-similar nor Houghton-like groups can be decomposed as generalized wreath products. Since shift-similar groups are very general, we suspect that the McDuff property cannot always be obtained even under our amenable near action hypothesis.

Our approach to Theorem \ref{mainthm:Gamma} is inspired by Jolissaint's proof that Thompson's group $F$ yields a McDuff factor (see \cite{jolissaint1998central}). In some ways, the classical Houghton groups are like Thompson's group $F$, although they are amenable---even elementary amenable---whereas $F$ is not elementary amenable, and it remains famously open whether $F$ is amenable. Given their similarity to $F$, and the fact that $H_n$ is elementary amenable, we conjectured that Houghton-like groups are inner amenable perhaps by virtue of yielding property Gamma factors. This is indeed the case when $G$ and $\sigma$ act amenably on $\mathbb{N}$, and to our surprise such shift-similar groups turn out to be abundant.

Jolissaint's approach is to view $F$ as acting on the Cantor space $\mathfrak{C} := \{0,1\}^\mathbb{N}$ of infinite binary words by homeomorphisms that preserve the natural lexicographic order and act locally by finite-prefix replacements. Then he considered the normal subgroup $D$ of elements of $F$ fixing a neighborhood of the infinite word $111 \dots$. The quotient $F/D$ is isomorphic to $\mathbb{Z}$, and the associated short exact sequence splits, giving us the semidirect product decomposition $F \cong D \rtimes \mathbb{Z}$. This translates to a crossed product decomposition of the von Neumann algebra: $L(F) \cong L(D) \rtimes \mathbb{Z}$. Then he showed $L(D)$ is McDuff and that the action of the amenable group $\mathbb{Z}$ on $L(D)$ is centrally free, from which it follows $L(F)$ is a McDuff factor. His proof, however, does not explicitly construct the central sequences in $L(F)$.

Working by analogy, the natural counterpart of $D$ for Houghton-like groups is the eventual pointwise stabilizer subgroup $\estab(R_k)$ of the $k$-th ray $R_k$. Apart from our proof of McDuffness, the analogy begins to break down once we look beyond this initial inspiration. When one tries to leverage $\estab (R_k)$ in the way that $D$ is used in the case of $F$, it becomes clear very quickly that different methods are required to prove property Gamma for our group factors.  First, the quotient $H_n(G)/\estab(R_k)$ is much more complicated than $F/D \cong \mathbb{Z}$. In fact, the quotient is generally non-amenable, so the criterion used by Jolissaint does not apply. Second, the short exact sequence does not split and hence we do not get an ordinary crossed product decomposition of $L(H_n(G))$ but rather a twisted crossed product decomposition. Third, the central sequences constructed in $L(D)$ used to prove central freeness are formed from the canonical unitaries, whereas for $L(H_n(G))$ it is not generally possible or straightforward to construct central sequences from canonical unitaries. We must instead use F\o lner sets to form certain (normalized) sums of permutations supported on them.

\subsubsection*{\textbf{Shift-Similar Closures of Residually Finite Groups}}

The shift-similar closure of a group is highly dependent on how it embeds into $\sym (\mathbb N)$. For example, in Section \ref{sec:Non-innerAmean}, we construct an embedding of the property (T) group $\psl_k(\mathbb{Z})$ using its natural action on the rational projective space $P^{k}(\mathbb{Q})$, and then show the shift-similar closure $G$ of $\psl_k(\mathbb{Z})$ with respect to this embedding inherits enough property (T) rigidity such that neither $G$ nor $H_n(G)$ is inner amenable.

\begin{introthm}[Theorem \ref{thm:NonInnerAmeanShift} and Theorem \ref{thm:NonInnerAmean}]\label{mainthm:NonInnerAmean}
    There exist shift-similar groups $G \leq \sym (\mathbb N)$ that do not act amenably on $\mathbb N$ such that neither $G$ nor $H_n(G)$ are inner amenable. As a consequence, neither $L(G)$ nor $L(H_n(G))$ has property Gamma. 
\end{introthm}

The proof of this result crucially relies on two facts: (1) that a property (T) group acts amenably if and only if the action has at least one finite orbit \cite[Lemma 4.2]{GlasnerMonod2007} and (2) that any group containing a relatively ICC property (T) group inherits the non-inner amenability of of the property (T) group (see Lemma \ref{lemma:relativelyICC}).

In light of $L(G)$ and $L(H_n(G))$  not having property Gamma, it is natural to wonder what sort of rigidity properties are possible (e.g. primeness or solidity). For $n \ge 2$, we show $L(H_n(G))$ cannot yield solid factors when $G$ is non-amenable. For these specific non-inner amenable shift-similar groups $G$ constructed from $\psl_k(\mathbb{Z})$, we can show neither $L(G)$ nor $L(H_1(G))$ is solid when $k \ge 4$, but we do not yet know whether they are prime factors. Naturally, for shift-similar groups $G$ in general, we wonder whether $L(G)$ and $L(H_n(G))$ can yield either prime or solid factors. 

In Section \ref{sec:ResidFiniteGrp}, we then use the residual finiteness of $\psl_k(\mathbb{Z})$ to get an action on the associated box space, and then enumerate the box space to get an embedding into $\sym (\mathbb N)$ whose shift-similar closure $G$ is inner amenable---in fact, both $L(G)$ and $L(H_n(G))$ have property Gamma. The following theorem actually applies to any residually finite group. 

\begin{introthm}[Theorem \ref{thm:ResidFiniteClosureGamma}]\label{mainthm:ResidFiniteClosureGamma}
    Every residually finite group admits an embedding with inner amenable shift-similar closure. 
\end{introthm}

The previous two theorems illustrate how the embedding of a group into $\sym(\mathbb N)$ determines which analytic properties are retained in taking the shift-similar closure. The class of shift-similar groups is, therefore, extremely large, rich, and diverse.

\subsubsection*{\textbf{Shift-Similar Wobbling Groups}}

Our basic recipe for producing property Gamma factors from shift-similar and Houghton-like groups is to find large uncountable subgroups $\Gamma$ of $\sym (\mathbb{N})$ which are simultaneously shift-similar and, together with $\sigma$, act amenably on $\mathbb{N}$. Then any countable group $G \le \Gamma$ has shift-similar closure $\overline{G} \le \Gamma$ and hence acts jointly amenably on $\mathbb{N}$ with $\sigma$. By Theorem \ref{mainthm:Gamma}, both $L(\overline{G})$ and $L(H_n(\overline{G}))$ are property Gamma factors and hence both $\overline{G}$ and $H_n(\overline{G})$ are inner amenable. If $\Gamma$ just acts amenably, then we can still conclude that $L(\overline{G})$ and $L(H_1(\overline{G}))$ have property Gamma. 

This strategy is quite effective for wobbling groups. Given a metric space $(X,d)$, the wobbling group is defined to be all of the permutations of the bounded displacement. That is,
\[W(X,d) :=\{g \in \sym(X) : \sup_{x \in X}d(g(x),x) < \infty \}.\]
Wobbling groups have become an important topic in coarse geometry and were used early on to prove non-amenability results, for example, in \cite{DeuberSimonovits1995} and \cite{Laczkovich1990}. More recently, they were used to prove amenability results and study different types of amenable actions (see for example \cite{JuschenkMonod2013} and \cite{Garban2020}).   

The wobbling group on the natural numbers with the euclidean metric appears as a canonical subgroup of $\sym (\mathbb{N})$ that is shift-similar closed. Indeed, applying a shift map to a permutation of bounded displacement still has bounded displacement with respect to the euclidean metric. Moreover, it is not hard to see that $W(\mathbb N, d_{euc})$ acts jointly amenability with $\sigma$ on $\mathbb N$, witnessed using the standard intervals $\{1, \dots , m\}$. Thus, the strategy described above applies to any subgroup of $W(\mathbb N, d_{euc})$. 

Given that any countable metric space is in bijection with $\mathbb N$, two natural questions arise that are interesting in their own right.
The first is, given a countable metric space $(X,d)$, when does $W(X,d)$ act amenably on $X$? When $X$ is a Cayley graph for a group with respect to the word metric, a classical result (see \cite{@paradoxicaldecomp}) actually tells us the amenability of the action is equivalent to the amenability of the group. More generally, for any uniformly discrete metric space of bounded geometry (UDGB) this equivalence also holds (see Lemma \ref{lemma:AmeanMetric/Wobbling}). The second question is:

\begin{question}\label{ques:ShiftSimWobble}
    For a countable metric space $(X,d)$, when does $W(X,d)$ admit a shift-similar embedding into $\sym (\mathbb N)$?
\end{question}

Given how natural the shift-similarity of $W(\mathbb N,d_{euc})$ is, Question \ref{ques:ShiftSimWobble} is much more approachable than asking when an arbitrary countable group is shift-similar, although it is far from trivial. The rest of the paper is devoted to making considerable progress on this question, which is interesting in its own right outside of the application of Theorem \ref{mainthm:Gamma}.

Our first result in this direction actually gives a complete characterization of when $W(X,d)$ is shift-similar in terms of the type of bijection from $\mathbb N$ to $X$ and the relationship to the chosen metric. 

\begin{introthm}[Corollary \ref{cor:WobblingShiftSim}]\label{mainthm:WobblingShiftSim}
Let $\phi : \mathbb{N} \to X$ be an enumeration of a metric space $(X,d)$, and define $T: X \setminus \{\phi (1) \} \to X$ by 
\[
T(x) = \phi \sigma^{-1}\phi(x) = \phi (\phi^{-1}(x)-1).
\]
Then the wobbling group $W(X,d)$ admits a shift-similar embedding via $\phi$ if and only if $T$ is bornologous with respect to $d$.
\end{introthm}

In contrast to the euclidean metric case, there are metrics, even on $\mathbb N$, where $W(X,d)$ is not shift-similar closed (see Section \ref{sec:NonShiftMetricExamples}). 

For a more refined answer to  Question \ref{ques:ShiftSimWobble}, we turn our attention to classifying graphs that satisfying Theorem \ref{mainthm:WobblingShiftSim}. Strikingly, using Sekanina's cube theorem in finite graph theory, we show all locally finite, connected graphs $\mathcal{G} = (V,E)$ with one end admit a bounded-step enumeration $\phi : \mathbb{N} \to V$ of its vertices. This implies $\phi^{-1} W(\mathcal{G}) \phi$ is shift-similar, and if $\mathcal{G}$ is additionally amenable, then $\phi^{-1} W(\mathcal{G}) \phi$ and $\sigma$ act amenably on $\mathbb{N}$. This includes all finitely generated groups with one end. Using this result on one-ended graphs, we show two-ended graphs which are ``coarsely" symmetric also admit a shift-similar embedding, and that the joint amenability condition holds when the graph is amenable. Thus, we have the follow results.

\begin{introthm}[Theorem \ref{thm:OneEndedShiftSim} and Theorem \ref{thm:TwoEndsShiftSim}]\label{mainthm:EndsOfGraphShiftSim}
Let $\mathcal{G} = (V,E)$ be a locally finite, connected graph. If $\mathcal{G}$ is either a one-ended or is a two-ended graph which admits a coarse involution swapping the ends, then there exists an enumeration $\phi : \mathbb{N} \to V$ such that $\phi^{-1} W(\mathcal{G},d_{\mathcal{G}}) \phi \le \sym (\mathbb{N})$ is shift-similar. If $\mathcal{G}$ is amenable, then $\phi^{-1} W(\mathcal{G}, d_{\mathcal{G}}) \phi$ together with $\sigma$ acts amenably on $\mathbb{N}$. Consequently, given any countable group $G \le \phi^{-1} W(\mathcal{G},d_{\mathcal{G}}) \phi$, both $L(\overline{G})$ and $L(H_n(\overline{G}))$ are property Gamma factors by Theorem \ref{mainthm:Gamma}.
\end{introthm} 

\looseness=-1 The initial motivation for studying wobbling groups was a way to produce countable groups whose shift-similar closure satisfies Theorem \ref{mainthm:Gamma}, without actually computing the shift-similar closure. Compared with Theorem \ref{mainthm:ResidFiniteClosureGamma} we are left with several questions that we address is future work.

Which groups admit embeddings into a shift-similar wobbling group $W(X,d)$ that acts amenably on $(X,d)$? Does every finitely generated group admit such an embedding? This would imply every finitely generated group has an inner amenable shift-similar closure. It would be interesting to find a finitely generated group that cannot have \emph{any}  inner amenable shift-similar closure. Finally, there are results that address (to some extent) which groups embed into the wobbling group of different non-amenable groups for which Theorem \ref{mainthm:Gamma} is no longer available (see for example \cite{AFS2025}).

\subsection*{Structure of the Paper}
Apart from the introduction and acknowledgments, the paper is divided into 6 sections. In Section \ref{sec:vNBackground} we give some background on von Neumann algebras, particularly related to the twisted cross product construction that is relevant to the proof of Theorem \ref{mainthm:Gamma}. In Section \ref{sec:ShiftSimandHoughtonLike} we introduce shift-similar groups, Houghton-like groups, and the more general notion of almost permutation groups. In Section \ref{sec:GammaMcDuff} we prove Theorem \ref{mainthm:Gamma} and Theorem \ref{mainthm:McDuff}. In Section \ref{sec:Non-innerAmean} we prove Theorem \ref{mainthm:NonInnerAmean}. In Section \ref{sec:ResidFiniteGrp} we prove Theorem \ref{mainthm:ResidFiniteClosureGamma}. In Section \ref{sec:WobblingGroups} we introduce wobbling groups, the relevant background from coarse geometry, and prove Theorem \ref{mainthm:WobblingShiftSim} and Theorem \ref{mainthm:EndsOfGraphShiftSim}.

\section*{Acknowledgments}
During the preparation of this manuscript, LLMs were utilized solely as an exploratory tool and to assist in proofreading. The authors assume full responsibility for all content.

P. DeBonis was supported in part by the NWO grant OCENW.M.24.252 ‘Quantum Markov semigroups: from geometry to classification’.

\section{Background on von Neumann Algebras}\label{sec:vNBackground}

A von Neumann algebra $M$ is a unital $*$-subalgebra of operators acting on some Hilbert space that is closed in the weak (equivalently, strong) operator topology. We will deal exclusively with separable Hilbert spaces. Every separable von Neumann algebra decomposes as a direct integral of factors---those von Neumann algebras with trivial center, which are the simple objects in the theory. Factors come in three flavors called type I, type II, and type III. These in turn have a further nice classification, but the theory essentially reduces to so-called II$_1$ factors, those factors which are infinite-dimensional and admit a faithful normal tracial state. They fundamentally represent the final frontier in the classification of von Neumann algebras. 

A large and important source of II$_1$ factors comes from countable discrete groups, via their unitary representations and actions on measure spaces, and this is the direction the modern theory of von Neumann algebras has leaned towards. An especially natural unitary representation associated to any group $G$ is the left regular representation on $\ell^2(G)$, which generates the group von Neumann algebra, denoted $L(G)$. It is a classical fact due to Murray and von Neumann that $L(G)$ is a II$_1$ factor if and only if $G$ has the infinite conjugacy class (ICC) property. 

\looseness = -1 Given a von Neumann algebra $M$ with faithful normal tracial state $tr$, we can define the $2$-norm on $M$ as $\norm{x}_2 = \sqrt{tr(x^*x)}$ for $x \in M$. Among other things, this norm can be used to measure when an operator-norm bounded sequence of elements in $M$ asymptotically commutes with every element of $M$. More precisely, a central sequence in $M$ is a sequence $(x_m)_{m \ge 1} \subset M$ uniformly bounded in the operator norm that is not equivalent to a scalar sequence, meaning
\[
\limsup_{m \to \infty} \norm{x_m - tr(x_m)1}_2 \neq 0,
\]
and such that 
\[
\lim_{m \to \infty} \norm{x_m y - y x_m}_2 = 0
\]
for every $y \in M$. A II$_1$ factor $M$  has property Gamma if it has a central sequence, and it is McDuff if it has at least two central sequences which do not asymptotically commute. The McDuff property is equivalent to tensorially absorbing the hyperfinite II$_1$ factor $\mathcal{R}$ (i.e., $M \otimes \mathcal{R} \cong M$) and to the central sequence algebra of $M$ being non-abelian, whereas property Gamma is equivalent to the central sequence algebra being non-trivial \cite{mcduff1970central}. These two properties are related to inner amenability of a group $G$, which is the existence of a conjugation invariant and finitely additive probability measure on $\mathcal{P}(G \setminus \{e\})$. More precisely, in \cite{Effros1975}, Effros proved that if the group von Neumann algebra $L(G)$ has property Gamma, then $G$ is inner amenable. 

Somewhat on the other side of this spectrum are prime and solid factors. A II$_1$ factor $M$ is said to be \emph{prime} if it cannot be decomposed into a tensor product of II$_1$ factors. A strengthening of this is solidity, first introduced by Ozawa in \cite{ozawa2004solid} where he proved that the group von Neumann algebra of a hyperbolic group is solid, answering a question of Ge in \cite{Ge98}. A II$_1$ factor $M$ is \emph{solid} if every diffuse subalgebra $P \subseteq M$ has amenable relative commutant $P' \cap M$. Although primeness and McDuff are completely incompatible, there are examples of property Gamma factors which are also prime, which we regard as exotic and unusual (see \cite{Stalder2006} and \cite{fima2010}).
 
\subsection{Twisted Crossed Product}
We recall the details of a very important way of constructing von Neumann algebras---namely, the twisted crossed product construction. See \cite{Sutherland1980} for a much more extensive discussion of the topic. First, we need to recall the notion of a cocycle action of a group on a von Neumann algebra, then the construction of a twisted crossed product from this cocycle action, and finally how short exact sequences of groups give rise to a twisted crossed product decomposition of a group von Neumann algebra. 

Let $(M, tr)$ be a tracial von Neumann algebra acting on a Hilbert space $\mathcal{H}$, and let $G$ be a countable discrete group. A cocycle action of $G$ on $M$ is a pair $(\alpha, v)$ with $\alpha : G \to \aut (M,tr)$ and $v : G \times G \to \mathcal{U}(M)$ satisfying the following for all $g,h, k \in G$:
\begin{enumerate}
\item $\alpha_{g} \alpha_{h} = ad(v(g,h)) \alpha_{gh}$
\item $v(g,h) v(gh,k) = \alpha_{g}(v(h,k)) v(g,hk)$
\item $v(g,e) = v(e,g) = 1$
\end{enumerate}

From this data, we can form a new von Neumann algebra denoted $M \rtimes_{(\alpha , v)} G$ (we will often drop the subscript $(\alpha ,v)$) defined to be the von Neumann algebra generated by $\pi (M)$ and $\lambda_{v}(G)$. It acts on the Hilbert space $\ell^2(G, \mathcal{H})$ of all functions $\xi : G \to \mathcal{H}$ with $\sum_{g \in G } \norm{\xi (g)}^2  < \infty$. The maps $\lambda_v : G \to B(\ell^2(G, \mathcal{H}))$ and $\pi : M \to B(\ell^2(G, \mathcal{H}))$ are defined as
\[
(\lambda_v(g) \xi )(h) = v(h^{-1},g) \xi (g^{-1} h), \qquad
(\pi (x) \xi )(g) = \alpha_{g^{-1}}(x) \xi (g).
\]
The former is a faithful unitary representation, while the latter is a faithful normal $*$-homomorphism, and they satisfy the commutation relation
\[
\lambda_v (g) \pi (x) = \pi (\alpha_g(x)) \lambda_v(g)
\]
for all $x \in M$ and $g \in G$. That $M \rtimes G$ is generated by $\pi (M)$ and $\lambda_v(G)$ means that we take arbitrary finite alternating products of things from these two sets, then take finite sums, and then finally take the WOT closure of this $*$-subalgebra. Because of the commutation relation, all finite products can be written as $\pi(x) \lambda_v(g)$ for some $x \in M$ and $g \in G$, and therefore the sums become $\sum_{g \in G} \pi(x_g) \lambda_v(g)$. We extend the trace on $M$ to $M \rtimes_{(\alpha , v)} G$: 
\[
tr (\sum_{g \in G} \pi (x_g) \lambda_v(g)) = tr (x_e).
\]
We note that the $*$-isomorphism class of $M \rtimes_{(\alpha , v)} G$ is independent of the choice of separable Hilbert space on which we represent $M$.

We now collect some useful lemmas about central sequences in twisted crossed product factors that we use to prove Theorem \ref{mainthm:Gamma} and Theorem \ref{mainthm:McDuff}. See \cite{Ocneanu1985, Matsumoto1988, Bisch1990, jolissaint1998central} for more on property Gamma and the McDuff property in crossed products, particularly when the acting group is amenable.  

If $(x_m)_{m \ge 1}$ is a sequence in $M$, we say it is \emph{asymptotically invariant} under the cocycle action if 
\[
\lim_{m \to \infty} \norm{\alpha_g(x_m) - x_m}_2 = 0
\]
for every $g \in G$. For a given $g \in G$ we say $\alpha_g$ is \emph{centrally trivial} if 
\[\lim_{m \to \infty} \norm{\alpha_g(x_m) - x_m}_2 = 0\]
for every central sequence $(x_m)_{m \geq 1}$ in $M$. We say the cocycle action of $G$ on $M$ is \emph{centrally free} if $\alpha_g$ is not centrally trivial for all $g \in G \setminus \{e\}$.
 
Our first lemma  gives us a condition on when an asymptotically commutative sequence in a subalgebra lifts to an asymptotically commutative sequence in the twisted crossed product of it. It is most likely well known to experts. 

\begin{lemma}\label{lemma:CentSeqLift1}
Let $S$ be a generating set for some group $G$. If $(x_m)_{m \ge 1}$ is an asymptotically commutative sequence in $M$ that is asymptotically invariant under the cocycle action of the generators $S$, then $(\pi(x_m))_{m \ge 1}$ is an asymptotically commutative sequence in $M \rtimes_{(\alpha , v)} G$. 
\end{lemma}
\begin{proof}
First, we show $(x_m)_{m \ge 1}$ is asymptotically invariant under the action of any element in $G$. Let $g,h \in S$. By the first cocycle relation, we have 
\[
\alpha_g(\alpha_h(x)) = v(g,h) \alpha_{gh}(x) v(g,h)^*
\]
for every $x \in M$. Hence,
\begin{align*}
\norm{\alpha_{gh}(x_m) - x_m}_2 &= \norm{v(g,h) \alpha_{gh}(x_m)v(g,h)^* - v(g,h) x_m v(g,h)^*}_2 \\
&= \norm{\alpha_g(\alpha_h(x_m)) - v(g,h) x_m v(g,h)^*}_2  \\
&\le \norm{\alpha_g(\alpha_h(x_m)) - \alpha_g(x_m)}_2 + \norm{\alpha_g(x_m) - x_m}_2 + \norm{x_m - v(g,h) x_m v(g,h)^*}_2 \\
&\le \norm{\alpha_g(\alpha_h(x_m) - x_m)}_2 + \norm{\alpha_g(x_m) - x_m}_2  + \norm{x_m - v(g,h) x_m v(g,h)^*}_2 \\
&\le \norm{ \alpha_h(x_m) - x_m}_2 + \norm{\alpha_g(x_m) - x_m}_2  + \norm{[x_m, v(g,h)]}_2.
\end{align*}
Note that the first two terms go to $0$ as $m \to \infty$ by virtue of asymptotic invariance under the generators, and the last norm goes to $0$ because $(x_m)_{m \ge 1}$ is an asymptotically commutative sequence in $M$. Now we show asymptotic invariance under inverses.

First, note that from the first cocycle relation, it is easy to deduce that $\alpha_e$ must be the identity map on $M$. Hence,

\begin{align*}
\norm{\alpha_{g^{-1}}(x_m) - x_m}_2 &= \norm{\alpha_g(\alpha_{g^{-1}}(x_m) - x_m)}_2 \\
&= \norm{\alpha_g(\alpha_{g^{-1}}(x_m)) - \alpha_g(x_m)}_2 \\
&= \norm{v(g,g^{-1})^* \alpha_g(\alpha_{g^{-1}}(x_m)) v(g,g^{-1}) - v(g,g^{-1})^* \alpha_g(x_m) v(g,g^{-1})}_2 \\
&= \norm{\alpha_{gg^{-1}}(x_m) - v(g,g^{-1})^* \alpha_g(x_m) v(g,g^{-1})}_2 \\
&= \norm{v(g,g^{-1}) x_m v(g,g^{-1})^* - \alpha_g(x_m)}_2 \\
&\le \norm{v(g,g^{-1}) x_m v(g,g^{-1})^* - x_m}_2 + \norm{x_m - \alpha_g(x_m)}_2 \\
&\le \norm{[v(g,g^{-1}),x_m]}_2 + \norm{\alpha_g(x_m) - x_m}_2 
\end{align*}
which converges to $0$ as $m \to \infty$ for the same reasons. This establishes asymptotic invariance under the inverse of generators. Inductively, we can extended to arbitrary finite products in $S$ which represents any element in $G$.

Finally, we show $(\pi(x_m))_{m \ge 1}$ asymptotically commutes with everything in the crossed product. It suffices to show $(\pi(x_n))_{m \ge 1}$ asymptotically commutes with an element of the form $\pi(x) \lambda_v(g)$, as the asymptotic commutativity can be extended to finite sums and then by density to infinite sums. Recall that $\norm{xy}_2 \le \norm{x}_{\infty} \cdot \norm{y}_2$. Then given $g \in G$, we have
\begin{align*}
\norm{\pi(x_m) \pi (x) \lambda_v(g) - \pi (x) \lambda_v(g) \pi (x_m)}_2 &= \norm {\pi (x_m x) \lambda_v(g) - \pi (x) \pi(\alpha_g(x_m)) \lambda_v(g)}_2 \\
&= \norm{(\pi (x_m x) - \pi(x\alpha_g(x_m))) \lambda_v(g)}_2 \\
&= \norm{\pi(x_m x - x \alpha_g(x_m))}_2 \\
&= \norm{x_m x - x \alpha_g(x_m)}_2 \\
&\le \norm{x_m x - x x_m}_2 + \norm{x x_m - x \alpha_g(x_m)}_2 \\
&= \norm{x_m x - x x_m}_2  + \norm{x(x_m - \alpha_g(x_m))}_2 \\
&\le \norm{x_m x - x x_m}_2  + \norm{x}_{\infty} \norm{\alpha_g(x_m) - x_m}_2 
\end{align*}
all of which converges to $0$ as $m \to \infty$. 
\end{proof}

If $(x_m)_{m \ge 1}$ is further assumed to be a central sequence, then Lemma \ref{lemma:CentSeqLift1} can be phrased alternatively in terms of ultrapowers. It says that the fixed point algebra $(M' \cap M^\omega )^G$ under the induced ultrapower cocycle action of $G$ naturally sits inside the central sequence algebra of $M \rtimes_{(\alpha , v)} G$. 

Now let us recall how a short exact sequence gives rise to a twisted crossed product decomposition of a group von Neumann algebra. Let 
\[
1 \longrightarrow N \xrightarrow{\iota} G \xrightarrow{\rho} Q \longrightarrow 1
\]
be a short exact sequence, where $N$ is a normal subgroup of $G$, $\iota$ the canonical embedding, and $\rho$ a surjective homomorphism such that $N = \image \iota = \ker \rho$. Since $\rho$ is surjective, we can find $s(q) \in G$ such that $\rho (s(q)) = q$, and for simplicity we always choose $s$ such that $s(e) = e$. Define $\alpha : Q \to \aut (L(N),\tau)$ by $\alpha_q(x) = \lambda_{G}(s(q)) x \lambda_{G}(s(q))^*$, where $\lambda_G$ is the usual left regular representation of $G$. This is well-defined because $N$ is a normal subgroup of $G$. We define $v : Q \times Q \to \mathcal{U}(L(N))$ by $v(q,r) = \lambda_{G}(s(q) s(r) s(qr)^{-1})$, which does in fact define a unitary in $L(N)$. 

Now, let $s_1, s_2 : Q \to G$ be two sections associated to $\rho : G \to Q$, both satisfying $s_1(e)=e$ and $s_2(e)=e$, and let $(\alpha_1,v_1)$ and $(\alpha_2,v_2)$ be the associated cocycle actions. Given $q \in Q$, we have $\rho(s_1(q) s_2(q)^{-1}) = e$. Hence, for every $q \in Q$, there exists $n_q \in \ker \rho = N$ such that $s_1(q) = n_q s_2(q)$. Hence, for $x \in L(N)$, 
\begin{align*}
\alpha_{1,q}(x) &= \lambda (s_1(q)) x \lambda (s_1(q))^* \\
&= \lambda (n_q s_2(q))  x \lambda (n_q s_2(q))^* \\
&= \lambda (n_q) \lambda (s_2(q)) x \lambda (s_2(q))^* \lambda (n_q)^* \\
&= i_{n_q} \alpha_{2,q}(x) 
\end{align*}
where $i : Q \to \text{Inn}(L(N))$ with $i(n_q) := i_{n_q}$ is the inner automorphism given by conjugating by $\lambda (n_q)$. Therefore, any two sections associated to a given homomorphism, their induced cocycle actions differ by an inner automorphism.

Now, suppose $(x_m)_{m \ge 1}$ is asymptotically commutative sequence in $L(N)$ that is asymptotically invariant under $(\alpha_2,v_2)$. Then for any $q \in Q$
\begin{align*}
\norm{\alpha_{1,q}(x_m) - x_m}_2 &= \norm{i_{n_q} \alpha_{2,q}(x_m) - x_m}_2 \\
&= \norm{\alpha_{2,q}(x) - i_{n_q^{-1}}(x_m)}_2 \\
&\le \norm{\alpha_{2,q}(x_m) - x_m}_2 + \norm{x_m - i_{n_q^{-1}}(x_m)}_2 \\
&= \norm{\alpha_{2,q}(x_m) - x_m}_2 + \norm{[x_m, \lambda (n_q)]}_2
\end{align*}
which converges to $0$ due to asymptotic invariance with respect to $\alpha_{2}$ and asymptotic commutativity of $(x_m)_{m \ge 1}$ in $L(N)$. Hence, if $(x_m)_{m \ge 1}$ is asymptotically invariant with respect to some cocycle action constructed from one section, it is asymptotically invariant with respect to any. Thus, any convenient lift $q \in Q$ suffices to check asymptotic invariance. 

Finally, we will use \cite[Theorem 3.1, Corollary 3.2]{Matsumoto1988}. Mastumoto only states Corollary 3.2 of his paper for a crossed product von Neumann algebra, but it also holds for a twisted crossed product von Neumann algebra. 

\begin{corollary}(\cite[Corollary 3.2]{Matsumoto1988})
    Let $M$ be a II$_1$ factor and $(\alpha, v)$ a free cocycle action of a countable discrete amenable group $G$ on $M$. If $M$ is a McDuff II$_1$ factor, then $M \rtimes_{(\alpha, v)}G$ is also a McDuff factor. 
\end{corollary}

\begin{proof}
    By the commutation relation $\lambda_v(g)\pi(x)\lambda_v(g)^*=\pi(\alpha_g(x))$, so we can view $G$ as a subgroup of $\mathcal N(M)$. Since $G$ and $M$ generate $M \rtimes_{(\alpha, v)}G$, the proof follows from \cite[Theorem 3.1]{Matsumoto1988}.
\end{proof}

\section{Shift-similar and Houghton-like Groups}\label{sec:ShiftSimandHoughtonLike}

In this section, we define shift-similar (Definition \ref{def:ShiftSim}) and Houghton-like groups (Definition \ref{def:HoughtonLikeGroups}) from a slightly different perspective than that given in \cite{mallery25}, but we prove both descriptions are equivalent. It turns out that many properties of Houghton-like groups are true for a much more general class subgroups of the permutation group acting on the disjoint union of $n$ copies $\mathbb{N}$. Thus, we start in this setting. 

\subsection{Basics on Permutation Groups} Given a set $X$, we write $\sym (X)$ for the group of permutations on $X$, $S_{\infty}(X)$ for the subgroup of finitely supported permutations on $X$, and $A_{\infty}$(X) for the subgroup of all finitely supported even permutations on $X$, which is an infinite simple group (when $|X| \ge 5$). Given a subset $Y \subseteq X$, we will naturally identify $\sym (Y)$ with the subgroup of permutations in $\sym (X)$ supported on $Y$ and acting as the identity outside of $Y$. When $X = \mathbb{N}$, for convenience we define $S_{\infty} := S_{\infty}(\mathbb{N})$ and $A_{\infty} := A_{\infty}(\mathbb{N})$.  

Our first lemma is that if $G$ is a group of permutations on $X$ containing $A_{\infty}(X)$, then it must have the infinite conjugacy class (ICC) property, meaning that the conjugacy class of any non-trivial element is infinite. We also prove $A_{\infty}(X)$ has trivial centralizer in $\sym (X)$. 

\begin{lemma}
Let $X$ be an infinite set, and let $G \le \sym (X)$ be a subgroup containing $A_{\infty}(X)$. Then every non-trivial element of $G$ has infinitely many $A_{\infty}(X)$-conjugates. Furthermore, the centralizer of $A_{\infty}(X)$ in $\sym(X)$ is trivial.
\end{lemma}
\begin{proof}
An element of a group has infinitely many $A_{\infty}(X)$-conjugates if and only if the centralizer of $g$ in $A_{\infty}(X)$ is infinite index. Suppose $g \in G$ has finite-index centralizer in $A_{\infty}(X)$. We will argue that $g=1$. Since the finitary alternating group $A_{\infty}(X)$ is an infinite simple group, the only finite-index subgroup of $A_{\infty}(X)$ is $A_{\infty}(X)$ itself. It follows that centralizer of $g$ in $A_{\infty}(X)$ is equal to all of $A_{\infty}(X)$, meaning that $g$ must commute with every element of $A_{\infty}(X)$. The fact that $g$ commutes with any $f \in A_{\infty}(X)$ means that $g$ and $f$ must stabilize each other's fixed point sets and, consequently, their supports. Since every subset of size three is the support of some element in $A_{\infty}(X)$, it follows $g$ must stabilize every subset of size three in $X$. Hence, if $x,y,z \in X$ with $x$ arbitrary and $y$ and $z$ distinct from each other and from $x$ and $g^{-1}(x)$, then $g$ must stabilize the set $\{g^{-1}(x),y,z\}$, meaning 
\[
\{g^{-1}(x),y,z\} = g(\{g^{-1}(x),y,z\}) = \{x,g(y),g(z)\}.
\]
Since $x \neq y,z$, this set equality forces $g^{-1}(x) = x$ and hence $g(x) = x$. Since $x \in X$ was arbitrary, $g$ must be the identity.

As for the trivial centralizer claim, let $g \in \sym (X)$ be a non-trivial element commuting with every element of $A_{\infty}$. Then there exist $x,y \in X$ distinct such that $g(x) = y$. Pick any $z \in X \setminus \{x,y, g(y)\}$. The $f = (x~y~z)$ defines an element in $A_{\infty}(X)$, so $g$ must commute with $f$. That is, $g = fgf^{-1}$. Applying this equation to $z$, we get 

$$g(y) = fgf^{-1}(y) = fg(x) = f(y) = z,$$
which is a contradiction. 
\end{proof}

For $n \in \mathbb{N}$, define $X_n := \{1, \dots , n \} \times \mathbb{N}$ as the disjoint union of $n$-copies of $\mathbb{N}$. Given $k= 1, \dots , n$, we define the $k$-th ray as 
\[
R_k := \{k\} \times \mathbb{N},
\]
and given $j \in \mathbb{N}$, define the $j$-th tail $R_k(j)$ as 
\[
R_k(j) := \{k \} \times \{j,j+1,j+2, \dots \}.
\]
We will be studying subgroups of $\sym (X_n)$, particularly Houghton-like groups. But for our purposes, the group $\sym(X_n)$ is too large and does not have the properties we want. Instead, we consider subgroups of $\sym_{ends}(X_n)$, the group of permutations $g$ of $X_n$ that eventually stabilize the ends or rays, meaning eventually along each ray, the bijection does not move points across rays. More precisely, a permutation $g$ of $X_n$ belongs to $\sym_{ends}(X_n)$ if and only if 
\[
|g(R_k) \triangle R_k| < \infty
\]
for every $k=1, \dots , n$. Clearly $\sym_{ends}(X_n)$ contains the identity, is closed under inverses, and is closed under products by the triangle inequality.  First, we prove a lemma about $\sym_{ends}(X_n)$ which makes this idea of stabilizing the ends of $X_n$ more explicit.

\begin{lemma}\label{lemma:threshold}
A permutation $g$ of $X_n$ belongs to $\sym_{ends}(X_n)$ if and only if there exists $i_g \in \mathbb{N}$ such that
\[
g(R_k(i_g)) \subseteq R_k
\]
for every $k=1, \dots , n$
\end{lemma}
\begin{proof}
Let $g \in \sym_{ends}(X_n)$. Then for every $k=1, \dots , n$,
\[
|g(R_k) \triangle R_k|
\]
is finite. In particular, 
\[
g(R_k) \setminus R_k
\]
is finite, meaning that only finitely many points on $R_k$ are mapped off of $R_k$ by $g$. Therefore, we can pick $i_k \in \mathbb{N}$ large enough so that every point on $R_k$ whose second coordinate is at least $i_k$ is mapped onto $R_k$. Defining $i_g$ as $i_g := \max \{i_1, \dots , i_n\}$, it is clear we have 
\[
g(R_k(i_g)) \subseteq R_k
\]
for each $k=1, \dots , n$.

Now, let $g$ be a permutation of $X_n$ for which there exists $i_g \in \mathbb{N}$ such that 
\[
g(R_k(i_g)) \subseteq R_k
\]
for every $k=1, \dots , n$. First, suppose $x \in g(R_k) \setminus R_k$. Then $x = g(k,j)$ for some $j \in \mathbb{N}$. If $j \ge i_g$, then $(k,j) \in R_k(i_g)$ and therefore 
\[
x = g(k,j) \in g(R_k(i_g)) \subseteq R_k
\]
which is a contradiction. Hence, $j < i_g$ which means there are at most finitely many choices for $j$ and hence for $x$. This means $g(R_k) \setminus R_k$ is finite. 

Now, suppose $x \in g^{-1}(R_k) \setminus R_k$. Then $x = g^{-1}(k,j)$ or $g(x) = (k,j)$ for some $j \in \mathbb{N}$. Let $x = (i, \ell)$ for some $i=1, \dots , n$ and $\ell \in \mathbb{N}$. Because $x \notin R_k$, $i \neq k$. If $\ell \geq i_g$, then because 
\[
g(R_i(i_g)) \subseteq R_i
\]
the first coordinate of $g(i, \ell)$ is $i$, which makes the equation $g(i, \ell ) = (k, j)$ impossible. Hence, $\ell < i_g$. This shows that there are at most finitely many choices for $\ell$ and $i$, so there at most finitely many choices for $x$. Hence, it follows $g^{-1}(R_k) \setminus R_k$ is finite and therefore $R_k \setminus g(R_k)$ is finite. Because
\[
g(R_k) \triangle R_k = (g(R_k) \setminus R_k) \cup (R_k \setminus g(R_k)) 
\]
it follows the symmetric difference is finite. Hence, $g \in \sym_{ends}(X_n)$. 
\end{proof}

Now we define a homomorphism on $\sym_{ends}(X_n)$ that tracks just the eventual behavior of permutations in $\sym_{ends}(X_n)$ along each ray. Before we can define this map, we recall the notion of an almost permutation of $\mathbb{N}$ and the group one can form from this concept.

We call $f$ an \emph{almost permutation} of $\mathbb{N}$ if there exist cofinite sets $A,B \subseteq \mathbb{N}$ such that $f$ is a bijection from $A$ to $B$. We say two almost permutations $f$ and $g$ of $\mathbb{N}$ are equivalent if they agree on the intersection of their domains (equivalently, they agree on a subset of their domains which is cofinite in $\mathbb{N}$). Given an almost permutation $f$, we let $[f]$ denote its equivalence class and $\sym^* (\mathbb{N})$  the set of all these equivalence classes. We multiply two germs of almost permutations $[f], [g] \in \sym^* (\mathbb{N})$ via the formula $[f] \cdot [g] := [f \circ g]$, where the composition is defined from $\dom(g) \cap g^{-1}(\dom(f))$ to $f(\image (g) \cap \dom(f))$. This is well-defined because the intersection of two cofinite sets is cofinite, and the image and preimage of cofinite sets under a bijection are cofinite.  

There is a natural homomorphism from $\sym^*(\mathbb{N})$ onto the integers, analogous to the index of a Fredholm operator. Given $[f] \in \sym^*(\mathbb{N})$, we define  $ind : \sym^*(\mathbb{N}) \to \mathbb{Z}$ by
\[
\ind([f]) = |\mathbb{N} \setminus B| - |\mathbb{N} \setminus A|,
\]
where $f : A \to B$ is an almost permutation. Under this definition, we have  $\ind([\sigma ]) = 1$. Hence, this is a surjective homomorphism and therefore splits. Its kernel consists precisely of those classes $[f]$ such that $f$ can be completed (in a non-unique way) to a full permutation of $\mathbb{N}$.

Given $g \in \sym_{ends}(X_n)$ and $k=1, \dots , n$, define 
\[
D_{g,k} = \{j \in \mathbb{N} : g(k,j) \in R_k\}
\]
and 
\[
C_{g,k} = \{j \in \mathbb{N} : (k,j) \in g(R_k)\}.
\]
Observe that these are cofinite sets in $\mathbb{N}$. Indeed, if $j \notin D_{g,k}$, then $g(k,j) \in g(R_k) \setminus R_k \subseteq g(R_k) \triangle R_k$ and therefore $(k,j) \in R_k \triangle g^{-1}(R_k)$. Note that $R_k \triangle g^{-1}(R_k)$ is finite because $g$ is a bijection, so there are at most finitely many choices for $j$. Hence $D_{g,k}$ is cofinite. By a similar argument, $C_{g,k}$ is also cofinite.

Now, given $j \in D_{g,k}$, define $f_{g,k}(j)$ as the unique natural number such that $g(k,j) = (k,f_{g,k}(j))$. This exists and is unique because $j \in D_{g,k}$ and $g$ is a function, so $f_{g,k}$ is well-defined. Moreover, by construction $f_{g,k}(j) \in C_{g,k}$. If $f_{g,k}(j_1) = f_{g,k}(j_2)=n$, then $g(k,j_1) = (k,n) = g(k,j_2)$, and injectivity of $g$ implies $j_1 = j_2$, so $f_{g,k}$ is injective. If $n \in C_{g,k}$, then $(k,n) \in g(R_k)$, so $g(k,m) = (k,n)$ for some $m \in \mathbb{N}$. Since $g(k,m) \in R_k$, we have $m \in D_{g,k}$, and by definition we have $f_{g,k}(m)=n$. Hence, $f_{g,k}$ is surjective. This proves $f_{g,k}$ is a well-defined bijection from $D_{g,k}$ to $C_{g,k}$, and therefore an almost permutation of $\mathbb{N}$.

Now, define $\phi_k : \sym_{ends}(X_n) \to \sym^*(\mathbb{N})$ by 
\[
\phi_k(g) = [f_{g,k}].
\]
We claim this is a group homomorphism. First, recall that if $f$ is an injective map between cofinite subsets of $\mathbb{N}$, then $\lim_{x \to \infty} f(x) = \infty$. Now, to show the homomorphism property, let $g,h \in \sym_{ends}(X_n)$. It suffices to show $f_{gh,k}$ and $f_{g,k} \circ f_{h,k}$ agree on a cofinite subset of $\mathbb{N}$. Because all three maps are injections between cofinite subsets, there exists $j_0 \in \mathbb{N}$ large enough so that whenever $j \ge j_0$, we have $j \in D_{gh,k}$, $j \in D_{h,k}$, and $f_{h,k}(j) \in D_{g,k}$. Hence, $f_{gh,k}(j)$ and $f_{g,k}(f_{h,k}(j))$ are well-defined. By definition (and since $j \ge j_0$), we have $h(k,j) = (k,f_{h,k}(j))$ and $g(k,f_{h,k}(j)) = (k,f_{g,k}(f_{h,k}(j)))$, and therefore
\[
gh(k,j) = g(k,f_{h,k}(j)) = (k,f_{g,k}(f_{h,k}(j))).
\]
But $j \in D_{gh,k}$, so by definition, we have 
\[
gh(k,j) = (k,f_{gh,k}(j)).
\]
Therefore, comparing second coordiantes, we have $f_{gh,k}(j) = f_{g,k}(f_{h,k}(j))$ for every $j \ge j_0$. Since $f_{gh,k}$ and $f_{g,k} \circ f_{h,k}$ agree on a cofinite set, they represent the same equivalence class of almost permutations. Hence,
\[
[f_{gh,k}] = [f_{g,k} \circ f_{h,k}] = [f_{g,k}] \cdot [f_{h,k}]
\]
which is equivalent to 
\[
\phi_k(gh) = \phi_k(g) \phi_k(h).
\]

Now, a permutation $g \in \sym_{ends}(X_n)$ belongs to the kernel if and only if $\phi_k(g) = [\id_{\mathbb{N}}]$ if and only if $f_{g,k}(j) = \id_{\mathbb{N}}(j) = j$ for all sufficiently large $j$. By definition, this means
\[
g(k,j) = (k,f_{g,k}(j)) = (k,j)
\]
for all sufficiently large $j$. Hence, $g$ belongs to the kernel if and only if $g$ eventually pointwise stabilizes the $k$-th ray $R_k$. Because of this, we denote the kernel by $\estab(R_k)$. 

We summarize the preceding discussion as the following lemma.

\begin{lemma}\label{lemma:PhiHom}
Let $k=1, \dots , n$. The map $\phi_k : \sym_{ends}(X_n) \to \sym^* (\mathbb{N})$ given by 
\[
\phi_k(g) = [f_{g,k}]
\]
is a well-defined group homomorphism, where $f_{g,k}$ is the almost permutation of $\mathbb{N}$ defined between the cofinite sets 
\[
D_{g,k} = \{j \in \mathbb{N} : g(k,j) \in R_k\},
\]
the domain, and
\[
C_{g,k} = \{j \in \mathbb{N} : (k,j) \in g(R_k)\},
\]
the codomain. Given $j \in D_{g,k}$, $f_{g,k}(j)$ is the unique value such that 
\[
g(k,j) = (k,f_{g,k}(j)).
\]
The kernel is $\estab(R_k)$, the group of all permutations in $\sym_{ends}(X_n)$ that eventually pointwise stabilize the $k$-th ray. 
\end{lemma}

The maps $\phi_1, \dots, \phi_n$ can be combined into the direct product map $\Phi : \sym_{ends}(X_n) \to \sym^*(\mathbb{N})^n$ defined as 
\[
\Phi (g) = (\phi_1(g), \dots , \phi_n(g))
\]
with kernel
\[
\ker \Phi = \bigcap_{k=1}^{n} \estab(R_k) = S_{\infty}(X_n).
\]
Given $g \in \sym_{ends}(X_n)$, define sets $A$ and $B$ as
\[
A = \bigsqcup_{k=1}^{n} (\{k\} \times D_{g,k}), \quad B = \bigsqcup_{k=1}^{n} ( \{k\} \times C_{g,k}).
\]
By construction, $g(A) = B$ so $g(X_n \setminus A) = X_n \setminus B$. Therefore, $|X_n \setminus A| = |X_n \setminus B|$ which translates to
\[
\sum_{k=1}^{n} |\mathbb{N} \setminus D_{g,k}| = \sum_{k=1}^{n} |\mathbb{N} \setminus C_{g,k}|
\]
and finally
\[
\sum_{k=1}^{n} \ind (\phi_k(g)) = 0
\]
Hence,
\[
\image \Phi = \{(\alpha_1, \dots \alpha_n) \in \sym^*(\mathbb{N})^n : \sum_{k=1}^{n} \ind (\alpha_k) = 0 \}.
\]

\subsection{Shift-similar and Houghton-like Groups} Now we restrict our attention to certain subgroups of $\sym_{ends}(X_n)$---namely, Houghton-like groups. Recall that the classical Houghton group $H_n$ acts on $X_n$ by permutations which are eventually translations along each ray (either forwards or backwards). More precisely, a permutation $g \in \sym (X_n)$ belongs to $H_n$ if and only if there exists a threshold $i_g \in \mathbb{N}$ and there exist integers $d_k \in \mathbb{Z}$ with $\sum_{k=1}^{n} d_k = 0$ such that 
\[
g(k,j) = (k,j+d_k)
\]
for all $k=1, \dots , n$ and $j \ge i_g$. Lemma \ref{lemma:threshold} immediately implies $H_n$ is a subgroup of $\sym_{ends}(X_n)$. 

Before we introduce Houghton-like groups, we define shift-similar groups. Intuitively, a group $G \le \sym (\mathbb{N})$ is \emph{shift-similar} if it is closed under the following operation. Take any permutation $g \in G$ and think of $g$ as a diagram of arrows: two horizontal, parallel copies of $\mathbb{N}$ with an arrow connecting $j \in \mathbb{N}$ in the bottom copy of $\mathbb{N}$ to $g(j) \in \mathbb{N}$ in the top copy. Given $j \in \mathbb{N}$, we delete $j$ in the bottom copy of $\mathbb{N}$, $g(j)$ in the top copy of $\mathbb{N}$, delete the arrow connecting $j$ and $g(j)$, and then relabel the top and bottom to be $\mathbb{N}$ again. If the resulting permutation is in $G$ for any $g \in G$ and any $j \in \mathbb{N}$, we say $G$ is shift-similar. We can make this more precise as follows.

\looseness=-1 Let $j \in \mathbb{N}$, and define $\mathbb{N}_j := \mathbb{N} \setminus \{j\}$. Consider $s_j : \mathbb{N} \to \mathbb{N}_j$ defined by sending $i$ to $i$ if $1 \le i < j$ and $i$ to $i+1$ if $i \ge j$. Define the $j$-th deletion-shift map $\psi_j : \sym(\mathbb{N}) \to \sym (\mathbb{N})$ by
\[
\psi_j(g) := s_{g(j)}^{-1} \circ g \big|_{\mathbb{N}_j} \circ s_j.
\]
Since $g$ is injective, we know that $\lim_{x \to \infty} g(x) = \infty$, so for all sufficiently large $x$, we have $g(x+1) > g(j)$ and therefore $\psi_j(g)(x) = g(x+1)-1$. Although these maps are not group homomorphisms, they satisfy something like a cocycle condition. Indeed, given $g,h \in \sym (\mathbb{N})$, we have 
\[
\psi_j(g \circ h) = s_{gh(j)}^{-1} \circ (g \circ h) |_{\mathbb{N}_j} \circ s_j = s_{gh(j)}^{-1} \circ g|_{\mathbb{N}_{h(j)}} \circ s_{h(j)} \circ s_{h(j)}^{-1} \circ h|_{\mathbb{N}_j} \circ s_j = \psi_{h(j)}(g) \circ \psi_j(h).
\]

The deletion-shift map $\psi_j$ implements the operation described above that intuitively characterizes shift-similarity.

\begin{definition}\label{def:ShiftSim}
A group $G \le \sym (\mathbb{N})$ is \emph{shift-similar} if $\psi_j(G) \subseteq G$ for every $j \in \mathbb{N}$. Given a group $G$ not necessarily acting on $\mathbb{N}$ a priori, we say $G$ is shift-similar if it admits an embedding into $\sym (\mathbb{N})$ such that its image under the embedding is shift-similar in the aforementioned sense. 
\end{definition}

We summarize the common examples of shift-similar groups presented in \cite{mallery25}. 

\begin{examples}\label{ex:ShiftSimEx}

First, the simplest infinite shift-similar group is $S_{\infty}$. Indeed, thinking of $S_n$ as supported on $\{1, \dots , n\} \subset \mathbb{N}$, $S_n$ is shift-similar (in fact, these are the only finite shift-similar groups), and since shift-similarity is closed under directed unions, it follows $S_{\infty} = \bigcup_{n \ge 1} S_n$ is shift-similar (\cite[Lemma 3.6, Example 3.7]{mallery25}). Moreover, an infinite shift-similar group must contain $S_{\infty}$ as a subgroup (see \cite[Theorem 3.12]{mallery25}). 

There is a natural isomorphism from $\sym (X_n) \to \sym (\mathbb{N})$ under which the image of $H_n$ is shift-similar (see \cite[Corollary 3.10]{mallery25}), so $H_n$ is shift-similar. This contrasts with the Higman-Thompson groups which cannot be self-similar because they are not residually finite.

Further examples include eventually periodic permutations. Given $n \in \mathbb{N}$, let $E_n \le \sym (\mathbb{N})$ denote the group of all eventually $n$-periodic permutations, that is, all $g$ such that 
\[
g(j+n) = g(j)+n
\]
for all sufficiently large $j \in \mathbb{N}$.

\looseness=-1 A permutation $g \in \sym (\mathbb{N})$ has \emph{half-finite support} if $g(j) = j$ for all sufficiently large even $j$ or odd $j$. If $\Gamma$ is a group of half-finite support permutations, then $G := \langle \Gamma, E_2 \rangle$ is shift-similar. 

Now, if $\Gamma$ is any countable group and $\phi : \Gamma \to \mathbb{N}$ any bijection, we can identify the disjoint union $\Gamma \coprod \Gamma$ with $\mathbb{N}$ by mapping the first copy of $\Gamma$ to $2 \mathbb{N}$ via $g \mapsto 2 \phi (g)$, and sending the second copy of $\Gamma$ to $2 \mathbb{N} -1$ via $g \mapsto 2 \phi (g)-1$. Define $G$ to be the group of permutations on $\mathbb{N}$ generated by the copy of $\Gamma$ acting on $2 \mathbb{N}$ by left multiplication and trivially on $2 \mathbb{N} -1$, and the permutations $(1~2)$, $\alpha = (\dots 6~4~2~1~3~5~7 \dots )$, and $\lambda = (1~2)(3~4)(5~6) \dots $. This group is shift-similar, and if $\Gamma$ is finitely generated, so is $G$. 

Finally, as alluded to in the introduction, we can produce examples of shift-similar groups from shift-similar closures. Given a group $G \le \sym (\mathbb{N})$, the \emph{shift-similar closure}, denoted $\overline{G}$, is the smallest group containing $G$ that is closed under the deletion-shift maps $\psi_j$. This will produce a wealth of examples, and will be our primary source.
\end{examples}
\hfill$\blacktriangleleft$

We will almost exclusively be interested in infinite shift-similar groups, and non-amenable ones at that. This is because the Houghton-like group $H_n(G)$ (defined below) is (elementary) amenable if and only if $G$ is (elementary) amenable (see \cite[Theorem 5.9]{mallery25}). In this case, $L(H_n(G))$ is the hyperfinite II$_1$ factor and therefore McDuff, which makes our result trivial.  

Using this cocycle type condition, it turns out that we do not need all the deletion-shift maps in order to determine or characterize when as an infinite group is shift-similar as the next two lemmas shows.

\begin{lemma}
For $g \in \sym (\mathbb{N})$, the following are equivalent:
\begin{enumerate}
    \item $g \in S_{\infty}$;
    \item $\psi_j(g) \in S_{\infty}$ for every 
    $j \in \mathbb{N}$;
    \item $\psi_j(g) \in S_{\infty}$ for some $j \in \mathbb{N}$. 
\end{enumerate}
\end{lemma}
\begin{proof}
$(1)$ implies $(2)$ is true by virtue of $S_{\infty}$ being shift-similar; and clearly statement $(2)$ implies $(3)$. Hence, all that remains is to prove $(3)$ implies $(1)$. 

Suppose $\psi_j(g) \in S_{\infty}$ for some $j \in \mathbb{N}$. Then there exists $k \in \mathbb{N}$ such that $\psi_j(g)(i) = i$ for every $i \ge k$. Let $k' := \max \{k+1,j+1,g(j)+1\}$. If $i \ge k'$, then $i-1 \ge k$, $i-1 \ge g(j)$, and finally $i-1 \ge j$. Together these three inequalities jointly entail that 
\[
\psi_j(g)(i-1) = i-1
\]
or
\[
g \big|_{\mathbb{N}_j} (s_j(i-1)) = s_{g(j)}(i-1)
\]
or
\[
g(i) = i
\]
for every $i \ge k'$, thereby proving $g \in S_{\infty}$. 
\end{proof}

Now we prove that, for infinite groups, the deletion-shift map $\psi_1$ contains all the information needed to characterize shift-similarity. 

\begin{lemma}
Let $G \le \sym (\mathbb{N})$ an infinite group. Then $G$ is shift-similar if and only if $G$ contains $S_{\infty}$ and is closed under $\psi_1$. 
\end{lemma}
\begin{proof}
The forward direction is trivial. Suppose $G$ contains $S_{\infty}$ and is closed under $\psi_1$. Let $g \in G$ and $j \in \mathbb{N}$ be arbitrary. Note that $(1~j) \in S_{\infty} \le G$. Hence, the cocycle type condition implies
\[
\psi_j((1~j) \circ g) = \psi_{g(1)}(1 ~ j) \circ \psi_j(g).
\]
Solving for $\psi_j(g)$
\[
\psi_j(g) = \psi_1((1~j) \circ g) \circ \psi_{g(1)}(1~j)^{-1} \in G S_{\infty} = G.
\]
\end{proof}

Thus, we will be exclusively interested in the map $\psi_1$ and hence the shift map $s_1 : \mathbb{N} \to \mathbb{N} \setminus \{1\}$, which we  henceforth denote as $\sigma$.

We present the Houghton-like groups as defined in \cite{mallery25}, and then derive a description closer in spirit to that of the classical Houghton groups. First, given a finite set $M \subset X_n$ and $1 \le k \le n$, we define the $k$-th quasi-ray $QR_k(M)$ formed from $M$ as
\[
QR_k(M) = R_k \setminus M.
\]
Enumerating the points in $QR_k(M)$ as
\[
QR_k(M) = \{(k,i_1),(k,i_2), \dots \}
\]
gives us a canonical bijection $\beta_{k,M} : \mathbb{N} \to QR_k(M)$ defined by 
\[
\beta_{k,M}(j) = (k,i_j).
\]
\looseness=-1 that enumerates in $QR_k(M)$ in increasing order. If $(k,j) \in QR_k(M)$ has $p$ points in $QR_k(M)$ strictly preceding it, then its index in the numeration is $p+1$; that is, $\beta_{k,M}(p+1) = (k,j)$.

Finally, we now define Houghton-like groups.

\begin{definition}\label{def:HoughtonLikeGroups}
Let $G \le \sym (\mathbb{N})$ be a shift-similar group, and $n \in \mathbb{N}$. The Houghton-like group $H_n(G)$ consists of all permutations $h$ of $X_n$ constructed as follows:
\begin{enumerate}
\item Take two finite subsets $M_+, M_- \subset X_n$ with the same cardinality. This gives us two partitions of $X_n$ into the finite sets $M_+$ and $M_-$ with the corresponding quasi-rays \newline $QR_1(M_{+}), \dots , QR_n(M_{+})$ and $QR_1(M_{-}), \dots , QR_n(M_{-})$
\item Pick a bijection $\tau$ from $M_+$ to $M_-$, and pick shift-similar group elements $h_1, \dots , h_n \in G$. Then $h$ maps $X_n$ bijectively to itself by sending $M_+$ to $M_-$ via $\tau$, and sending $QR_k(M_+)$ to $QR_k(M_-)$ via $\beta_{k,M_-} \circ h_k \circ \beta_{k,M_+}^{-1}$.
\end{enumerate}
\end{definition}

The idea of a Houghton-like group is to decompose the translation on each ray into a backward shift and a forward shift that balance out globally, and then intertwine these shifts with a shift-similar group element for each ray. More precisely, we will show $g \in \sym (X_n)$ belongs to $H_n(G)$ if and only if there exists $i_g \in \mathbb{N}$, positive integers $s_k,r_k$ for $k=1, \dots , n$ with $\sum_{k=1}^{n} s_k = \sum_{k=1}^{n} r_k$, and shift-similar group elements $g_1, \dots , g_n \in G$ such that 
\[
g(k,j) = (k,g_k(j-s_k) + r_k)
\]
for all $k=1, \dots , n$ and $j \ge i_g$. It is clear that $H_n(G)$ contains the classical Houghton group by taking $g_k=e$ and defining $d_k =  r_k - s_k$, and $H_n(G)$ reduces to $H_n$ when $G=\{e\}$.

Before we prove these two ways of describing a Houghton-like group are equivalent, we first need to prove the following lemma.

\begin{lemma}
Let $k=1, \dots , n$ and $M \subset X_n$ a finite subset. Then there exist  $i_k \in \mathbb{N}$ and $s_k \in \mathbb{N}_0$ such that $\beta_{k,M}(i) = (k,i+s_k)$ for all $i \ge i_k$.
\end{lemma}
\begin{proof}
Because $M$ is finite, we can form the number $s_k := |R_k \cap M|$ for each $k=1, \dots , n$. If $R_k \cap M$ is not empty, define $t_k := \max \{j \in \mathbb{N} : (k,j) \in M \}$, the largest second coordinate of points in $R_k \cap M$. Otherwise, define $t_k = 0$. 

Now, take any $j > t_k$. Then clearly $(k,j) \notin R_k \cap M$ and consequently $(k,j) \in QR_k(M)$. The number of points on $QR_k(M)$ strictly preceding $(k,j)$ equals the number of all points on $R_k$ strictly preceding $(k,j)$ minus those that belong to $M$. The number of points on $R_k$ strictly preceding $(k,j)$ is clearly $j-1$. Among these points, exactly $s_k$ of them lie in $M$ because all points of $M \cap R_k$ have second coordinate at most $t_k$, and $t_k < j$. Thus, the number of points in $QR_k(M)$ strictly preceding $(k,j)$ is 
\[
(j-1)-s_k.
\]
Because $\beta_{k,M}$ enumerates $QR_k(M)$ in increasing order, the index of $(k,j)$ in this enumeration is
\[
((j-1)-s_k) +1 = j-s_k
\]
and hence 
\[
\beta_{k,M}(j-s_k) = (k,j)
\]
for every $j > t_k$. Let $i = j-s_k$. Since $j > t_k$ is equivalent to $i > t_k - s_k$, we obtain 
\[
\beta_{k,M}(i) = (k,i+s_k).
\]
Define $i_k := \max \{1, t_k-s_k+1\}$. Then for every $i \ge i_k$, we have $i > t_k  - s_k$ and consequently
\[
\beta_{k,M}(i) = (k,i+s_k),
\]
which completes the proof. 
\end{proof}

Now we can prove the equivalence of these two descriptions of Houghton-like groups.

\begin{theorem}\label{thm:AlternateH_n(G)}
Let $G \le \sym (\mathbb{N})$ be a shift-similar group. Then $g \in H_n(G)$ if and only if there exist non-negative integers $s_1, \dots , s_n; r_1, \dots r_n$ with $\displaystyle \sum_{k=1}^{n} s_k = \sum_{k=1}^{n} r_k$, a natural number $i_g \in \mathbb{N}$, and shift-similar group elements $g_1, \dots , g_n \in G$ such that 

$$g(k,j) = (k, g_k(j-s_k) + r_k)$$
for all $k=1, \dots , n$ and $j \ge i_g$.
\end{theorem}
\begin{proof}
First, let $g \in H_n(G)$. As in the proof of the above lemma, we define the translation parameters $s_k$ and $r_k$ as $s_k = |R_k \cap M_+|$ and $r_k = |R_k \cap M_-|$. The sum balance condition $\sum_{k=1}^{n} s_k = \sum_{k=1}^{n} r_k$ holds because $|M_+| = |M_-|$. By the above lemma, there exists non-negative numbers $i_k^+, i_k^-$ such that
\[
\beta_{k,M_+}(j) = (k, j + s_k) \qquad \text{ for all } j \ge i_k^+
\]
and
\[
\beta_{k,M_-}(j) = (k,j + r_k) \qquad \text{ for all } j \ge i_k^-.
\]
For all sufficiently large $j$, the map $j \mapsto g_k(j-s_k)$  is an injection from a tail of $\mathbb{N}$ to $\mathbb{N}$. This implies $\displaystyle \lim_{j \to \infty} g_k(j-s_k) = \infty$ and hence there exists a natural number $N_k \in \mathbb{N}$ such that 
\[
g_k(j-s_k) > i_{k}^{-}
\]
for every $j \ge N_k$. Now, define 
\[
i_g := \max_{1 \le k \le n} \{i_k^+ + s_k, i_k^{-}, N_k\},
\]
and let $j \ge i_g$. Then $j \ge i_k^+ + s_k$ or equivalently $j  - s_k \ge i_k^+$. Replacing $j$ by $j-s_k$ in our formula for $\beta_{k,M+}$, we get 
\[
\beta_{k,M_+}(j-s_k) = (k, j - s_k + s_k) = (k,j)
\]
and hence 
\[\beta_{k,M_+}^{-1}(k,j) = j-s_k.
\]
Hence, 
\[
g(k,j) = (\beta_{k,M_-} \circ g_k \circ \beta_{k,M_+}^{-1})(k,j) = \beta_{k,M_-}(g_k(j-s_k)) = (k,g_k(j-s_k)+r_k)
\]
for all $j \ge i_g$ and all $k=1, \dots , n$. 

Now, suppose have a permutation $g$ as described in the statement of this theorem. Define 
\[
M_+ = \bigcup_{k=1}^{n} \{(k,j) :  1 \le j \le s_k \} \quad \text{ and } \quad M_- := \bigcup_{k=1}^{n} \{(k,j) : 1 \le j \le r_k \}.
\]
\looseness=-1 The balance condition $\sum_{k=1}^{n} s_k = \sum_{k=1}^{n} r_k$ implies $|M_+| = |M_-|$. Let $\tau$ be any bijection from $M_+$ to $M_-$, and define a bijection $h$ to map $M_+$ to $M_-$ via $\tau$ and the $k$-th quasi-ray $QR_k(M_+)$ to $QR_k(M_-)$ via $\beta_{k,M_-} \circ g_k \circ \beta_{k,M_+}^{-1}$. By definition, $h$ is an element of $H_n(G)$. Pick $J \in \mathbb{N}$ large enough so that for all $j \ge J$, we have the following conditions simultaneously hold:
\begin{enumerate}
\item $j \ge i_g$ so that the tail formula for $g$ holds
\item $\displaystyle j > \max \{s_k : k=1, \dots , n\}$ so that $(k,j) \notin M_+$
\item $j - s_k \ge i_k^+$ so that $\beta_{k,M_+}^{-1}(k,j) = j - s_k$
\item $g_k(j-s_k) \ge i_k^-$ so that $\beta_{k,M_-}(g_k(j-s_k)) = (k,g_k(j-s_k)+r_k)$
\end{enumerate}
Hence, if $j \ge J$, putting the above conditions together allows us to conclude 
\[
h(k,j) = \beta_{k,M_-}(g_k(\beta_{k,M_+}^{-1}(k,j))) = (k,g_k(j-s_k) + r_k) = g(k,j).
\]
\looseness =-1 This shows $g$ and $h$ agree on every quasi-ray except possibly the set finite set $M_+$, meaning that $h^{-1}g \in S_{\infty}(X_n)$. Since $ S_{\infty}(X_n) \le H_n(G)$, closure tells us that $h \cdot h^{-1}g = g$ belongs to $H_n(G)$. 
\end{proof}

Using this alternative description of Houghton-like groups, it immediately follows from Lemma \ref{lemma:threshold} that they are subgroups of $\sym_{ends}(X_n)$. Hence, the homomorphism $\phi_k : \sym_{ends}(X_n) \to \sym^* (\mathbb{N})$ from Lemma \ref{lemma:PhiHom} restricts to $H_n(G)$. We will compute the image and kernel. For convenience, we use $\phi_k$ to denote this restriction instead of $\phi_k |_{H_n(G)}$. First, we consider the case when $n \ge 2$.

First, we determine the image. Given $h \in H_n(G)$, by Theorem \ref{thm:AlternateH_n(G)}, there exist $i_h \in \mathbb{N}$, there exist $s_1, \dots , s_n; r_1, \dots , r_n \ge 0$, and there exist shift-similar group elements $g_1, \dots , g_n \in G$ such that 
\[
h(k,j) = (k,g_k(j-s_k)+r_k)
\]
for all $k=1, \dots , n$ and $j \ge i_h$. Hence, if $j \ge i_h$, we have
\[
h(k,j) = (k,g_k(j-s_k)+r_k) = (k,\sigma^{r_k} g_k \sigma^{-s_k}(j))
\]
and therefore
\[
\phi(h) [\sigma^{r_k} g_k \sigma^{-s_k}] = t^{r_k} [g_k] t^{-s_k},
\]
where $t := [\sigma ]$. Note that $t$ has infinite order because for every $n \in \mathbb{N}$, we have $\sigma^n(j) = j+n$ for all sufficiently large $j$, meaning $\sigma^n$ does not equal $id_{\mathbb{N}}$ on any cofinite subset of $\mathbb{N}$. Also, given $g \in G$, for all sufficiently large $j \in \mathbb{N}$, we have $g(j+1) > g(1)$ and therefore
\[
\psi_1(g)(j) = s_{g(1)}^{-1}g s_1(j) = s_{g(1)}^{-1}g(j+1) = g(j+1)-1 = \sigma^{-1}g \sigma (j)
\]
and hence
\[
t^{-1} [g] t = [\psi_1(g)].
\]
Thus, shift-similarity is equivalent to the condition that $t^{-1}[G]t \subseteq [G]$.

The above calculation shows that the image of $\phi_k$ lands in the group generated by $t$ and $[G] := \{[g] : g \in G \}$, the image of $G$ under the natural embedding $\sym (\mathbb{N}) \to \sym^*(\mathbb{N})$.  In fact, $\image \phi_k$ equals this group for $n \ge 2$. First, for any $\ell \in \{1, \dots , n\} \setminus \{k\}$, define $\tau_{k,\ell}$ by $\tau_{k, \ell}(\ell,j) = (\ell , j-1)$ for $j \ge 2$, $f_{k,\ell}(\ell, 1) = (k,1)$, $f_{k,\ell}(k,j) = (k,j+1)$, and let it fix all other points. Clearly $\tau_{k, \ell} \in H_n \le H_n(G)$, and furthermore $\phi_k(\tau_{k,\ell}) = [\sigma ] = t$ because on the $k$-th ray $f_{k, \ell}$ acts by $\sigma$. Now, given $g \in G$, define $h \in H_n(G)$ so that it acts on $R_k$ via $g$ and identity on all other rays. That is, $h_k = g$, $h_i =e$ for all $i \neq k$, and $s_i = r_i = 0$ for all $i=1, \dots , n$. Thus, $\phi_k(h) = [g]$. Hence, both types of generators lie in $\image \phi_k$, and therefore $\image \phi_k$ equals the group generated by $t$ and $[G]$. 

Therefore, the image of $\phi_k$ is an ascending HNN extension of $[G]$ by the stable letter $t = [\sigma]$. We can see more explicitly by looking at the group of germs $\mathcal{G}(G) := G/S_{\infty}$ with the monomorphism $\psi_{\infty} : \mathcal{G}(G) \to \mathcal{G}(G)$ defined by $\psi_{\infty}(gS_{\infty}) := \psi_1(g)S_{\infty}$. Recall the ascending HNN extension $\mathcal{G}(G) \ast_{\psi_{\infty}}$ has the presentation 
\begin{align*}
\mathcal{G}(G)~ \ast_{\psi_{\infty}} &:= \langle \mathcal{G}(G), T | T^{-1} g S_{\infty} T = \psi_{\infty}(gS_{\infty}) \text{ for all } g \in G \rangle \\
&= \langle \mathcal{G}(G), T | T^{-1} g S_{\infty}T  = \psi_1(g) S_{\infty} \text{ for all } g \in G \rangle. 
\end{align*}
The assignments $T \mapsto t$ and $g S_{\infty} \mapsto [g]$ extend to a well-defined isomorphism. 

Now, let us determine the kernel of $\phi_k$. First, it is not hard to see that if $g$ is a permutation of $\mathbb{N}$ for which there exists $d \in \mathbb{Z}$ such that $g(j) = j + d$ for all sufficiently large $j$, then necessarily $d=0$. Then $h \in H_n(G)$ belongs to the kernel if and only if $\phi_k(h) = [\id_{\mathbb{N}}]$ if and only if $[\sigma^{r_k} g_k \sigma^{-s_k}] = [\id_{\mathbb{N}}]$ if and only if $[g_k] = [\sigma^{s_k-r_k}]$ if and only if 
\[
g_k(j) = \sigma^{s_k-r_k}(j) = j + s_k -r_k
\]
for all sufficiently large $j \in \mathbb{N}$. This forces $s_k = r_k$, and hence $g_k(j) = j$ for all sufficiently large  $j \in \mathbb{N}$. Thus, $h$ belongs to the kernel if and only if $s_k = r_k$ and $g_k \in S_{\infty}$.  

Now, for the $n=1$ case, the sum balance condition implies $s:= s_1 = r_1$ and hence
\[
\phi_1(h) = t^s [g] t^{-s}.
\]
On the other hand, given $t^s [g] t^{-s}$ for some $s \ge 0$ and $g \in G$, we can define $h \in H_1(G)$ as
\[
h(j)=
\begin{cases}
j, & 1\le j\le s,\\[2mm]
g(j-s)+s, & j>s
\end{cases}
\]
and clearly $\phi_1(h) = t^s [g] t^{-s}$. Hence, the image is the directed union
\[
\image \phi_1 = \bigcup_{s \ge 0} t^s [G] t^{-s}.
\]
The shift-similar group $G$ is the subgroup of $H_1(G)$ where $i_h= 1$ and $s = 0$. Hence, the image of $G$ under $\phi_1$ is simply $[G]$. Note the element $[\sigma]$ does not lie in the image of either $G$ or $H_1(G)$ under $\phi_1$. Indeed, $\phi_1 (h) = [\sigma]$ if and only if $h(j) = \sigma (j) = j+1$ for all sufficiently large $j \in \mathbb{N}$, which is impossible because $h$ is a permutation of $\mathbb{N}$. As for the kernel, it is clearly just $S_{\infty}$ for either $G$ or $H_1(G)$. Let us now record these calculations.
\begin{theorem}
Let $G \le \sym (\mathbb{N})$ be a shift-similar group. For $n \ge 2$, the image of $\phi_k$ is the ascending HNN extension of $[G]$ with the torsion-free element $t := [\sigma]$, and the kernel is $\estab (R_k)$, the eventual pointwise stabilizer of the $k$-th ray. When $n=1$, the image of $\phi_1$ is the directed union
\[
\image \phi_1 = \bigcup_{s \ge 0} t^s [G] t^{-s}
\]
and the kernel is $S_{\infty}$. 
\end{theorem}
The homomorphism $\phi_k$ gives us a short exact sequence of $\estab(R_k)$, $H_n(G)$, and $\image \phi_k$ and hence a twisted crossed product decomposition of $L(H_n(G))$:
\[
L(H_n(G)) \cong L(\estab (R_k)) \rtimes \image \phi_k.
\]
  
Although our focus is shift-similar and Houghton-like groups, we will prove some results about a much broader class of subgroups of $\sym_{ends}(X_n)$, which we now define. Given subgroups $G_1, \dots , G_n \le \sym^*(\mathbb{N})$ of germs of almost permutations on $\mathbb{N}$, define 
\[
\Gamma_n(G_1, \dots , G_n) := \Phi^{-1}(G_1 \times \dots , G_n) = \{g \in \sym_{ends}(X_n) : \phi_k(g) \in G_k \text{ for all } k=1, \dots, n \}.
\]
If $G:=G_1 = \dots  = G_n$, we simply write $\Gamma_n(G)$ for $\Gamma_n(G, \dots , G)$. 

If $G_1, \dots, G_n \le \sym(\mathbb{N})$ are shift-similar groups, we write 
\[
H_n(G_1, \dots , G_n) := \Gamma_n(\langle [G_1], t \rangle , \langle [G_2],t \rangle , \dots , \langle [G_n], t \rangle).
\]
Thus, we can form a Houghton-like group where we use a different shift-similar group for each ray.  If $G := G_1 = \dots G_n$, we simply write $H_n(G) = \Gamma_n(\langle [G],t \rangle)$, which coincides with the Houghton-like groups introduced in \cite{mallery25}. We will also refer to $H_n(G_1, \dots , G_n)$ as Houghton-like groups. 

\subsection{Amenable (Near) Actions of (Almost) Permutation Groups}
Given a group $G \le \sym^*(\mathbb{N})$, it admits a well-defined natural action on the Boolean algebra $\mathcal{P}(\mathbb{N})/\fin$ defined as 
\[
[f] \cdot [A] = [f(A \cap \text{dom } f)]
\]
for $[f] \in G$ and $[A] \in \mathcal{P}(\mathbb{N})/\fin$. This is called a ``near" action of $G$ on $\mathbb{N}$ (see \cite{cornulier19} for a wealth of information on near actions). This near action is amenable if there exists a mean on the Boolean algebra $\mathcal{P}(\mathbb{N})/\fin$. That is, a map $m : \mathcal{P}(\mathbb{N})/\fin \to [0,1]$ such that $m ([\emptyset]) = 0$, $m ([\mathbb{N}]) = 1$, $m([A] \vee [B]) = m([A]) + m([B])$ whenever $[A] \wedge [B] = [\emptyset]$, and $m([f] \cdot [A]) = m([A])$ for $[f] \in G$ and $[A] \in \mathcal{P}(\mathbb{N})/\fin$.

\begin{theorem}\label{thm:amenablenearaction}
Let $\Gamma \le \sym_{ends}(X_n)$ be a countable group containing $S_{\infty}(X_n)$. The following are equivalent:
\begin{enumerate}[leftmargin=*]
\item The action of $\Gamma$ on $X_n$ is amenable.
\item The induced action of $\Gamma$ on $\beta X_n \setminus X_n$ preserves a regular Borel probability measure. 
\item For some $k=1, \dots , n$, the induced action of $\image \phi_k$ on $\beta \mathbb{N} \setminus \mathbb{N}$ preserves a regular Borel probability measure.
\item For some $k=1, \dots , n$, the natural near action of $\image \phi_k$ on $\mathbb{N}$ is amenable; that is, there is a $\image \phi_k$-invariant mean on $\mathcal{P}(\mathbb{N})/\fin$.
\item For some $k=1, \dots , n$, there exist non-empty finite sets $F_m$ such that $\min F_m \to \infty$, $|F_m| \to \infty$, and 
\[
\frac{|f(F_m) \triangle F_m|}{|F_m|} \to 0 
\]
as $m \to \infty$ for every almost permutation $f$ with $[f] \in \image \phi_k$. 
\end{enumerate}
\end{theorem}
\begin{proof}
$(1) \implies (2)$ Since $\Gamma$ acts amenably on $X_n$, there exists a $\Gamma$-invariant mean on $\ell^{\infty}(X_n)$. Because $S_{\infty}(X_n) \le \Gamma$, this mean must vanish on indicator functions of finite sets and hence it vanishes on all finitely supported functions, which are dense in $c_0(X_n)$. Thus, this mean induces a state on the quotient 
\[
\ell^{\infty}(X_n)/c_0(X_n) \cong C(\beta X_n \setminus X_n). 
\]
The Riesz representation theorem applied to this state yields a $\Gamma$-invariant regular Borel probability measure on $\beta X_n \setminus X_n$ for the induced action of $\Gamma$ on $\beta X_n \setminus X_n$. 

$(2) \implies (3)$ Note that $\beta X_n \setminus X_n$ can be identified with 
\[
\bigsqcup_{k=1}^{n} \beta \mathbb{N} \setminus \mathbb{N}.
\]
Since $\Gamma \le \sym_{ends}(X_n)$, each copy of $\beta \mathbb{N} \setminus \mathbb{N}$ is invariant under the induced action of $\Gamma$, and at least one copy, say the $k$-th copy, must have positive measure. Restrict the regular Borel probability measure to this copy and normalize it. If $g, g' \in \Gamma$ satisfy $\phi_k(g) = \phi_k(g')$, then $g$ and $g'$ agree on a cofinite subset of the $k$-th ray, and hence induce the same homeomorphism of the $k$-th copy of $\beta \mathbb{N} \setminus \mathbb{N}$. Thus, the induced action of $\Gamma$ factors through $\phi_k$ on this copy of $\beta \mathbb{N} \setminus \mathbb{N}$ and therefore the action of $\image \phi_k$ preserves this regular Borel probability measure on $\beta \mathbb{N} \setminus \mathbb{N}$.

$(3) \implies (4)$ Let $\mu$ be a regular Borel probability measure on $\beta \mathbb{N} \setminus \mathbb{N}$ invariant under the action of $\image \phi_k$. The map $[A] \mapsto A^* := \overline{A} \cap (\beta \mathbb{N} \setminus \mathbb{N})$ is an equivariant Boolean algebra isomorphism from $\mathcal{P}(\mathbb{N})/\fin$ to the Boolean algebra of clopen sets in $\beta \mathbb{N} \setminus \mathbb{N}$. Then 
\[
\nu ([A]) := \mu (A^*)
\]
defines a $\image \phi_k$-invariant mean on $\mathcal{P}(\mathbb{N})/\fin$.

$(4) \implies (5)$ Since $\image \phi$ is countable, \cite[Corollary 8.A.10]{cornulier19} gives us non-empty pairwise disjoint finite sets $B_j \subset \mathbb{N}$ satisfying the F\o lner condition for $\image \phi_k$. Define 
\[
F_m := \bigcup_{j=m}^{2m-1} B_j.
\]
Disjointness implies
\[
|F_m| = \sum_{j=m}^{2m-1} |B_j| \ge m
\]
so $|F_m| \to \infty$ as $m \to \infty$. If $M \in \mathbb{N}$, then because the $B_j$ are pairwise disjoint, each point of $\{1, \dots M \}$ belongs to at most one $B_j$. Hence, only finitely many $B_j$ meet $\{1, \dots , M \}$. Let $m_0$ be larger than any of these indices. Then if $m \ge m_0$, we have $B_j \cap \{1, \dots , M \} = \emptyset$ and thus $F_m \cap \{1, \dots , M\} = \emptyset$. Hence, $\min F_m > M$ for every $m \ge m_0$, proving $\min F_m \to \infty$ as $m \to \infty$. 

Finally, let $f$ be an almost permutation with $[f] \in \image \phi_k$. For all sufficiently large $m \in \mathbb{N}$, $f(F_m)$ is defined and hence 
\[
\frac{|f(F_m) \triangle F_m|}{|F_m|} \le \frac{\sum_{j=m}^{2m-1}|f(B_j) \triangle B_j|}{\sum_{j=m}^{2m-1} |B_j|} \le \max_{m \le j \le 2m-1} \frac{|f(B_j) \triangle B_j|}{|B_j|} 
\]
which converges to $0$ as $m \to \infty$. 

$(5) \implies (1)$ Define $E_m := \{k\} \times F_m$. Given $g \in \Gamma$, Lemma \ref{lemma:PhiHom} says there exists an almost permutation $f_{g,k}$ such that $g(k,j) = (k,f_{g,k}(j))$ for all sufficiently large $j$, and by definition $[f_{g,k}] \in \image \phi_k$. Since $\min F_m \to \infty$ as $m \to \infty$, we can pick $m_0 \in \mathbb{N}$ large enough so that if $m \ge m_0$, $g(k,j) = (k,f_{g,k}(j))$ for every $j \in F_m$, and therefore $g(E_m) = \{k\} \times f_{g,k}(F_m)$ which implies 
\[
\frac{|g(E_m) \triangle E_m|}{|E_m|} = \frac{|f_{g,k}(F_m) \triangle F_m|}{|F_m|} \to 0 
\]
as $m \to \infty$, thereby proving $\{E_m\}_{m \ge 1}$ is a F\o lner sequence for $\Gamma$ acting on $X_n$. 
\end{proof}

Condition (5) is precisely the condition we use to construct a central sequence in $L(\Gamma)$ to prove property Gamma. We call the sets in (5) F\o lner sets for the natural near action of $\image \phi_k$ on $\mathbb{N}$. 

The joint amenability condition of $G$ and $\sigma$ on $\mathbb{N}$ is exactly amenability of the near action of $\langle [G],[\sigma] \rangle$ on $\mathbb{N}$.

\section{Property Gamma and McDuff Factors}\label{sec:GammaMcDuff}

\subsection{Property Gamma Factors} 

In this section, we prove the following theorem for which Theorem \ref{mainthm:Gamma} is a special case.

\begin{theorem}[Theorem \ref{mainthm:Gamma}]\label{thm:PropertyGamma}
Let $\Gamma \le \sym_{ends}(X_n)$ be a countable ICC group (e.g., it contains $S_{\infty}(X_n)$). If the natural near action of $\image \phi_k$ on $\mathbb{N}$ is amenable for some $k=1, \dots , n$, then $L(\Gamma)$ is a property Gamma factor. 
\end{theorem}

As we hinted at in the introduction, we prove Theorem \ref{thm:PropertyGamma} by leveraging the subgroup $\estab(R_k)$. We first point out that it is rather easy to construct non-commuting central sequences in $L(\estab(R_k))$. Indeed, let $(a_m)_{m \ge 1}$ be a sequence of elements finitely supported on $R_k$ whose support slides out to infinity along $R_k$ as $m \to \infty$ (e.g., $a_m$ can transpose the points $(k,m)$ and $(k,m+1)$), and let $(b_m)_{m \ge 1}$ be a similar sequence where the supports of $a_m$ and $b_m$ have a single common point (e.g., $b_m$ could transpose $(k,m+1)$ and $(k,m+2)$). It is not hard to verify that these give non-commuting central sequences in $L(\estab (R_k))$. Note that this all holds even without the amenability hypothesis. Let us record this observation.

\begin{observation}
The von Neumann algebra $L(\estab(R_k))$ is a McDuff factor for $k=1, \dots , n$.
\end{observation}

However, these very natural and easily constructed central sequences in $L(\estab (R_k))$ do not lift to central sequences in $L(\Gamma)$. To see this, let $\Gamma = H_n$ with $n \ge 2$ and $k=1$, and define $f \in H_n$ to be the permutation which shifts every element on $R_{2}$ down one unit, sends $(2,1)$ to $(1,1)$, and shifts every element on $R_1$ up one element. If $(\lambda (a_m))_{m \ge 1}$ were to asymptotically commute with $\lambda (f)$, then $a_m$ would have to algebraically commute with $f$ for all sufficiently large $m \in \mathbb{N}$ and hence every power of $f$. Hence, for $m \in \mathbb{N}$ sufficiently large, this would imply 
\[
\supp (a_m) = \supp (f^j a_m f^{-j}) = f^j(\supp (a_m))
\]
for every $j \in \mathbb{Z}$. Given $x \in \supp (a_m)$, we would have $f^j(x) \in \supp (a_m)$ for every $j \in \mathbb{N}$. Since the points $f^j(x)$ are pairwise distinct for every $j \in \mathbb{N}$, this would cause the support of $a_m$ to blow up to infinity, contradicting the fact that $a_m$ is finitely supported. 

Let $\Gamma$ be as in Theorem \ref{thm:PropertyGamma}. By Theorem \ref{thm:amenablenearaction}, amenability of the natural near action of $\image \phi_k$ on $\mathbb{N}$ is equivalent to the existence of F\o lner sets $\{F_m\}_{m \ge 1}$ for this natural near action. Define 
\[
E_m := \{k\} \times F_m,
\]
and define $\sym_2(E_m)$ as all the two-cycles (transpositions) supported on $E_m$. Note that
\[
|\sym_2(E_m)| = \binom{|E_m|}{2} = \binom{|F_m|}{2}.
\]
For each $m \in \mathbb{N}$, define $x_m$ as
\begin{equation}\label{eq:UnboundedACsequence}
    x_m = |\sym_2(E_m)|^{-1/2} \sum_{a \in \sym_2(E_m)} \lambda (a).
\end{equation}
Because $x_m$ is essentially a sum over involutions, $x_m$ is self-adjoint, and it is not hard to see that 
\[
\norm{x_m}_2 = 1 
\]
for every $m \in \mathbb{N}$. This is clearly an element of $L(\estab (R_k))$, and it asymptotically commutes with everything in $L(\estab (R_k))$. Indeed, given $g \in \estab(R_k)$, the element $g$ acts by identity on $E_m$ for all sufficiently large $m \in \mathbb{N}$ because $\min F_m \to \infty$ as $m \to \infty$. Hence, $g$ will have support disjoint from everything in $\sym (E_m)$ and hence in $\sym_2(E_m)$. Therefore, $(x_m)_{m \ge 1}$ clearly asymptotically commutes with $\lambda (g)$, so $(x_m)_{m \ge 1}$ asymptotically commutes with everything in $L(\estab (R_k))$. 

Although $(x_m)_{m \ge 1}$ will asymptotically commute with everything in $L(H_n(G))$ and $L(G)$, this is not the central sequence we construct for $L(H_n(G))$. Unfortunately, it is not operator norm bounded. Indeed, on the one hand, it is not hard to see
\[
\norm{x_m}_{\infty} \le |\sym_2(E_m)|^{1/2} = \sqrt{ \binom{|F_m|}{2}}
\]
Now, consider the sequence of unit vectors
\[
\xi_m = |\sym (E_m)|^{-1/2} \sum_{a \in \sym(E_m)} \delta_a.
\]
For every $a \in \sym_2(E_m)$, we have $\lambda (a) \xi_m = \xi_m$, so $x_m \xi_m = |\sym_2(E_m)|^{1/2} \xi_m$ and therefore 
\[
\norm{x_m}_{\infty} \ge \norm{x_m \xi_m} = |\sym_2(E_m)|^{1/2},
\]
from which it follows that
\[
\norm{x_m}_{\infty} = |\sym_2(E_m)|^{1/2} = \sqrt{\frac{|F_m| \cdot (|F_m|-1)}{2}} \to \infty
\]
as $m \to \infty$. 

\looseness=-1 Although $(x_m)_{m \ge 1}$ is not operator norm bounded, it does asymptotically commute with everything in $L(\Gamma)$ (see Lemma \ref{lemma:ACsequenceforPropertyGamma} below). We use this to prove that, for a suitable $t \in \mathbb{R}$ chosen later, the exponentiated sequence
\begin{equation}\label{eq:AcSequenceBounded}
    u_m := e^{it x_m},
\end{equation}
also asymptotically commutes with everything in $L(\Gamma)$. Each $u_m$ is well-defined because $x_m$ is self-adjoint. Since exponentiation of self-adjoint elements produces unitaries, $(u_m)_{m \ge 1}$ is trivially operator norm bounded. Moreover, because von Neumann algebras are closed Borel functional calculus, it is a sequence in $L(\estab(R_k))$. We claim that this is a central sequence in $L(\Gamma)$. All we need to show is that it asymptotically commutes with everything in $L(\Gamma)$ and that it is not equivalent to a scalar sequence with respect to the two-norm.

To prove asymptotic commutativity, we need the following lemma.

\begin{lemma}\label{lemma:expACseq}
    Let $(M,tr)$ be a tracial von Neumann algebra and $t \in \mathbb{R}$. If $(x_m)_{m \geq 1}$ is a sequence of self-adjoint elements in $M$ that asymptotically commutes with everything in $M$, then the sequence $(u_m)_{m \ge 1}$ with $u_m := e^{itx_m}$ asymptotically commutes with everything in $M$.
\end{lemma}
\begin{proof}
Let $x,y$ be self-adjoint elements in $M$. First, we derive the inequality 
\[
\norm{e^{itx} - e^{ity}}_2 \le |t| \norm{x-y}_2.
\]
It is enough to prove it for $t \ge 0$ because if $t < 0$, then
\[
\norm{e^{itx} - e^{ity}}_2 = \norm{(e^{itx} - e^{ity})^*}_2 = \norm{e^{-itx}-e^{-ity}}_2. 
\]
For $0 \le s \le t$, define the differentiable operator valued function $f(s) = e^{isx} e^{i(t-s)y}$. By the product rule and chain rule, we have 
\begin{align*}
f'(s) &= \frac{d}{ds} [e^{isx}] e^{i(t-s)y} + e^{isx} \frac{d}{ds} [e^{i(t-s)y}] \\
&= i x e^{isx} e^{i(t-s)y} + e^{isx} (-iy e^{i(t-s)y}) \\
&= ie^{isx} x e^{i(t-s)y} -i e^{isx} y e^{i(t-s)y} \\
&= ie^{isx}(x - y) e^{i(t-s)y}. 
\end{align*}
Note we used the fact that $x$ and $e^{isx}$ commute, which can easily be seen from the power series representation of $e^{itx}$, or more generally, it follows from the fact that if $f$ is a complex-valued Borel function defined on the spectrum of the self-adjoint element $x$, then $x$ and $f(x)$ commute.  Integrating $f'(s)$ from $0$ to $t$ yields 
\[
e^{itx} - e^{ity} = i \int_{0}^{t} e^{isx} (x-y) e^{i(t-s)y} ds. 
\]
Taking the two-norm, we get
\begin{align*}
\norm{e^{itx} - e^{ity}}_2 &= \norm{i \int_{0}^{t} e^{isx} (x-y) e^{i(t-s)y} ds}_2 \\
&\le \int_{0}^{t} \norm{e^{isx} (x-y) e^{i(t-s)y}}_2 ds \\
&= \int_{0}^{t} \norm{x-y}_2 \\
&= |t| \cdot \norm{x-y}_2. 
\end{align*}

Now, let $(x_m)_{m \ge 1}$ be our asymptotically commutative sequence, and define $u_m := e^{itx_m}$. Because a von Neumann algebra is spanned by its unitary elements, it suffices to check asymptotic commutativity against a unitary. Let $v \in \mathcal{U}(M)$. Then 
\begin{align*}
v u_m v^{*} &= v e^{i t x_m} v^{*} \\
&=  v \sum_{k=0}^{\infty} \frac{(itx_m)^k}{k!} v^{*} \\
&= \sum_{k=0}^{\infty} \frac{(i t v x_m v^{*})^k}{k!} \\
&= e^{it v x_m v^{*}}. 
\end{align*}
Hence, 
\begin{align*}
\norm{v u_m v^{*} - u_m}_2 &= \norm{e^{it v x_m v^{*}} - e^{itx_m}}_2 \\
&\le |t| \cdot \norm{v x_m v^{*} - x_m}_2 
\end{align*}
which converges to $0$ as $m \to \infty$.
\end{proof}

Now, we argue that $(x_m)_{m \ge 1}$ asymptotically commutes with everything in $L(\Gamma)$.

\begin{lemma}\label{lemma:ACsequenceforPropertyGamma}
    The sequence $(x_m)_{m \ge 1}$ defined in Equation \eqref{eq:UnboundedACsequence} is a central sequence in $L(\Gamma)$.   
\end{lemma}

\begin{proof}
We already pointed out above that $(x_m)_{m \ge 1}$ asymptotically commutes with everything in $L(\estab (R_k))$, so by our lifting lemma (Lemma \ref{lemma:CentSeqLift1}), showing $(x_m)_{m \ge 1}$ asymptotically commutes with everything in $L(\Gamma)$ only requires checking that it is asymptotically invariant under the cocycle action of $\image \phi_k$. 

To this end, let $f$ be an almost permutation of $\mathbb{N}$ with $[f] \in \image \phi_k$, and let $g \in \Gamma$ be any lift of $[f]$ via $\phi_k$. From the definition of $\phi_k$, we have
\[
g(k,j) = (k,f(j)) 
\]
for all sufficiently large $j \in \mathbb{N}$. Since $\min F_m \to \infty$ as $m \to \infty$, we have that for all sufficiently large $m$ the equation $g(k,j) = (k,f(j))$ holds for every $j \in F_m$. Since $g$ is a lift of $[f]$ via $\phi_k$, the automorphism $\alpha_{[f]}$ is induced by conjugation by $\lambda (g)$. Hence,

\begin{align*}
\norm{\alpha_{[f]}(x_m) - x_m}_2^2 &= \norm{\lambda (g) x_m \lambda (g^{-1}) - x_m}_2 \\
&= \norm{x_m}_2^2 - 2 \text{Re } tr(\lambda (g) x_m \lambda (g^{-1}) x_m) + \norm{x_m}_2^2 \\
&= 2 \left(1 - |\sym_2(E_m)|^{-1} \sum_{a_1,a_2 \in \sym_2(E_m)} tr(\lambda (g a_1 g^{-1} a_2)) \right).
\end{align*}

Note $tr(\lambda (g a_1 g^{-1} a_2)) \neq 0$ (and hence equals $1$) if and only if $f a_1 f^{-1} a_2 = e$ if and only if $a_2 = fa_1 f^{-1}$, so we are summing over all $a$ such that 
\begin{align*}
a & \in \sym_2(E_m) \cap f \sym_2(E_m) f^{-1} \\
&= \sym_2(E_m) \cap \sym_2(f(E_m)) \\
&= \sym_2(E_m \cap f(E_m)) \\
&= \sym_2(\{k\} \times (F_m \cap f(F_m)). 
\end{align*}
Hence,
\begin{align*}
\norm{\alpha_{[f]}(x_m) - x_m}_2^2 &= 2 \left( 1 - \frac{\sym_2(\{1\} \times (F_m \cap f(F_m))}{|\sym_2(F_m)|} \right) \\ 
&= 2 \left(1 - \frac{ \binom{ |\{1\} \times (F_m \cap f(F_m))|}{2}}{ \binom{|F_m|}{2}} \right) \\
&= 2 \left(1 - \frac{|F_m \cap f (F_m)| \cdot (|F_m \cap f (F_m)|-1) }{|F_m| \cdot (|F_m|-1)} \right) \\
&= 2 \left( 1 - \frac{|F_m \cap f (F_m)|}{|F_m|} \cdot \frac{\frac{|F_m \cap f (F_m)|}{|F_m|} - \frac{1}{|F_m|}}{1 - \frac{1}{|F_m|}} \right). 
\end{align*}
Recall that if $f$ is an almost permutation of $\mathbb{N}$ with $F \dom (f)$ finite, then, since $f$ is injective,
\[
|f(F) \triangle F| = 2\bigl(|F| - |F \cap f(F)|\bigr).
\]
Consequently, for a sequence $\{F_m\}_{m\ge1}$ of finite subsets of $\mathbb{N}$ with $\min F_m \to \infty$ and $|F_m| \to \infty$, we have $F_m \subseteq \dom (f)$ for all sufficiently large $m$, and
\[
\frac{|F_m \cap f(F_m)|}{|F_m|} \to 1
\quad\Longleftrightarrow\quad
\frac{|f(F_m) \triangle F_m|}{|F_m|} \to 0.
\]
Therefore, because $\{F_m\}_{m \ge 1}$ are F\o lner sets for the natural near action of $\image \phi_k$ on $\mathbb{N}$, we have
\[
\frac{|F_m \cap f(F_m)|}{|F_m|} \to 1
\]
as $m \to \infty$. Hence, 
\[
\lim_{m \to \infty} \norm{\alpha_{[f]}(x_m) - x_m}_2 = 0.
\] 
\end{proof}

Finally, we prove $(u_m)_{m \ge 1}$ is not equivalent to any scalar sequence. 

\begin{lemma}\label{lemma:NonTivialACseq}
    The sequence $(u_m)_{m \ge 1}$ from Equation \eqref{eq:AcSequenceBounded} is not equivalent to any scalar sequence in $L(\Gamma)$. 
\end{lemma}

\begin{proof}
If $k$ is an odd integer, then modulo a scalar, $x_m^k$ is a sum over $\lambda (g)$ where each $g$ is a product of an odd number of two-cycles. Since a product of an odd number of two-cycles is never the identity, we have $tr(x_m^k)=0$ for every odd integer $k$. Hence, 
\begin{align*}
tr(u_m) &= tr(e^{itx_m}) \\
&= tr \left( \sum_{k=0}^{\infty} \frac{(it x_m)^k}{k!} \right) \\
&= \sum_{k=0}^{\infty} \frac{(i t)^k tr(x_m^k)}{k!} \\
&= \sum_{k=0}^{\infty} \frac{(-1)^k t^{2k} tr(x_m^{2k})}{(2k)!} \\
&= tr \left(\sum_{k=0}^{\infty} \frac{(-1)^{k} (t x_m)^{2k}}{(2k)!} \right) \\
&= tr(\cos(tx_m)). 
\end{align*}
Recall that for all $x \in \mathbb{R}$, we have 
\[
1 - \frac{x^2}{2} \le \cos x \le 1 - \frac{x^2}{2} + \frac{x^4}{24} 
\]
and therefore 
\[
1 - \frac{t^2 x_m^2}{2} \le \cos (tx_m) \le 1 - \frac{t^2 x_m^2}{2} + \frac{t^4 x_m^4}{24}.
\]
Taking the trace, we get 
\[
1 - \frac{t^2}{2} \le tr(\cos (tx_m)) \le 1 - \frac{t^2}{2} + \frac{t^4}{24} tr(x_m^4),
\]
where we used the fact that $tr(x_m^2) = 1$. Preferably, we pick $t$ such that $1 - \frac{t^2}{2} > 0$ or, equivalently, $|t| < \sqrt{2}$. Let us compute the fourth moment $tr(x_m^4)$: 
\begin{align*}
tr(x_m^4) &= |\sym_2(E_m)|^{-2} \sum_{a,b,c,d \in \sym_2(E_m)} tr(\lambda (abcd)) \\ \\
&= \frac{| \{(a,b,c,d) : abcd = e \}|}{ \binom{|F_m|}{2}^2}. 
\end{align*} 
To count all quadruples $(a,b,c,d) \in \sym_2(E_m)^4$ with $abcd=e$, we split the counting into two cases: (1) $a=b$, and (2) $a \neq b$. The latter case has two subcases. First, if $a=b$, then $abcd = e$ becomes $cd=e$ or $c = d^{-1} = d$. Hence, $b$ and $d$ are totally determined by $a$ and $c$, respectively. Once $a$ and $c$ are selected, there is exactly one option for $b$ and $d$, and $a$ and $c$ are independently chosen. Hence, there are $|\sym_2(E_m)| = \binom{|F_m|}{2}$ choices for each. Therefore, the number of quadruples $(a,b,c,d)$ falling under this case is  
\[
\binom{|F_m|}{2}^2
\]
all of which are guaranteed to satisfy the condition $abcd=e$. 

Now, suppose $a \neq b$. There are two possibilities: (1) the two-cycles $a$ and $b$ are disjoint, or (2) they have exactly one point in common. The two-cycle $a$ is of the form $((k,i)~(k,j))$ for $i,j \in F_m$ distinct, but we can clearly identify it with $(i~j)$ to aid in dealing with these subcases and counting. More formally, there is a clear isomorphism from $\sym (E_m) \to \sym(F_m)$ coming from the projection map which preserves cycle lengths.

First, suppose $a$ and $b$ are disjoint. Then $b = (k~\ell)$ for $k, \ell \in F_m$ distinct and different from $i,j \in F_m$. Because $a$ and $b$ are disjoint, they commute. Since $ab$ transposes $i$ and $j$ and $k$ and $\ell$, we get the equations 
\[
dc(i) = j, \quad dc(j)=i, \quad  dc(k)=\ell, \text{ and } \quad dc( \ell) = k.
\] 
Since $c$ only moves two points, and we are dealing with four points that are moved, it follows $c$ must fix at least one point. If $c(i)=i$, then the first equation above implies $d(i) = j$. Since $i$ and $j$ are distinct, this means $d$ transposes $i$ and $j$ and therefore $d= (i~j) = a$. From the equation $ab=dc$, we get $b=c$. If $c(j)=j$, then from the second equation above we get $d(j) = i$ and therefore $d= (i~j) = a$, from which we get $b=c$. If $c(k) = k$, then from the third equation above, we get $d(k) = \ell$ and hence $d = (k ~ \ell) = b$, from which we get $a=c$. Finally, if $c(\ell) = \ell$, then from the last equation above, we get $d(\ell) = k$ and hence $d=b$ and $a=c$. In any case, we get that if $a$ and $b$ are disjoint and $abcd=e$, then $\{a,b\} = \{c,d\}$. Conversely, if $a$ and $b$ are disjoint (hence commute) and $\{a,b\} = \{c,d\}$, then $abcd = e$. Hence, there are 
\[
\frac{|F_m| \cdot (|F_m| - 1)}{2}
\]
choices for $a$, and
\[
\frac{(|F_m|-2)(|F_m|-3)}{2} 
\]
choices for $b$. Finally, because $\{a,b\} = \{c,d\}$, there are only two choices for $(c,d)$, namely, $(c,d) = (a,b)$ or $(c,d) = (b,a)$. Hence, there are a total number of 
\begin{align*}
    2 \cdot \frac{|F_m| \cdot (|F_m| - 1)}{2} \cdot \frac{(|F_m|-2)(|F_m|-3)}{2} &=  \frac{|F_m| \cdot (|F_m| - 1)(|F_m|-2)(|F_m|-3)}{2} \\
    &= 12 \binom{|F_m|}{4}
\end{align*}
possible quadruples under this case. 

Finally, we consider the case when $a$ and $b$ are distinct but share exactly one point in their support. If $a$ and $b$ share exactly one point in common, then the product $ab$ is a unique three-cycle, and the equation $abcd=e$ can be written as $ab = dc$, meaning that both products $ab$ and $dc$ must evaluate to the same three-cycle (the product $dc$ is forced to be a three-cycle because $ab$ is). There are $\binom{|F_m|}{3}$ ways to pick a three-element subset of $F_m$, and two ways to build a three-cycle from it. Hence, there are 
\[
2 \cdot \binom{|F_m|}{3}
\] 
distinct three-cycles. Any given three-cycle can be factored as a product of two two-cycles in exactly three ways. Thus, for each three-cycle, there are three valid ordered pairs $(a,b)$ such that the product $ab$ equals that given three-cycle. Similarly, $dc$ must equal that same given three-cycle, so there are exactly three valid ordered pairs $(c,d)$. Because the choices for $(a,b)$ and $(c,d)$ are independent for a given fixed three-cycle, each three-cycle contributes $3 \times 3 = 9$ quadruples satisfying the condition. Multiplying this by the total number of three-cycles, we conclude that under this case, there are  
\[
9 \cdot 2 \cdot \binom{|F_m|}{3} = 18 \binom{|F_m|}{3}
\]
possible quadruples $(a,b,c,d)$ with $a$ and $b$ sharing exactly one point in common and $abcd=e$.

These are disjoint cases, so 
\[
|\{(a,b,c,d) : abcd=e\}| = \binom{|F_m|}{2}^2 + 12 \binom{|F_m|}{4} + 18 \binom{|F_m|}{3}.
\]
Now, let $N_m := |F_m|$, so $N_m \to \infty$ as $m \to \infty$. Dividing by $\binom{N_m}{2}^2$ gives us the trace of $x_m^4$: 
\begin{align*}
tr(x_m^4) &= 1 + \frac{12 \binom{N_m}{4} + 18\binom{N_m}{3}}{\binom{N_m}{2}^2} \\
&= 1 + \frac{12 \left[ \frac{N_m(N_m-1)(N_m-2)(N_m-3)}{24} \right] + 18 \left[ \frac{N_m(N_m-1)(N_m-2)}{6} \right]}{\left[ \frac{N_m(N_m-1)}{2} \right]^2} \\
&= 1 + \frac{\frac{1}{2}N_m(N_m-1)(N_m-2)(N_m-3) + 3N_m(N_m-1)(N_m-2)}{\frac{1}{4}N_m^2(N_m-1)^2} \\
\end{align*}
\begin{align*}
&= 1 + \frac{\frac{1}{2}N_m(N_m-1)(N_m-2) \left[ (N_m-3) + 6 \right]}{\frac{1}{4}N_m^2(N_m-1)^2} \\
&= 1 + \frac{2(N_m-2)(N_m+3)}{N_m(N_m-1)} \\
&= 1 + \frac{2N_m^2 + 2N_m - 12}{N_m^2 - N_m}.
\end{align*}
Taking the limit as $m \to \infty$, we obtain
\[
\lim_{m \to \infty} tr(x_m^4) = 3.
\]
Now, note that 
\[
\norm{u_m - tr(u_m)1}_2^2 =  1 - |tr(u_m)|^2. 
\]
Also recall from above that $tr(u_m) = tr(\cos (tx_m))$ and that
\[
1 - \frac{t^2}{2} \le tr(\cos(tx_m)) \le 1 - \frac{t^2}{2} + \frac{t^4}{24} tr(x_m^4).
\]
Since $|t| < \sqrt{2}$, the lower bound is positive and hence $|tr(u_m)| = tr(u_m) = tr(\cos(tx_m))$. Define $g(s) := 1 - \left(1 - \frac{t^2}{2} + \frac{t^4}{24}s\right)^2$, so that the computation above reads
\[
\norm{u_m - tr(u_m)1}_2^2
= 1 - |tr(u_m)|^2
= 1 - |tr(\cos(tx_m))|^2
\ge g(tr(x_m^4))
\]
for every $m \in \mathbb{N}$. Since $g$ is continuous and $tr(x_m^4) \to 3$, taking the limit inferior gives
\[
\liminf_{m \to \infty} \norm{u_m - tr(u_m)1}_2^2
\ge \liminf_{m \to \infty} g(tr(x_m^4))
= g(3)
= 1 - \left(1 - \frac{t^2}{2} + \frac{t^4}{8}\right)^2.
\]
If $t=1$, then
\[
\liminf_{m \to \infty} \norm{u_m - tr(u_m)1}_2 \ge \frac{\sqrt{39}}{8} > 0,
\]
which proves that $(u_m)_{m \ge 1}$ is not equivalent to a scalar sequence in the two-norm.

\end{proof}

Here are some immediate corollaries we obtain.

\begin{corollary}\label{cor:almostpermgamma}
Let $G_1, \dots , G_n \le \sym^*(\mathbb{N})$ be countable groups of germs of almost permutations on $\mathbb{N}$. If the natural near action of $G_k$ on $\mathbb{N}$ is amenable for some $k=1, \dots , n$, then $L(\Gamma_n(G_1, \dots , G_n))$ is a property Gamma factor. 
\end{corollary}
\begin{proof}
It is clear that $S_{\infty}(X_n) \le \Gamma_n(G_1, \dots , G_n)$. By assumption, there exists $k=1, \dots , n$ such that the natural near action of $G_k = \image \phi_k |_{\Gamma_n(G_1, \dots , G_n)}$ on $\mathbb{N}$ is amenable. Hence, the conclusion follows immediately from Theorem \ref{thm:PropertyGamma}.
\end{proof}

\begin{corollary}\label{cor:shiftsimilargamma}
\leavevmode
\begin{enumerate}
\item Let $G_1, \dots , G_n \le \sym(\mathbb{N})$ be shift-similar groups such that $G_k$ and $\sigma$ act amenably on $\mathbb{N}$ for some $k=1, \dots , n$. Then 
\[
L(H_n(G_1, \dots , G_n))
\]
is a property Gamma factor. 

\item If $G \le \sym (\mathbb{N})$ is an infinite shift-similar group such that $G$ and $\sigma$ act amenably on $\mathbb{N}$, then 
\[
L(H_n(G))
\]
is a property Gamma factor. 

\item Finally, if $G \le \sym (\mathbb{N})$ is an infinite shift-similar group acting amenably on $\mathbb{N}$, then 
\[
L(G) \: \: \text{and} \: \: L(H_1(G))
\]
are property Gamma factors. 
\end{enumerate}
\end{corollary}
\begin{proof}
Proof of (1): If $\mu$ is a finitely additive probability measure on $\mathbb{N}$ invariant under $G_k$ and $\sigma$, then 
\[
\mu (\mathbb{N}) = \mu (\{1\} \cup \sigma (\mathbb{N})) = \mu (\{1\}) + \mu (\sigma (\mathbb{N})) = \mu (\{1\}) + \mu (\mathbb{N}) 
\]
\looseness = -1 and hence $\mu (\{1\}) = 0$. By induction, $\mu (\{n\}) = 0$ for every $n \in \mathbb{N}$, and hence $\mu (F) = 0$ for every finite set $F \subset \mathbb{N}$. Therefore, $\mu$ gives to a mean on $\mathcal{P}(\mathbb{N})/\fin$ preserved by $\langle [G_k],[\sigma] \rangle$. Since 
\[
H_n(G_1, \dots , G_n) = \Gamma_n( \langle [G_1],t\rangle, \dots , \langle [G_n], t \rangle ),
\]
the preceding corollary implies
\[
L(H_n(G_1, \dots , G_n))
\]
is a property Gamma factor. 

Proof of (2): This follows immediately from part (1) by taking $G = G_1 = \dots G_n$.

Proof of (3): If $G$ is infinite and acts amenably on $\mathbb{N} = X_1$, then we can immediately conclude $L(G)$ has property Gamma. Moreover, since $S_{\infty} \le G$, $S_{\infty}$-invariance implies the measure vanishes on all finite subsets of $\mathbb{N}$. Hence, the $G$-invariant measure on $\mathbb{N}$ gives rise to a $[G]$-invariant mean on $\mathcal{P}(\mathbb{N})/\fin$. Because $[G]$ preserves a mean on $\mathcal{P}(\mathbb{N})/\fin$, each $t^s[G]t^{-s}$ also preserves a translation of this mean.

Shift-similarity implies that if $r \ge s$, then $t^s [G] t^{-s} \subseteq t^r [G] t^{-r}$. If $M_s$ denotes the set of $t^s[G]t^{-s}$-invariant means on $\mathcal{P}(\mathbb{N})/\fin$, then $M_{s+1} \subseteq M_s$. These are non-empty compact subsets of the function space $[0,1]^{\mathcal{P}(\mathbb{N})/\fin}$, so their intersection is non-empty. A mean in the intersection is invariant under the directed union
\[
\bigcup_{s \ge 0} t^s [G]t^{-s} = \image \phi_1.
\]
Hence, the natural near action of $\image \phi_1$ on $\mathbb{N}$ is amenable and therefore $L(H_1(G))$ is a property Gamma factor.
\end{proof}

Given that the hypotheses apply to such a broad class of groups, it seems unreasonable to expect $L(\Gamma)$ to have the McDuff property in general. Indeed, even the Houghton-like groups $H_n(G)$ alone are very broad and general. Initially, we attempted to prove groups satisfying the conditions of Theorem \ref{thm:PropertyGamma} yield McDuff factors, but they resisted every attempt at constructing a central sequence which does not commute with the one we constructed above. We do know there are cases where this is in principle possible. For example, the classical Houghton group $H_n$ is amenable and ICC, so $L(H_n)$ is the hyperfinite II$_1$ factor whose central sequence algebra has trivial center. Thus, for abstract reasons, we know there exists a central sequence in $L(H_n)$ which does not commute with the central sequence $(u_m)_{m \ge 1}$ we constructed above. In the next section, we show that if $\image \phi_k$ is amenable for some $k=1, \dots , n$, then $L(\Gamma)$ is McDuff. However, this is not an ideal condition for the case of $H_n(G)$ because it would imply $H_n(G)$ itself is amenable, which puts us in the trivial case. We leave it to future work to find a alternative, non-trivial conditions which guarantee that $L(H_n (G))$, or \(L(\Gamma)\) more generally, is McDuff.

As we noted in the introduction, one might wonder how these results relate to those obtained in \cite{kerr25}, particularly their Theorem C. In that paper, Kerr and Petrakos proved that if a group $\Gamma$ acts amenably on a set $Y$ with no finite orbits, and $H$ is any non-trivial ICC group, then the generalized wreath product 
\[
H \wr_{Y} \Gamma := \bigoplus_{y \in Y } H_y \rtimes \Gamma,
\]
where $H_y = H$ for every $y \in Y$, yields a type II$_1$ McDuff factor. In our setup, our groups contain the finitely supported permutations and therefore act transitively on $X_n$, so there is a single infinite orbit. Additionally, we have an amenable action condition. We can, however, prove our groups cannot decompose as a generalized wreath product, thereby proving our results are distinct from those in \cite{kerr25}. 

\begin{theorem}\label{thm:NonWreathProduct}
Let $X$ be any infinite set, and let $G$ be any subgroup of $\sym (X)$ containing $A_{\infty}(X)$. Then $G$ cannot be decomposed as a non-trivial generalized wreath product (i.e., not of the form $H \wr_{Y} \Gamma$ with $|Y| \ge 2$ or $\Gamma$ non-trivial).
\end{theorem}
\begin{proof}
By way of contradiction, suppose $G$ is isomorphic to $H \wr_{Y} \Gamma$ for some groups $H$ and $\Gamma$, and some set $Y$ with a $\Gamma$-action. Then $N := \bigoplus_{y \in Y} H_y$ is a normal subgroup of $G$. Since $A_{\infty}(X)$ is infinite simple, the normal subgroup $A_{\infty}(X) \cap N$ of $A_{\infty}(X)$ is either trivial or equals all of $A_{\infty}(X)$. First, suppose $A_{\infty}(X) \cap N$ is trivial. Let $n \in N$ be a non-trivial element of $N$, and let $g \in A_{\infty}(X)$ be arbitrary. On the one hand, $[n,g] = (n^{-1} g^{-1}n)g$ shows that $[n,g]$ belongs to $A_{\infty}(X)$ because it is a normal subgroup. On the other hand, $[n,g] = n^{-1}(g^{-1} n g)$ shows that $[n,g]$ belongs to $N$ for the same reason. Hence, $[n,g] \in A_{\infty}(X) \cap N = \{e\}$ and therefore $n$ commutes with all of $A_{\infty}(X)$. But this is impossible because $A_{\infty}(X)$ has trivial centralizer in $\sym (X)$.

Now, let us consider the case when $A_{\infty}(X) \cap N = A_{\infty}(X)$. This implies $A_{\infty}(X)$ is a subgroup of $N = \bigoplus_{y \in Y} H_y$. This implies $N$ is non-trivial, and therefore $H$ and hence each $H_y$ is non-trivial. Let $y \in Y$ be arbitrary. Because $A_{\infty}(X)$ and $H_y$ are both normal in $N$, we have $[A_{\infty}(X),H_y] \le A_{\infty}(X)$. If $A_{\infty}(X) \cap H_y = \{e\}$, then $[A_{\infty}(X),H_y] = \{e\}$ and therefore everything in the non-trivial group $H_y$ commutes with $A_{\infty}(X)$, which is not possible because the centralizer of $A_{\infty}(X)$ in $\sym (X)$ is trivial. Hence, there exists $y_0 \in Y$ such that $A_{\infty}(X) \cap H_{y_0} \neq \{e\}$. This means that $A_{\infty}(X) \cap H_{y_0}$ is a non-trivial normal subgroup of the infinite simple group $A_{\infty}(X)$, therefore it must be the case that $A_{\infty}(X) = A_{\infty}(X) \cap H_{y_0} \le H_{y_0}$. This says that $A_{\infty}(X)$ is entirely contained in a single summand of $N = \bigoplus_{y \in Y} H_y$. Because the summands in a direct sum pairwise commute, and because all of the other summands are non-trivial (indeed, they are all isomorphic to $H \neq \{e\}$), this forces $Y = \{y_0\}$. This means the direct sum degenerates to the single group $H_{y_0}$. But since $Y = \{y_0\}$, and the only action of a group on a one point space is trivial, this means the wreath product further degenerates to the direct product $H_{y_0} \times \Gamma$. Since the factors in a direct product commute, we know $\Gamma$ commutes with $H_{y_0}$ and hence its subgroup $A_{\infty}(X)$. But this contradicts the centralizer fact about $A_{\infty}(X)$ unless $\Gamma$ is the trivial group. Hence, $\Gamma$ must necessarily be the trivial group, which shows that the generalized wreath product completely degenerates to $H \wr_{Y} \Gamma = H$, which is a trivial generalized wreath product. 
\end{proof}

\subsection{McDuff Factors}

Given a group countable subgroup $\Gamma$ of $\sym_{ends}(X_n)$ containing $S_{\infty}(X_n)$, if we further assume $\image \phi_k$ is amenable for at least one $k=1, \dots , n$, then we can actually conclude that $L(\Gamma)$ is a McDuff factor as we prove in this section. 

Much more is known about central sequences in crossed product von Neumann algebras when the acting group is amenable. Indeed, we will appeal to \cite{jolissaint1998central} and prove a stronger property. Given a II$_1$ factor $M$, we say $M$ has the \emph{relative McDuff property} with respect to a subfactor $N \subset M$ if $M'\cap N^\omega$ is non-abelian. 

\begin{theorem}[Theorem \ref{mainthm:McDuff}]\label{thm:McDuff}
Let $\Gamma \le \sym_{ends}(X_n)$ be a countable group containing $S_{\infty}(X_n)$. If $\image \phi_k$ is amenable for some $k=1, \dots , n$, then $L(\estab(R_k)) \subseteq L(\Gamma)$ has the relative McDuff property. In particular, $L(\Gamma)$ is a McDuff factor.  
\end{theorem}
\begin{proof}
As before, we have a short exact sequence involving $\estab(R_k)$, $\Gamma$, and $\image \phi_k$, and therefore we get the following crossed product decomposition:
\[
L(\Gamma) = L(\estab(R_k)) \rtimes \image \phi_k
\]
Because $\image \phi_k$ is amenable, \cite[Proposition 2.6]{jolissaint1998central} tells us it suffices to prove the cocycle action of $\image \phi_k$ on $L(\estab(R_k))$ is centrally free to prove the relative McDuff property. 

To this end, let $q \in \image \phi_k$ be non-trivial, and let $g \in \Gamma$ be a lift of $q$ via $\phi_k$. Since $q$ is non-trivial,  $g$ is not in $\estab(R_k)$. Hence, there exists an infinite set of distinct points $\{j_m\}_{m \ge 1} \subseteq \mathbb{N}$ such that $g(k,j_m) \neq (k,j_m)$. Since $g \in \sym_{ends}(X_n)$, we can remove finitely many initial points (if necessary) so that $g(k,j_m) \in R_k$ for every $m \in \mathbb{N}$. 

Now, let $p_m := (k,j_m)$ and pick any $q_m \notin \{p_m,g^{-1}(p_m)\}$ further along $R_k$ than $x_m$. Define $a_m := (p_m ~ q_m)$ and $x_m := \lambda (a_m)$. By the remarks in the previous section, $(x_m)_{m \ge 1}$ clearly forms a central sequence in $L(\estab(R_k))$. Moreover, 
\[
g a_m g^{-1} = g (x_m ~ y_m) g^{-1} = (g(x_m) ~ g(y_m))
\]
which cannot be equal to $a_m = (x_m ~ y_m)$ otherwise $x_m = g(x_m)$ or $x_m = g(y_m)$. Hence,
\[
\lim_{m \to \infty} \norm{\alpha_{q}(x_m) - x_m}_2 = \lim_{m \to \infty} \norm{\lambda (g) \lambda (a_m) \lambda (g^{-1}) -  \lambda (a_m)}_2 = \sqrt{2} \neq 0
\]
\end{proof}

The above theorem applies to the groups $\Gamma_n(G_1, \dots , G_n)$ for $G_1, \dots , G_n \le \sym^*(\mathbb{N})$ with at least one $G_k$ amenable, but the other groups can be as wild as one wishes. These will give McDuff factors. In particular, if $G_1, \dots , G_n \le \sym (\mathbb{N})$ are actually shift-similar groups with at least one amenable, then $H_n(G_1, \dots,G_n)$ yields a McDuff factor. 

\section{Non-inner Amenable Shift-similar and Houghton-like Groups}\label{sec:Non-innerAmean}

Before we look at examples of shift-similar groups $G$ where $G$ acts amenably on $\mathbb{N}$ with $\sigma$, let us explore what happens when a shift-similar group $G$ and $\sigma$ no longer satisfy this joint amenability condition. In this section, we construct a shift-similar group $G$ which cannot act amenably on $\mathbb{N}$ and is such that neither $G$ nor $H_n(G)$ is inner amenable. As a consequence, neither $L(G)$ nor $L(H_n(G))$ has property Gamma.

The main ingredient in constructing such a shift-similar group is a property (T) group $\Gamma$ acting on some countable set $X$ without finite orbits  with the property that for any finite-index subgroup $H$ in $\Gamma$, the centralizer of $H$ in $\sym (X)$ is trivial. For property (T) groups, acting  amenably is equivalent to having at least one finite orbit \cite[Lemma 4.2]{GlasnerMonod2007}. Then take any bijection from $X$ to $\mathbb{N}$ which will give rise to an isomorphism from $\sym (X)$ to $\sym (\mathbb{N})$. Then we push $\Gamma$ into $\sym (\mathbb{N})$ via this isomorphism, and define $G$ as the shift-similar closure of this copy of $\Gamma$ in $\sym (\mathbb{N})$. Because $\Gamma$ acts non-amenably $X$, the copy of $\Gamma$ in $\sym (\mathbb{N})$ cannot act amenably on $\mathbb{N}$ and hence neither can $G$. In particular, $G$ and $\sigma$ together cannot act amenably on $\mathbb{N}$. And we claim that neither $G$ nor $H_n(G)$ is inner amenable. An example of such a property (T) group with this property is the projective special linear group $\psl_k(\mathbb{Z})$ acting on the rational projective space $P^{k}(\mathbb{Q})$.  But before we get into the details of proving neither $G$ nor $H_n(G)$ are inner amenable, let us verify that $\psl_k(\mathbb{Z})$ acting on the rational projective space $P^{k}(\mathbb{Q})$ has the desired property. 

\begin{proposition}
    The natural action of $\psl_k(\mathbb{Z})$  on $P^{k}(\mathbb{Q})$
    has no finite orbits and any finite index subgroup $H$ of $\psl_k(\mathbb{Z})$ has trivial centralizer in $\sym (X)$.
 \end{proposition}

\begin{proof}
Let $H \le \psl_k(\mathbb{Z})$ be a finite-index subgroup, and let $\phi$ be an bijection of $P^{k}(\mathbb{Q})$ commuting with $H$. Let $N$ be the normal core of $H$. Then $N$ is also finite-index in $\psl_k(\mathbb{Z})$, and $\phi$ commutes with $N$. We will argue that $\phi$ must be the identity. 

To this end, let $x \in P^{k}(\mathbb{Q})$. This point can be represented by a primitive vector $v \in \mathbb{Z}^k$; i.e., $x = [v]$. The primitive vector $v$ can be extended to a $\mathbb{Z}$-basis for $\mathbb{Z}^k$. By Bezout's identity and column operations, we can construct a matrix $g \in \SL_k(\mathbb{Z})$ whose first column is $v$. This implies that $x = g \cdot [e_1]$, where $e_1$ is the first standard basis vector. We note in passing that this proves the action is transitive and therefore is without finite orbits.  

Now, let $U$ be all the matrices in $\SL_k(\mathbb{Z})$ of the form 
\[
\begin{bmatrix} 1 & u^T \\ 0 & I_{k-1} \end{bmatrix}
\]
for $u \in \mathbb{Z}^{k-1}$. 

The natural projection $ \SL_k(\mathbb{Z}) \to \psl_k(\mathbb{Z})$ restricts to an injection on $U$. So, we will identify $U$ with its image under the projection. Moreover, $U$ is isomorphic to $\mathbb{Z}^{k-1}$ by simply sending the matrix of the above form to $u \in \mathbb{Z}^{k-1}$. Since it is clear that every matrx in $U$ stabilizes $e_1$, we have $U \le  \stab([e_1])$.

More than that, $[e_1]$ is the only point that $U$ stabilizes. Indeed, let $[w] \in P^{k}(\mathbb{Q})$ be a point stabilized by $U$. By clearing denominators, we can assume $w \in \mathbb{Z}^k$. Then for every $A \in U$, we have $A \cdot [w] = [w]$ or $[Aw] = [w]$ which is equivalent to $Aw = \lambda w$ for some $\lambda \in \mathbb{Q}$ non-zero. Decompose $w$ as 
\[
w = (w_1, \tilde{w})^T,
\]
where $w_1 \in \mathbb{Z}$ and $\tilde{w} \in \mathbb{Z}^{k-1}$. Define $A \in U$ with $u = \tilde{w}$. Then $Aw = \lambda w$ gives us the equations 
\[
w_1 + \norm{\tilde{w}}_2^2 = \lambda w_1 \quad \text{ and } \quad \tilde{w} = \lambda \tilde{w}.
\]
The latter equation implies either $\lambda = 1$ or $\tilde{w}=0$, and $\lambda = 1$ combined with the first equation implies $\tilde{w}=0$. In either case, we get $\tilde{w} = 0$. and therefore $w = w_1 e_1$ which implies $[w] = [e_1]$.  

Because 
\[
|U : U \cap N | \le |PSL_k(\mathbb{Z}) : N | \le |\psl_k(\mathbb{Z}) : H|! < \infty,
\]
we have that $U \cap N$ is a finite-index subgroup of $U \cong \mathbb{Z}^{k-1}$. Every finite-index subgroup of $\mathbb{Z}^{k-1}$ contains a subgroup of the form $m \mathbb{Z}^{k-1}$ for some $m \in \mathbb{N}$. Hence, $U \cap N$ contains the image of $m \mathbb{Z}^{k-1}$ under the isomorphism from $\mathbb{Z}^{k-1}$ to $U$. The image of $m \mathbb{Z}^{k-1}$, which we denote $U_m$, consists of all matrices of the form 
\[
\begin{bmatrix} 1 & mu^T \\ 0 & I_{k-1} \end{bmatrix}
\]
for $u \in \mathbb{Z}^{k-1}$. Then $U_m \le N \cap U \le H \cap U$, so $U_m$ is a subgroup of both $\stab([e_1])$ and $H$. Conjugating by $g$, 
\[
g U_m g^{-1} \le gNg^{-1} \cap gUg^{-1} \le H
\]
and every element of $gU_m g^{-1}$ fixes $x = g \cdot [e_1]$ because 
\[
g \stab([e_1]) g^{-1} = \stab(g \cdot [e_1]) = \stab(x).
\]
Just as with $U$, the same argument given above shows that $[e_1]$ is the unique fixed point of $U_m$. Hence, $x = g \cdot [e_1]$ is the unique fixed point of $gU_mg^{-1}$.

Since $\phi$ commutes with $H$, it in particular commutes with $g U_m g^{-1}$. That is, $\phi gAg^{-1} = gAg^{-1} \phi$ for every $A \in U_m$. Evaluating this at $x = g \cdot [e_1]$, we get 
\begin{align*}
gAg^{-1} \phi(x) &= \phi gAg^{-1}(x) \\
&= \phi gAg^{-1}(g \cdot [e_1]) \\
&= \phi gA([e_1]) \\
&= \phi(g \cdot [e_1]) \\
&= \phi (x) 
\end{align*}
Since $gAg^{-1}$ is an arbitrary element of $gU_mg^{-1}$ which has unique fixed point $x$, it follows that $\phi (x) = x$.
\end{proof}

Now, let $G \le \sym (\mathbb{N})$ shift-similar closure of any property (T) group $\Gamma$ acting on $\mathbb{N}$ by permutations with no finite orbits and with the finite-index centralizer property. Before we prove neither $G$ nor $H_n(G)$ is inner amenable, we need the following lemma on groups containing a relatively ICC property (T) subgroup. We say $H \leq G$ is \emph{relatively ICC} if $\{hgh^{-1} : h \in H \}$ is infinite for every $g \in G \setminus \{e\}$.

\begin{lemma}\label{lemma:relativelyICC}
Let $G$ be a group containing a relatively ICC property (T) group $H$. Then $G$ is non-inner amenable.
\end{lemma}
\begin{proof}
By way of contradiction, suppose $G$ is inner amenable. Then there exists a sequence of unit vectors $\{\xi_n \}_{n \ge 1}$ in $\ell^2(G \setminus \{e\})$ such that $\norm{ \pi(g) \xi_n - \xi_n }_2 \to 0$ as $n \to \infty$ for every $g \in G$, where $\pi : G \to \mathcal{U}(\ell^2(G  \setminus \{e\})$ is the unitary representation coming from the conjugation action of $G$ on $G \setminus \{e\}$. When $\pi$ is restricted to $H$, this sequence of vectors clearly form almost invariant vectors for $H$ with this unitary representation. Because $H$ has property (T), there must exist a bona fide invariant non-zero vector $\xi \in \ell^2(G \setminus \{e\})$ such that $\pi (h) \xi = \xi$ for every $h \in G$. 

Since $\xi$ is non-zero, there exists $x_0 \in G \setminus \{e\}$ such that $\xi (x_0) \neq 0$. Because $x_0$ is a non-trivial element $G$, it has infinitely many $H$-conjugates. The invariance condition that $\pi (g) \xi = \xi$ holds for every $g \in H$ implies $\xi (gx_0g^{-1}) = \xi (x_0)$ for every $g \in H$ or $\xi(y) = \xi (x_0)$ for every $y \in x_0^H$, which denotes the $H$-conjugacy class of $x_0$. Since $x_0$ has infinitely many $H$-conjugates, we have 
\[
\sum_{g \in G \setminus \{e\}} |\xi (g)|^2 \ge \sum_{y \in x_0^H} |\xi (y)|^2 = |\xi (x_0)|^2 \sum_{y \in x_0^H} = \infty,
\]
which contradicts square summability of $\xi$. Hence, $G$ cannot be inner amenable. 
\end{proof}

\begin{theorem}[Theorem \ref{mainthm:NonInnerAmean}]\label{thm:NonInnerAmeanShift}
    Let $G \le \sym(\mathbb{N})$ be the shift-similar closure of any property (T) group $\Gamma \le \sym (\mathbb{N})$ with the finite-index-centralizer property that acts on $\mathbb{N}$ without finite orbits. Then $G$ is non-inner amenable.
\end{theorem}
\begin{proof}
Let $g \in G \setminus \{e\}$. If $C_{\Gamma}(g)$ were finite-index in $\Gamma$, then $g \in C_{\sym (\mathbb{N})} (C_\Gamma (g)) = \{e\}$, which is a contradiction. Hence, $C_{\Gamma}(g)$ has infinite index in $\Gamma$ which is equivalent to $g$ having infinitely many $\Gamma$-conjugates. Non-inner amenability of $G$ now follows from Lemma \ref{lemma:relativelyICC}. Since this holds for any non-trivial element in $G$, we get something strictly stronger than the relative ICC property with respect to $\Gamma$. 
\end{proof}

In fact, Theorem \ref{thm:NonInnerAmeanShift} applies to any subgroup of $\sym (\mathbb N)$ containing $\Gamma$. Before we prove $H_n(G)$ is non-inner amenable, we need the following basic lemma. 

\begin{lemma}
If $G$ is a group that acts on $X$ without finite orbits, then any finite-index subgroup of $G$ acts on $X$ without finite orbits. 
\end{lemma}
\begin{proof}
We prove the lemma contrapositively. Let $H$ be a finite-index subgroup of $G$, and suppose there exists $x \in X$ such that the $H$-orbit of $x$ is finite. This is equivalent to the stabilizer index $|H : \stab_H(x)|$ being finite. Because $\stab_H(x) \subseteq \stab_G(x) \subseteq G$, the tower property of group indices allows us to bound the index of the $G$-stabilizer:
\[
|G : \stab_G(x)| \le |G : \stab_H(x)| = |G : H| \cdot |H : \stab_H(x)|.
\]
Since both factors on the right-hand side are finite by hypothesis, their product is finite.  Thus, $|G : \stab_G(x)| < \infty$, which means the $G$-orbit of $x$ is finite.
\end{proof}

Since $\psl_k(\mathbb{Z})$ acts transitively on $P^{k}(\mathbb{Q})$ (i.e., has a single infinite orbit), any finite-index subgroup of $\psl_k(\mathbb{Z})$ acts on $P^{k}(\mathbb{Q})$ without finite orbits. These facts about orbits can, of course, all be pushed over to $\mathbb{N}$ via any bijection.

Now, let us prove $H_n(G)$ is non-inner amenable. 

\begin{theorem}[Theorem \ref{mainthm:NonInnerAmean}]\label{thm:NonInnerAmean}
    Let $G \le \sym(\mathbb{N})$ be the shift-similar closure of any property (T) group $\Gamma$ with the finite-index-centralizer property that acts on $\mathbb{N}$ without finite orbits. Then the Houghton-like group $H_n(G)$ is non-inner amenable.
\end{theorem}

\begin{proof}
First, we define a map $\Delta^n : \Gamma \to \sym_{ends}(X_n)$ as $\Delta^n(g)(k,j) = (k,g(j))$ for all $k=1, \dots , n$ and $j \in \mathbb{N}$. This is an injective homomorphism whose image isomorphic to the usual $n$-fold diagonal group of $\Gamma$. We claim that $\Delta^n(\Gamma)$ has the finite-index centralizer property in $\sym_{ends}(X_n)$. Recall that the finite-index subgroups in the $n$-fold diagonal group $\Delta^n(\Gamma)$ are of the form $\Delta^n(H)$ for some finite-index subgroup $H$ in $\Gamma$. 

Now, suppose $g \in \sym_{ends}(X_n)$ centralizes $\Delta^n(H)$. We will argue that $g=e$. Let us decompose $g$ into two coordinate functions $g_1 : X_n \to \{1, \dots , n\}$ and $g_2 : X_n \to \mathbb{N}$, where $g_1$ determines which ray a given point lands on and $g_2$ determines the position of that given point on the ray. Hence, 
\[
g(k,j) = (g_1(k,j), g_2(k,j))
\]
for $(k,j) \in X_n$. Since $g$ commutes with all of $\Delta^n(H)$, we have $g \Delta^n(h) = \Delta^n(h)g$ for for every $h \in H$ which means that 
\[
g\Delta^n(h)(k,j) = g(k,h(j)) = \bigg( g_1(k,h(j)) , g_2(k,h(j)) \bigg)
\]
equals 
\[
\Delta^n(h)g (k,j) = \Delta^n(h) (g_1(k,j), g_2(k,j)) = \bigg(g_1(k,j), h(g_2(k,j)) \bigg)
\]
every $(k,j) \in X_n$ and every $h \in H$. This gives us two independent coordinate equations 
\[
g_1(k,h(j)) = g_1(k,j)
\]
and
\[
g_2(k,h(j)) = h(g_2(k,j))
\]
for every $k=1, \dots , n$ and $j \in \mathbb{N}$ and every $h \in H$. First, we show the element $g$ must completely preserve rays. That is, there is no initial block where $g$ can be some wild permutation moving points across rays. 

The first coordinate equation says that the map $g_1(k, \cdot) : \mathbb{N} \to \{1, \dots , n\}$ is invariant under the action of $H$. This means that $g_1(k, \cdot )$ is constant on the orbits of $H$ acting on $\mathbb{N}$. Because $H$ is finite-index in $\Gamma$, every orbit of $H$ acting on $\mathbb{N}$ must be infinite. Now, since $g \in \sym_{ends}(X_n)$, Lemma \ref{lemma:threshold} says there exists $i_g \in \mathbb{N}$ such that $g_1(k,j) = k$ for all $j \ge i_g$. This is because $g$ eventually must map points on the $k$-th ray to points on the $k$-th ray (it eventually stabilizes the $k$-th ray). Now, let $j \in \mathbb{N}$. Then $j$ must belong to one of the infinite $H$-orbits, and necessarily this infinite $H$-orbit must intersect the infinite tail $\{i_g, i_g +1, i_g +2, \dots \}$ on which $g_1(k, \cdot)$ constantly evaluates to $k$. Hence, this show $g_1(k,j) = k$ not just for all $j \ge i_g$ but in fact for all $j \in \mathbb{N}$. This means that $g$ bijectively maps the entire $k$-ray to itself. Now, we use the trivial centralizer fact to argue that $g$ must act trivially on each ray and therefore be the identity element.

To this end, let $\alpha_k : \mathbb{N} \to \mathbb{N}$ be defined as $\alpha_k (j) = g_2(k,j)$, which is essentially the restriction of $g$ to the $k$-th ray. The second coordinate equation can be written as 

$$\alpha_k(h(j)) = h( \alpha_k(j))$$
for every $j \in \mathbb{N}$ and $h \in H$. Hence, $\alpha_k \in C_{\sym (\mathbb{N})}(H) = \{e\}$ and therefore $g(k,j) = j$ for every $k=1, \dots , n$ and $j \in \mathbb{N}$. Finally, we see that

$$g(k,j) = (g_1(k,j), g_2(k,j)) = (k,j)$$
which proves $g$ is the trivial permutation (identity), and therefore
\[
C_{\sym_{ends}(X_n)}(\Delta^n(H)) = \{e\}.
\]
This shows that $\Delta^n(\Gamma)$ has the finite-index centralizer property. In particular, this shows that every non-trivial element of $H_n(G)$ has infinitely many $\Delta^n(\Gamma)$-conjugates. Since $\Delta^n(\Gamma) \cong \Gamma$, the group $\Delta^n(\Gamma)$ has property (T), so Lemma \ref{lemma:relativelyICC} tells us that $H_n(G)$ is not inner amenable. 
\end{proof}

\looseness =-1 Note that the above argument actually applies to any subgroup of $\sym_{ends}(X_n)$ containing $\Delta^n(\Gamma)$. 

Given that these examples are non-inner amenable, one might wonder whether they yield prime or even solid factors. Regarding $L(H_n(G))$ and solidity, we can actually prove $L(H_n(G))$ is not solid in most cases. First, if $G$ is amenable, then $H_n(G)$ is amenable and hence $L(H_n(G))$ is the hyperfinite II$_1$ factor, which is trivially solid because amenability passes to subalgebras. 

Now suppose $G$ is non-amenable and $n \ge 2$. Note that we can embed $H_{n-1}(G)$ into $H_n(G)$ by sending $g \in H_{n-1}(G)$ to the element which just acts the rays $R_1, \dots , R_{n-1}$ in the usual Houghton-like way, and is the identity on $R_n$. In other words, we identify it with all subgroup of elements where $g_n = 1$ and $s_n = r_n = 0$. Since $G$ is non-amenable, it follows $H_{n-1}(G)$ is non-amenable and therefore $L(H_{n-1}(G))$ is a non-amenable II$_1$ factor subalgebra. Now, because elements in $H_{n-1}(G)$ have support disjoint from elements in $S_{\infty}(R_n)$, the finitely supported permutations on $R_n$, the two subgroups commute. This translates to commuting subalgebras and hence the inclusion
\[
L(H_{n-1}(G)) \subseteq L(S_{\infty}(R_n))' \cap L(H_n(G)), 
\]
which means the relative commutant of the diffuse subalgebra $L(S_{\infty}(R_n))$ is non-amenable. 

As for whether $L(G)$ yields a solid factor for our specific example coming from the shift-similar closure of $\Gamma \cong \psl_k(\mathbb{Z})$, we can only rule out solidity for $k \ge 4$ and that is because assuming $k \ge 4$ gives us ``enough room" to construct a subalgebra whose relative commutant is non-amenable. To see this, let $V$ be the two dimensional subspace of $\mathbb{Q}^k$ spanned by $e_1$ and $e_2$. Let $P(V)$ be its image under the canonical projection map from $\mathbb{Q}^k$ to $P^{k-1}(\mathbb{Q})$, and let $Y$ be the image of $P(V)$ under the fixed bijection from $P^{k-1}(\mathbb{Q})$ to $\mathbb{N}$. Let $\stab_{\psl_k(\mathbb{Z}}(P(V))$ be the pointwise stabilizer of $P(V)$. The pointwise stabilizer $\stab_{\psl_k(\mathbb{Z})}(P(V))$ consists of all matrices of the form
\[
\begin{bmatrix}
I_2 & 0\\
0 & A
\end{bmatrix}
\]
for $A \in \psl_{k-2}(\mathbb{Z})$.  Hence, the pointwise stabilizer is isomorphic to $\psl_{k-2}(\mathbb{Z})$. Since $k -2 \ge 2$, the pointwise stabilizer is non-amenable. Let $S_{\infty}(Y)$ be the finitely supported permutations on the infinite set $Y$. The elements of $S_{\infty}(Y)$ and the embedded copy of the pointwise stabilizer $\stab_{\psl_k(\mathbb{Z}}(P(V))$ in $G$ have disjoint support. This implies that we have the inclusion 
\[
L(\stab_{\Gamma}(P(V)) \subseteq L(S_{\infty}(Y))' \cap L(G),
\]
and hence the relative commutant of $L(S_{\infty}(Y))$ is non-amenable. Because $L(G)$ is a subfactor of $L(H_1(G))$ and solidity passes to subfactors, $L(H_1(G))$ cannot be solid either for $k \ge 4$. 

For either $L(G)$ or $L(H_1(G))$, we do not know how to deal with the $k=3$ case yet, but we view all of this as extremely interesting. Shift-similar groups are quite rich and complex, and shift-similar closures inherit much of the complexity and structure of the underlying group.

Although we can completely rule out solidity for $L(H_n(G))$ for $n \ge 2$, and solidity for $L(G)$ and $L(H_1(G))$ in some cases, we cannot prove or disprove primeness in general. We do not know whether these examples yield prime factors. However, we can definitively rule out one obvious obstruction to primeness that occurs at the group level---namely, a non-trivial direct product. Indeed, by Theorem \ref{thm:NonWreathProduct} we know that neither $G$ nor $H_n(G)$ can be written as a non-trivial direct product. This is because direct products fall under the class of generalized wreath products by taking the acting group to act trivially on a point. Hence, if these von Neumann algebras were not prime, the tensor product decomposition could not come from subgroups.

Finally, one very natural question arises from this discussion. Given the construction of a shift-similar group $G$ such that neither $G$ nor $H_n(G)$ are inner amenable uses property (T) groups,  are there property (T) subgroups of $\sym (\mathbb{N})$ whose shift-similar closures act jointly amenably with the canonical shift map on $\mathbb{N}$, so that $L(G)$ and $L(H_n(G))$ are property Gamma factors and hence $G$ and $H_n(G)$ are inner amenable? We address this 
questions in the following two sections through two completely different sources of examples. 

\section{Examples from Residually Finite Groups}\label{sec:ResidFiniteGrp}

In this section, we construct a particular embedding of any infinite residually finite group such that the shift-similar closure gives us examples satisfies the conditions of Theorem \ref{mainthm:Gamma} and Theorem \ref{mainthm:McDuff}. In particular, Theorem \ref{thm:ResidFiniteClosureGamma} applies to property (T) groups satisfying Theorem \ref{mainthm:NonInnerAmean}. Thus, given the same property (T) group, $\psl_k(\mathbb{Z})$, we have have one embedding where the shift-similar closure is inner amenable and one where the shift-similar closure is non-inner amenable. These constructions underscore and highlight how rich and diverse shift-similarity is. 

In order to prove Theorem \ref{thm:ResidFiniteClosureGamma}, 
we characterize all collections of finite sets $\{F_m\}_{m \ge 1}$ of $\mathbb N$ which are F\o lner sets for the shift map $\sigma$; that is, all collections of finite sets $\{F_m\}_{m \ge 1}$ with $|F_m| \to \infty$ as $m \to \infty$ for which 
\[
\lim_{m \to \infty} \frac{|\sigma (F_m) \triangle F_m|}{|F_m|} = 0.
\]
Given a finite subset $F$ of natural numbers, we can uniquely decompose it as a union of maximally connected components:
\[
F = I_1 \sqcup I_2 \sqcup \dots \sqcup I_k.
\]
A subset $I \subseteq \mathbb{N}$ is connected if $a,b \in I$ and $a \le c \le b$, then $c \in I$. In this setting, connected subsets are simply contiguous intervals of natural numbers, and viewing $\mathbb{N}$ as a graph with edges between $j$ an $j+1$, the ``maximal connected components" are simply the vertices of standard connected components of the subgraph induced by $F$. We can therefore write $F$ as 
\[
F = \bigsqcup_{i=1}^{k} [a_i,b_i]
\]
with each $[a_i,b_i]$ a maximal contiguous intervals and with $b_i +1 < a_j$ whenever $i < j$. It is not hard to see that 
\[
F \setminus \sigma (F) = \{a_1, \dots , a_k\}
\]
and therefore 
\[
|\sigma (F) \triangle F| = 2|F \setminus \sigma (F)| = 2k.
\]

Now, given a collection $\mathcal{F} = \{F_m\}_{m \ge 1}$ of finite subsets of $\mathbb{N}$ with $|F_m| \to \infty$ as $m \to \infty$, each $F_m$ can be uniquely decomposed as a disjoint union of $k_m \in \mathbb{N}$ maximal contiguous intervals. Hence, 
\[
\frac{|\sigma (F_m) \triangle F_m|}{|F_m|} = \frac{2 k_m}{|F_m|}
\]
which converges to $0$, and therefore forms a F\o lner sequence for $\sigma$, if and only if 
\[
\lim_{m \to \infty} \frac{k_m}{|F_m|} = 0,
\]
which says that the number of intervals in each $F_m$ cannot be too large compared to the size of $F_m$. So, for example, if $F_m$ is a single interval growing with $m$, such as the standard initial intervals $F_m = [1,m] = \{1, \dots , m \}$, these will form F\o lner sets for $\sigma$. Other examples include the following.
\begin{enumerate}
 \setlength{\itemsep}{1.5ex}  
\item $F_m = \cup_{r=1}^{m} [ r^2 \cdot 2^m + 1, r^2 \cdot 2^m + 2^r]$ with $k_m = m$ and $|F_m| = 2^{m+1}-2$ 
\item $F_m = \bigcup_{r=0}^{2^m-1} [2mr+1, 2mr+m]$ with $k_m = 2^m$ and $|F_m| = m \cdot 2^m$ 
\item $F_m = \bigcup_{r=0}^{m-1} [rm^2 +1 , r m^2 + \floor{\sqrt{m}}]$ so $k_m = m$ and $|F_m| = m \floor{\sqrt{m}}$
\end{enumerate}

However, if each $F_m$ is a union of a bunch of small connected components, then this condition will fail. Non-examples include
\begin{enumerate}
\setlength{\itemsep}{1.5ex} 
\item $F_m = \{1,3, \dots , 2m-1\}$ with $k_m = m$ and $|F_m| = m$
\item $F_m = \{1,4,7, \dots , 3m-2\}$ so $k_m = m$ and $|F_m| = m$
\item $F_m = \bigcup_{r=0}^{m-1} [3r+1,3r+2]$ with $k_m = m$ and $|F_m| = 2m$
\end{enumerate}

Next, given a collection $\mathcal{F} = \{F_m\}_{m \ge 1}$ of F\o lner sets for $\sigma$, we can build a shift-similar group from the collection $\mathcal{F}$. The entire group, and hence any countable subgroup of it, satisfies the joint amenability condition. First, we can modify the family so that sets satisfy $\min F_m \to \infty$ as $m \to \infty$. Now, define $G^{(1)}_{\mathcal{F}}$ as 
\[
G_{\mathcal{F}}^{(1)} := \{g \in \sym (\mathbb{N}) : \lim_{m \to \infty} \frac{|g(F_m) \triangle F_m|}{|F_m|} = 0 \}.
\]
Clearly this contains the identity. Given $g,h \in G^{(1)}_{\mathcal{F}}$, the triangle inequality tells us 
\[
\frac{|gh(F_m) \triangle F_m |}{|F_m|} \le \frac{|gh(F_m) \triangle g(F_m) |}{|F_m|} + \frac{|g(F_m) \triangle F_m|}{|F_m|} = \frac{|h(F_m) \triangle F_m|}{|F_m|} + \frac{|g(F_m) \triangle F_m|}{|F_m|} \to 0 
\] 
as $m \to \infty$. Similarly, it is closed under inverses. Finally, given $g \in G^{(1)}_{\mathcal{F}}$ and $j \in \mathbb{N}$, because $\min F_m \to \infty$ as $m \to \infty$, we can pick $m$ sufficiently large so that for every $i \in F_m$ we have $i \ge j$ and $g(i+1) > g(j)$ for all $i \in F_m$. Hence, for $i \in F_m$, we have
\[
\psi_j(g)(i) = s^{-1}_{g(j)} g s_j(i) = g(i+1)-1 = \sigma^{-1}g \sigma (i),
\]
so $\psi_j(g) = \sigma^{-1} g \sigma$ on $F_m$ for all sufficiently large $m \in \mathbb{N}$, which means 
\begin{align*}
\frac{|\psi_1(g)(F_m) \triangle F_m|}{|F_m|} &= \frac{|\sigma^{-1} g \sigma (F_m) \triangle F_m|}{|F_m|} \\
& \le \frac{|\sigma (F_m) \triangle F_m|}{|F_m|} + \frac{|g(F_m) \triangle F_m|}{|F_m|} + \frac{|F_m \triangle \sigma (F_m)|}{|F_m|}
\end{align*}
which converges $0$ as $m \to \infty$. Hence, $G^{(1)}_{\mathcal{F}}$ is a subgroup which is shift-similar. There is actually an even smaller   
shift-similar subgroup of $G^{(1)}_{\mathcal{F}}$ that we can associate to $\mathcal{F}$:
\[
\mathcal{G}_{\mathcal{F}}^{(2)}:=\{g \in \sym (\mathbb{N}) : \exists C_g, m_g \: \text{ s.t. } |g(F_m) \triangle F_m| \le C_g \quad\forall m \ge m_g\}
\]

Now, given a residually finite group $\Gamma$, we can construct a sequence of sets $\mathcal{F} = \{F_m\}_{m \ge 1}$ from the residual chain that will be F\o lner sets for $\sigma$, and we can embed $\Gamma$ into the above subgroup $G_{\mathcal{F}}^{(2)}$. In particular, we have the following theorem.

\begin{theorem}\label{thm:ResidFiniteClosureGamma}
    Every infinite residually finite group $\Gamma$ admits an embedding into $\sym (\mathbb N)$ such that the shift-similar closure together with $\sigma$ acts amenably on $\mathbb{N}$.
\end{theorem}
\begin{proof}

Let $\Gamma$ be an infinite residually finite group. Then there exists a residual chain, that is, a decreasing sequence 
\[
\Gamma > \Gamma_1 > \Gamma_2 > \dots 
\]
of finite-index normal subgroups of $\Gamma$ such that $\bigcap_{m \ge 1} \Gamma_m = \{e\}$ and $N_m := |\Gamma : \Gamma_m| \to \infty$ as $m \to \infty$. Define $X$ as 
\[
X := \bigsqcup_{m \ge 1} \Gamma / \Gamma_m = \bigsqcup_{m \ge 1} Q_m.
\]

The group $\Gamma$ acts naturally on $X$ by permutations coming from the left multiplication action of $\Gamma$ on each finite quotient $Q_m = \Gamma / \Gamma_m$. Because $\bigcap_{m \ge 1} \Gamma_m = \{e\}$, the action is faithful, and the action of $\Gamma$ on each $Q_m$ is confined to $Q_m$, meaning no element of $\Gamma$ moves a coset from one coset space to another. The action on of $\Gamma$ on $Q_m$ is transitive, so the orbits are $Q_m$. 

Now, define the interval $I_1 = \{1, \dots , N_1\}$, and for $m \ge 2$, define the intervals 
\[
I_m := \{\sum_{i=1}^{m-1} N_i +  1, \dots , \sum_{i=1}^{m} N_i \}.
\]
These are adjacent, disjoint intervals partitioning $\mathbb{N}$ with $|I_m| = N_m = |Q_m|$. This means the interval $I_m$ can be put in  bijection with $Q_m$. So, for each $m \in \mathbb{N}$, fix any bijection $\phi_m$ from $Q_m$ to $I_m$. Let $\phi : X \to \mathbb{N}$ be the bijection which is the extension of all of them, and let $\widetilde{\phi} : \sym (X) \to \sym (\mathbb{N})$ the induced isomorphism. Then it is easy to see that $\widetilde{\phi}(\Gamma) \le G_{\mathcal{F}}^{(2)}$, where $\mathcal{F} = \{I_m\}_{m \ge 1}$. 

Define $G$ as the shift-similar closure of $\widetilde{\phi}(\Gamma) \le G^{(2)}_{\mathcal{F}}$. Since $G^{(2)}_{\mathcal{F}}$ is shift-similar, we have $G \le G^{(2)}_{\mathcal{F}}$, and hence $G$ with $\sigma$ satisfies the joint amenability condition. 
\end{proof}

Applying Corollary \ref{cor:shiftsimilargamma}, we immediately obtain the following.

\begin{corollary}
Let $\Gamma$ be a residually finite group. Embed it into $\sym(\mathbb{N})$ as outlined above, and let $G$ denote its shift-similar closure. Let $G_1, \dots , G_n \le \sym(\mathbb{N})$ be shift-similar groups with $G_k = G$ for some $k=1, \dots , n$. Then 
\[
L(G), \quad L(H_n(G)) \quad \text{ and } L(H_n(G_1, \dots , G_n))
\]
are property Gamma factors for all $n \ge 1$. 
\end{corollary}

\begin{remark}
    Since $\psl_k(\mathbb{Z})$ is residually finite and has property (T), Theorem \ref{thm:ResidFiniteClosureGamma} provides an example of a property (T) group with a shift-similar closure giving rise to a factor with property Gamma. In contrast, Theorem \ref{thm:NonInnerAmean} gives an embedding of  $\psl_k(\mathbb{Z})$ whose shift-similar closure in non-inner amenable.  
\end{remark}

As for applying Theorem \ref{thm:McDuff} to get McDuff factors, this process will not generally yield McDuff factors unless one of the other groups in the collection $G_1, \dots G_n$ is amenable. This is because the shift-similar closure $G$ might not necessarily be amenable even if $\Gamma$ is an amenable residually finite group.

The same action of a residually finite group on the union of its coset spaces can also be used to construct a McDuff factor, though it does not come from a shift-similar group.

\begin{proposition}
Let $\Gamma$ be a residually finite group, viewed as acting on 
\[
X = \bigsqcup_{m \ge 1} \Gamma / \Gamma_m
\]
by permutations induced by left multiplication on each finite quotient $\Gamma/\Gamma_m$. If $G_0 := \langle  S_{\infty}(X), \Gamma \rangle$, then $L(G_0)$ has property Gamma. If, moreover, $\Gamma$ is non-abelian, then $L(G_0)$ is McDuff. 
\end{proposition}

\begin{proof}
Let $\Gamma$ be a non-abelian residually finite group, and define $G_0 := \langle S_{\infty}(X), \Gamma \rangle$. Since $\Gamma$ is non-abelian, there exists $x, y \in \Gamma$ such that $[x,y] \neq e$, which means $[x,y]$ cannot lie in
\[
\bigcap_{m \ge 1} \Gamma_m = \{e\}.
\]
Hence, there exists $m_0 \in \mathbb{N}$ such that $[x,y] \notin \Gamma_{m_0}$, and hence $[x,y] \notin \Gamma_m$ for every $m \ge m_0$ because these subgroups form a decreasing sequence. 

We can remove the finitely many $\Gamma_m$ to which $[x,y]$ may possibly belong, and then we can re-index so that $[x,y] \notin \Gamma_m$ for every $m \in \mathbb{N}$. This means that $[x,y]$ survives the quotient $Q_m$ for every $m \in \mathbb{N}$, so every quotient is non-abelian. 

Now, pick any $p_m \in Q_m$, and define $r_{m,p_m}$ to be right multiplication by $p_m \in Q_m$ on the cosets and extend it to be the identity outside of $Q_m$. By construction, $r_{m,p_m} \in S_{\infty}(X)$. We claim $(\lambda (r_{m,p_m}))_{m \ge 1}$ is a central sequence in $L(G_0)$. Given any $g \in S_{\infty}(X)$, we have that $g$ acts by identity on $Q_m$ for all sufficiently large $m$. Since $r_{m,p_m}$ is supported on $Q_m$, they have disjoint support and hence commute. Now, recall that we are viewing $\Gamma$ as a subgroup of $\sym (X)$ via left multiplication on the coset spaces. It is not hard to see that $r_{m,q}$ commutes with every element of $\Gamma$. Hence, the sequence $(\lambda (r_{m,p_m}))_{m \ge 1}$ asymptotically commutes with every generator of $L(G_0)$, so $(\lambda (r_{m,p_m}))_{m \ge 1}$ is a central sequence in $L(G_0)$. 

To get a central sequence that does not commute with $(\lambda (r_{m,p_m}))_{m \ge 1}$, pick any $q_m \in Q_m$ that does not commute with $p_m$, which is possible thanks to $Q_m$ being non-abelian. Hence, $(\lambda (r_{m,q_m}))_{m \ge 1}$ is another central sequence, and it is not hard to show $p_m$ and $q_m$ not commuting implies $r_{m,p_m}$ and $r_{m,q_m}$ do not commute. Hence, $L(G_0)$ is McDuff. 
\end{proof}

Note that every element of $\Gamma$ acting on $X$ has infinite support, so the intersection of $S_{\infty}(X)$ and $\Gamma$ is trivial. Since $S_{\infty}(X)$ is normal and $G_0$ is generated by $S_{\infty}(X)$ and $\Gamma$, we can write $G_0$ as the internal semidirect product $G_0 = S_{\infty}(X) \rtimes \Gamma$, and therefore
\begin{equation}
    L(G_0) = L(S_{\infty}(X) \rtimes \Gamma) \cong L(S_{\infty}(X)) \rtimes \Gamma
\end{equation}

Hence, we get a crossed product with the hyperfinite factor by the action of a property (T) group. This can be contrasted with property (T) extensions of abelian groups that give rise to property Gamma factors that are not McDuff, such as [Cite Kida's non-McDuff,Gamma example].

Now, let us make a quick observation on groups generated in part by $S_\infty$. Let $G \le \sym (\mathbb{N})$ be an infinite group. Since $G \le \langle S_{\infty}, G \rangle$, we have the following inclusion of their shift-similar closures 
\[
\overline{G} \le \overline{\langle S_{\infty}, G \rangle}.
\]
Because $G$ is infinite, $\overline{G}$ contains $S_{\infty}$ as well as $G$. Hence, $\langle S_{\infty}, G \rangle \le \overline{G}$ and therefore $\overline{\langle S_{\infty}, G \rangle} \le \overline{G}$, and consequently, we have
\[
\overline{\langle S_{\infty}, G \rangle} = \overline{G}.
\]
Applying this to $\widetilde{\phi}(G_0)$, we have
\[
\overline{\widetilde{\phi}(G_0)} = \overline{\widetilde{\phi}(\langle S_{\infty}(X), \Gamma \rangle) } = \overline{\langle \widetilde{\phi}(S_{\infty}(X)), \widetilde{\phi}(\Gamma) \rangle} = \overline{\langle S_{\infty}, \widetilde{\phi}(\Gamma ) \rangle}= \overline{\widetilde{\phi}(\Gamma)}
\]
So, while
\[
L(\widetilde{\phi}(G_0)) \cong L(G_0)
\]
is McDuff, it is unclear whether 
\[
L(\overline{\widetilde{\phi}(G_0)}) = L(\overline{\widetilde{\phi}(\Gamma)})
\] 
is McDuff, although $L(\overline{\widetilde{\phi}(\Gamma)})$, as well as  $L(H_n(\overline{\widetilde{\phi}(\Gamma)}))$, has property Gamma. 

\section{Wobbling Groups}\label{sec:WobblingGroups}

The wobbling group of $\mathbb{N}$ endowed with the euclidean metric $W(\mathbb N, d_{euc})$ is a concrete subgroup of $\sym ( \mathbb N)$ coming from coarse geometry. It acts jointly amenably on $\mathbb N$ together with $\sigma$, and the standard intervals $F_m := \{1, \dots , m\}$ can be used to witness the amenability of this joint action. Additionally, it is not difficult to show $W(\mathbb N, d_{euc})$ is shift-similar. The wobbling group on its own is uncountable, 
but taking the shift-similar closure of any countable subgroup $G$ of $W(\mathbb N, d_{euc})$ gives a somewhat concrete example of a shift-similar group $\overline{G}$ that can be used to produce property Gamma factors using Theorem \ref{mainthm:Gamma} or Theorem \ref{mainthm:McDuff}. A key feature of using a wobbling group is that we can guarantee $\overline{G}$ and $\sigma$ act amenably on $\mathbb{N}$ without actually computing the shift-similar closure. As the constructions in Section \ref{sec:Non-innerAmean} and Section \ref{sec:ResidFiniteGrp} highlight, whether or not the shift-similar closure of a group acts amenably is highly dependent on the embedding into $\sym(\mathbb N)$. 

We provide sufficient conditions for $W(X,d)$ to act amenably on $(X,d)$ and subsequently characterize when the wobbling group of a countable metric space has a shift-similar embedding into $\sym (\mathbb{N})$ via some enumeration $\phi : \mathbb{N} \to X$. In reality, none of these spaces are truly external to $\mathbb{N}$. Indeed, if $(X,d)$ is any countable metric space, and $\phi : \mathbb{N} \to X$ is an enumeration, then clearly $d_{\mathbb{N}}(n,m) = d(\phi(n), \phi (m))$ defines a metric on $\mathbb{N}$ and
\[
W(\mathbb{N},d_{\mathbb{N}}) = \phi^{-1} W(X,d) \phi.
\]
Thus, characterizing shift-similar embeddings reduces to characterizing those metrics $d$ on $\mathbb{N}$ for which $W(\mathbb{N},d)$ is shift-similar. As it turns out, $W(\mathbb{N},d)$ being shift-similar is equivalent to the coarse geometric property of $\sigma^{-1}$ being bornologous with respect to $d$ (Theorem \ref{thm:shiftbornologous}). In Section \ref{sec:NonShiftMetricExamples} we then show there are natural metrics on $\mathbb N$ where the resulting wobbling group is not shift-similar. 

The structure of the underlying metric space can severely restrict the sort of groups that embed into its wobbling groups. For example, in \cite{juschenko15}, Juschenko and de la Salle showed that if $(X,d)$ has bounded geometry and uniform subexponential growth, then $W(X,d)$ cannot contain an infinite property (T), which applies to $W(\mathbb N, d_{euc})$. However, a result by Tessera included in \cite{juschenko15}, shows there are wobbling groups containing property (T) subgroups. Motivated by this result, we then study wobbling groups of locally finite connected graphs. Using the bornologous characterization of shift-similarity of a wobbling group from Theorem \ref{thm:shiftbornologous}, we prove every connected graph $\mathcal{G} = (V,E)$ with one-end admits a bijection $\phi : \mathbb{N} \to V$ such that $\phi^{-1} W(\mathcal{G},d_{\mathcal{G}}) \phi$ is shift-similar (Theorem \ref{thm:OneEndedShiftSim}). We then prove every two-ended graph which is ``coarsely" symmetric has a shift-similar embedding of its wobbling group (Theorem \ref{thm:TwoEndsShiftSim}). Without the coarsely symmetric assumption, we can construct a two-ended graph whose wobbling group does not admit a shift-similar embedding via any enumeration $\phi : \mathbb{N} \to V$. In both the one-ended and two-ended setting, if the graph is amenable, then the wobbling group acts amenably on $\mathbb{N}$ with $\sigma$.

We also prove there exists a one-ended graph whose wobbling group contains a countable property (T) subgroup. Moreover, this wobbling group is shift-similar and acts jointly amenably on $\mathbb N$ with $\sigma$. In particular, this gives a completely different construction from Theorem \ref{thm:ResidFiniteClosureGamma} of a property (T) group whose shift-similar closure gives rise to a property Gamma factor.   

\subsection{Amenable Actions of Wobbling Groups}

We recall some basic ideas from coarse geometry and the definition of a wobbling group. Unless otherwise specified, we work exclusively with countable metric space, and often times the metric spaces will also be uniformly discrete and of bounded geometry--such as Cayley graphs of finitely generated groups or locally finite, connected graphs more generally. Recall a metric space is \emph{uniformly discrete} if there exist $\delta > 0$ such that $d(x,y) \ge \delta$ for every $x,y \in X$. A metric space $(X,d)$ has \emph{bounded geometry} if for every $r > 0$, there exists $M(r) > 0$ such that 
\[
\sup_{x \in X} |B(x,r)| \le M(r) < \infty.
\]
Given a set $F$ in a metric space and $ r > 0$, define the $r$-boundary of $F$ as  
\[
\partial_r F = \{x \in X \setminus F : d(x,F) \le r\}.
\]
Recall that that a metric space $(X,d)$ is amenable if there exists a sequence of finite sets $\{F_m\}_{m \ge 1}$ in $X$ such that for every $r > 0$, we have
\[
\lim_{m \to \infty} \frac{|\partial_r F_m|}{|F_m|} = 0.
\]

Given a metric space $(X,d)$, the wobbling group of $(X,d)$, denoted as $W(X,d)$, is defined as all permutations $g$ of $X$ of bounded displacement. Equivalently, it is all permutations of $X$ that are coarsely close to the identity. More precisely, 
\[
W(X,d) = \{g \in \sym (X) : \sup_{x \in X} d(g(x),x) < \infty \}
\]

For a group $G$, if $S_1$ and $S_2$ are finite symmetric generating sets for $G$, then 
\[
W(G,d_{S_1}) = W(G,d_{S_2})
\]
because any two word metrics are bi-Lipschitz. Given any enumeration $\phi : \mathbb{N} \to X$, it is not difficult to see that  $W(X,d)$ acts amenably on $X$ if and only if $\phi^{-1} W(X,d) \phi$ acts amenably on $\mathbb{N}$.

Our primary interest in uniformly discrete metric spaces of bounded geometry (UDBG spaces) stems from the following lemma, which gives us a direct link between amenability of a metric space and amenability of $W(X,d)$ acting on $X$.

\begin{lemma}\label{lemma:AmeanMetric/Wobbling}
Let $(X,d)$ be a metric space. If $(X,d)$ is an amenable metric space, then the action of the wobbling group $W(X,d)$ on $X$ is amenable. If, moreover, $(X,d)$ is a UDBG space, then $W(X,d)$ acting amenably on $X$ implies $(X,d)$ is amenable as a metric space. 
\end{lemma}
\begin{proof}
First, suppose $(X,d)$ is an amenable metric space. Then there exists a sequence of sets $\{F_m\}_{m\ge 1}$ such that for every $r > 0$ we have 
\[
\lim_{m \to \infty} \frac{|\partial_r F_m|}{|F_m|} = 0.
\]
Given $g \in W(X,d)$, there exists $r > 0$ such that 
\[
d(g(x),x) \le r
\]
for every $x \in X$. First, we show $g(F_m) \setminus F_m \subseteq \partial_r F_m$. Indeed, let $y \in g(F_m) \setminus F_m$. Then there exists $x \in F_m$ such that $y=g(x)$ and $y=g(x) \notin F_m$. Now, because $d(y,F_m)$ is the infimum, we have 
\[
d(y,F_m) \le d(y,x) = d(g(x),x) \le r,
\]
which proves the inclusion. Hence,
\[
\frac{|g(F_m) \triangle F_m|}{|F_m|} = \frac{2|g(F_m) \setminus F_m|}{|F_m|} \le 2 \frac{|\partial_r F_m|}{|F_m|},
\]
which converges to $0$ as $m \to \infty$, thereby proving $W(X,d)$ acts amenably. 

Now, assume $(X,d)$ is a UDBG space and that $W(X,d)$ acts amenably on $X$. By way of contradiction, suppose $(X,d)$ is a non-amenable metric space. Recall Whyte proved in \cite{whyte99} that $(X,d)$ being a non-amenable is equivalent to the existence an action of a free group $\mathbb{F}_2$ on $X$ which is``translation-like," meaning it acts freely and $\{d_{X}(g(x),x) : x \in X \}$ is a bounded subset of $\mathbb{R}$ (this characterization of Whyte's result comes from \cite{seward14}). These conditions are precisely the conditions ensuring that $\mathbb{F}_2$ embeds into $W(X,d)$ and acts freely on $X$. Hence, if $W(X,d)$ acts amenably on $X$, then $\mathbb{F}_2$ acts amenably on $X$. But if a group acts freely and amenably on a set, then the group itself must be amenable, which is absurd since $\mathbb{F}_2$ is non-amenable. Hence, $(X,d)$ must be amenable. 
\end{proof}

Recall that a finitely generated group $G$ is amenable if and only if its Cayley graph is amenable as a metric space. Therefore, we have the following corollary that is well known to experts and appears, for example, in \cite{@paradoxicaldecomp}.

\begin{corollary}
A finitely generated group $G$ is amenable if and only if its wobbling group acts amenably on $G$.
\end{corollary}

\subsection{Shift-Similar Wobbling Groups}

We now turn our attention to characterizing for which metrics $d$ on $\mathbb{N}$ for which $W(\mathbb N, d)$ is shift-similar. Shift-similarity turns out to be
completely determined by when the inverse of the canonical shift map on $\mathbb N$ is bornologous with respect to $d$. We recall the definition of bornologous and related notations that we will need for the rest of the paper. 

Let $(X,d_X)$ and $(Y,d_Y)$ be metric space. A map $f : X \to Y$ is \emph{bornologous} if for every $R > 0$, there exists $S > 0$ such that whenever $d_{X}(x,y) \le R$, we have $d_{Y}(f(x),f(y)) \le S$. For a map $f$ between connected graphs, we have that $f$ is bornologous if and only if it is Lipschitz. Two maps $f,g : X \to Y$ are \emph{close} if 
\[
\sup_{x \in X } d_Y(f(x),g(x)) < \infty.
\]

A map $f : X \to Y$ is called a \emph{coarse equivalence} if it is bornologous and there exists a bornologous map $g : Y \to X$ such that $f \circ g$ and $g \circ f$ are close to $id_{X}$ and $id_{Y}$, respectively. We call $f$ and $g$ coarse inverses. As a consequence, any coarse equivalence \(f\) is automatically coarsely proper, meaning the preimage of any bounded set is bounded. Indeed, let $B \subset Y$ be a bounded set with diameter $D > 0$. Let $x,y \in f^{-1}(B)$. Then $f(x),f(y) \in B$ and therefore $d_{Y}(f(x),f(y)) \le D$. Since $g$ is bornologous, there exists $S_{D} > 0$ such that $d_{X}(g(f(x)), g(f(y)) \le S_D$. Also, since $g \circ f$ is close to $id_{X}$, say uniformly within $C > 0$ units, we have
\[
d_{X}(x,y) \le d_{X}(x, gf(x)) + d_{X}(gf(x),gf(y)) + d_{X}(gf(y),y) \le C + S_D + C,
\]
which proves $f^{-1}(B)$ is bounded. Likewise, $g$ is coarsely proper. 

If $\phi : X \to Y$ is a bijective coarse equivalence of metric spaces, then $\phi$ induces an isomorphism between the wobbling groups $W(X,d_X)$ and $W(Y,d_Y)$ via $g \mapsto \phi \circ g \circ \phi^{-1}$.

We now characterizes when the wobbling group $W(\mathbb{N},d)$ is shift-similar.

\begin{theorem}\label{thm:shiftbornologous}
Let $d$ be a metric on $\mathbb{N}$. Then $W(\mathbb{N},d)$ is shift-similar if and only if $\sigma^{-1}$ is bornologous with respect to the metric $d$ on $\mathbb{N}$.
\end{theorem}
\begin{proof}
Because the wobbling group with respect to any metric contains the finitely supported permutations, it is enough to show $W(\mathbb{N},d)$ is closed under $\psi_1$ to demonstrate it is shift-similar. First, suppose $\sigma^{-1}$ is bornologous with respect to $d$. Let $g \in W(\mathbb{N},d)$. Then there exists $R > 0$ such that 
\[
d(g(x),x) \le R
\]
for every $x \in \mathbb{N}$. Because $\sigma^{-1}$ is bornologous, there exists $S > 0$ such that whenever $d(x,y) \le R$, we have $d(\sigma^{-1} (x),\sigma^{-1} (y)) \le S$. The wobbling assumption on $g$ implies, in particular, that $d(g(x+1),x+1) \le R$ for every $x \in \mathbb{N}$. Hence, the bornologous condition on $\sigma^{-1}$ implies 
\[
S \ge d(\sigma^{-1}(g(x+1)), \sigma^{-1}(x+1)) = d(g(x+1)-1,x)
\]
\looseness = - 1For all sufficiently large $x \in \mathbb{N}$, we have $\psi_1(g)(x) = g(x+1)-1$, so the above can be written as 
\[
d(\psi_1(g)(x),x) \le S
\]
for all sufficiently large $x \in \mathbb{N}$. This proves that $\psi_1(g) \in W(\mathbb{N},d)$. 

We prove the other implication contrapositively. Suppose $\sigma^{-1}$ is not bornologous with respect to $d$. Then there exists $R > 0$ such that for every $n \in \mathbb{N}$, we can find points $x_n,y_n \in \mathbb{N}$ satisfying $d(x_n,y_n) \le R$ but $d(\sigma^{-1}(x_n), \sigma^{-1}(y_n) > n$ or $d(x_n - 1, y_n - 1) > n$. In particular, the latter condition implies that $x_n$ and $y_n$ are distinct for every $n \in \mathbb{N}$. We construct a wobbling group element $g$, defined as an infinite product of 2-cycles, such that $\psi_1(g)$ is not a wobbling element. We will construct $g$ by pairing up terms from the sequences $(x_n)_{n \ge 1}$ and $(y_n)_{n \ge 1}$ and defining $g$ to swap them. We break this construction into two cases.

First, suppose no natural number appears infinitely many times in either $\{x_n\}_{n \ge 1}$ or $\{y_n\}_{n \ge 1}$. From these sequences, we extract pairwise disjoint pairs. To start, let $a_1 = x_1$ and $b_1 = y_1$. Because $a_1$ and $b_1$ appear only finitely many times in either sequence, we can find a large enough index $n_2$ such that $(a_2,b_2) = (x_{n_2},y_{n_2})$ is distinct from the numbers in the previously chosen pair. Inductively, we can construct an infinite sequence of pairwise disjoint pairs $\{(a_k,b_k)\}_{k \ge 1}$ where $a_k$ comes from $\{x_n\}_{n \ge 1}$ and $b_k$ comes from $\{y_n\}_{n \ge 1}$, and we have 
\[
d(a_k-1,b_k-1) > n_k.
\]
Now, define the permutation $g$ to transpose $a_k$ and $b_k$, and fix all other natural numbers. On the one hand, when $x=a_k$, we get
\[
d(g(x),x) = d(g(a_k),a_k) = d(b_k,a_k) = d(x_{n_k}, y_{n_k}) \le R.
\]
Taking $x=b_k$ is symmetric, and if $x \neq a_k, b_k$ for any $k \ge 1$, then $d(g(x),x)=0$. Hence, we have $d(g(x),x) \le R$ for every $x \in \mathbb{N}$ and therefore $g \in W(\mathbb{N},d)$. On the other hand, for all sufficiently large $k$, we have
\[
d(\psi_1(g)(b_k-1),b_k-1) = d(g(b_k)-1,b_k-1) = d(a_k-1,b_k-1) > n_k \to \infty
\]
as $k \to \infty$, so $\psi_1(g) \notin W(\mathbb{N},d)$, which proves $W(\mathbb{N},d)$ is not shift-similar. 

Now, suppose some natural number $c$ appears in the sequences infinitely often. By the pigeonhole principle, $c$ must appear in either $\{x_n\}_{n \ge 1}$ or $\{y_n\}_{n \ge 1}$ infinitely often. Without loss of generality, assume it appears in $\{x_n\}_{n \ge 1}$ infinitely often. Then there exists a subsequence $\{x_{n_k}\}_{k \ge 1}$ with $x_{n_k} = c$ for every $k \in \mathbb{N}$. Hence,
\[
d(c-1,y_{n_k}-1) = d(x_{n_k}-1,y_{n_k}-1) > n_k \to \infty
\]
Therefore, the subsequence $\{d(c-1,y_{n_k}-1) \}_{k \ge 1}$ diverges to infinity. Now, if $\{y_{n_k} \}_{k \ge 1}$ consisted of only finitely many values, then $d(c-1, y_{n_k}-1)$ could only take on finitely many values, and this would make it impossible for the sequence to diverge to infinity. Hence, $\{y_{n_k}\}_{k \ge 1}$ consists of infinitely many distinct values, and we can extract a subsequence denoted $\{z_k\}_{k \ge 1}$ such that every value is distinct. Because we have $R \ge d(x_{n_k},y_{n_k}) = d(c,y_{n_k})$ for every $k \in \mathbb{N}$, we have $d(c,z_k) \le R$ for every $k \in \mathbb{N}$.  As we did above, we construct mutually pairwise disjoint pairs, except in this case we only have one sequence to extract the pairs from, and then we define a permutation which swaps them. First, let $a_1 = z_1$. Because $d(c-1,z_k-1) \to \infty$, there exists $b_1 := z_{k_1}$ such that
\[
d(c-1,b_1-1) > d(c-1,a_1-1) + 1
\]
and necessarily $b_1$ is distinct from $a_1$. By the reverse triangle inequality,  
\[
d(a_1-1, b_1-1) \ge d(c-1,b_1-1) - d(c-1,a_1-1) > d(c-1,a_1-1) + 1 - d(c-1,a_1-1) > 1
\]
Now, assume we have picked pairwise disjoint pairs $(a_1,b_1), \dots , (a_{m-1},b_{m-1})$ with
\[
d(a_i-1,b_i-1) > i
\]
for all $i=1, \dots , m-1$. Let $a_m$ be anything in $\{z_k\}_{k \ge 1}$ that is distinct from every number in the pairs we have already selected, which is possible thanks to the number of pairs previously selected being finite. Again, because $d(c-1,z_k-1) \to \infty$, we can pick $b_m$ distinct from the pairs we have already chosen and such that
\[
d(c-1,b_m-1) > d(c-1,a_m-1) + m
\]
and again our choice forces $b_m$ to be distinct from $a_m$. Again, by the reverse triangle inequality, we get
\[
d(a_m-1,b_m-1) \ge d(c-1,b_m-1) - d(c-1,a_m-1) > m
\]
Hence, we have inductively constructed disjont sequences $\{a_m\}_{m \ge 1}$ and $\{b_m\}_{m \ge 1}$ such that all the pairs $(a_m,b_m)$ are pairwise disjoint and $d(a_m-1,b_m-1) > m$ for every $m \in \mathbb{N}$.

Now, define $g$ to transpose $a_m$ and $b_m$, and fix all other natural numbers. Then if $x = a_m$, we get
\[
d(g(a_m),a_m) = d(b_m,a_m) \le d(b_m,c) + d(c,a_m) \le 2 R
\]
By symmetry, $x=b_m$ gives us the same inequality, and if $x \neq a_m,b_m$, then $d(g(x),x) = 0$. Hence, $d(g(x),x) \le 2 R$ for every $x \in \mathbb{N}$ and therefore $g \in W(\mathbb{N},d)$. However,
\[
d(\psi_1(g)(a_m-1),a_m-1) = d(g(a_m)-1,a_m-1) = d(b_m-1,a_m-1) > m
\]
which demonstrates $\psi_1(g)$ is not in $W(\mathbb{N},d)$, and therefore $W(\mathbb{N},d)$ is not shift-similar. 
\end{proof}

The following corollary is a somewhat more ``intrinsic" way of determining whether the wobbling group $W(X,d)$ embeds into $\sym (\mathbb{N})$ via some enumeration that naturally follows from the above theorem. 

\begin{corollary}[Theorem \ref{mainthm:WobblingShiftSim}]\label{cor:WobblingShiftSim}
Let $(X,d)$ be a countable metric space and $\phi : \mathbb{N} \to X$ a bijection. Define $T : X \setminus \{\phi (1)\} \to X$ by 
\[
T(x) = \phi \sigma^{-1} \phi^{-1}(x) = \phi (\phi^{-1}(x)-1).
\]Then $W(X,d)$ has a shift-similar embedding via $\phi$ if and only if $T$ is bornologous with respect to $d$.  
\end{corollary}
\begin{proof}
It is not difficult to see that $T$ is bornologous with respect to $d$ if and only if $\sigma^{-1}$ is bornologous with respect to the metric 
\[
d_{\mathbb{N}}(n,m) = d(\phi (n), \phi (m))
\]
for $n,m \in \mathbb{N}$ which by the above theorem is equivalent to 
\[
\phi^{-1} W(X,d) \phi = W(\mathbb{N}, d_{\mathbb{N}}) 
\]
being shift-similar. 
\end{proof}

Now we provide a sufficient condition for an enumeration $\phi : \mathbb{N} \to X$ of a countable metric space $(X,d)$ to yield a shift-similar embedding of $W(X,d)$ into $\sym (\mathbb{N})$. First, we call a bijection $\phi : \mathbb{N} \to X$ from $\mathbb N$ to a countable metric space $(X,d)$ a \emph{bounded-step enumeration} if 
\[
\sup_{n \in \mathbb{N}} d(\phi (n), \phi (n+1)) < \infty.
\]
Bounded-step enumerations turn out to be very useful as we will see in the next lemma and the following sections.

\begin{lemma}\label{lemma:boundeStepBij}
Let $(X,d)$ be a countable metric space and $\phi : \mathbb{N} \to X$ a bounded-step enumeration.
Then $T : X \setminus \{\phi(1)\} \to X$ defined as $T(x) = \phi (\phi^{-1}(x)-1)$ is bornologous with respect to $d$. Equivalently, $\phi W(X,d)\phi^{-1}$ is shift-similar in $\sym (\mathbb{N})$. 
\end{lemma}
\begin{proof}
By Corollary \ref{cor:WobblingShiftSim}, proving $\phi W(X,d) \phi^{-1}$ is shift-similar is equivalent to proving $T : X \setminus \{\phi (1)\} \to X$ is bornologous, which is how we proceed. Let $R > 0$ and $x,y \in X \setminus \{\phi (1)\}$ satisfy $d(x,y) \le R$. Then there exist $n,m \in \mathbb{N} \setminus \{1\}$ such that $x = \phi (n)$ and $y = \phi (m)$. Hence,
\begin{align*}
d(T(x),T(y)) & = d(\phi (\phi^{-1}(x)-1), \phi (\phi^{-1}(y)-1)) \\
&= d(\phi (\phi^{-1} \phi (n) -1), \phi (\phi^{-1} \phi (m)-1)) \\
&= d(\phi (n-1),\phi (m-1) ) \\
&\le d(\phi(n), \phi (n-1)) + d(\phi (n),\phi (m)) + d(\phi (m),\phi (m-1)) \\
&\le 2 \sup_{n \in \mathbb{N} } d(\phi (n),\phi (n+1)) + R
\end{align*}
Thus, if we define $S := 2 \sup_{n \in \mathbb{N} } d(\phi (n),\phi (n+1)) + R$, we can conclude $T$ is bornologous.
\end{proof}

\subsection{Examples of Shift-Similar and Non Shift-Similar Wobbling Groups}\label{sec:NonShiftMetricExamples}

Examples of metrics on $\mathbb N$ satisfying Theorem \ref{thm:shiftbornologous} (and hence produce shift-similar wobbling groups) include translation-invariant metrics and ``coarsely" translation-invariant metrics. The standard Euclidean metric is translation-invariant and one can induce many more using non-decreasing functions. Indeed, if $d$ any additive metric and $f : \mathbb{N}_0 \to [0, \infty)$ any non-decreasing function satisfying $f(0) = 0$, $f(n) > 0$ for every $n > 0$, and $f(a+b) \le f(a) + f(b)$, then $\rho (x,y) := f(d(x,y))$ is also a translation-invariant metric. Coarsely translation-invariant metrics are metrics $d$ for which there exist $C > 0$ such that
\[
|d(x+1,x+1) -d(x,y)| \le C
\]
for every $x \in \mathbb{N}$. Examples of coarsely translation-invariant metrics are $d(x,y) = |\sqrt{x}-\sqrt{y}|$, $d(x,y) = |x^p - y^p|$ for $p \in (0,1)$, $d(x,y) = |\ln x - \ln y|$, and. As for a metric which is not coarsely translation-invariant but $\sigma^{-1}$ is still bornologous with respect to it, we have $d(x,y) = |x^2 - y^2|$. 

There are also metrics on $\mathbb N$ which fail the conditions of Theorem \ref{thm:shiftbornologous}, which we now explore.

\begin{example}[p-adic Valuation Wobbling Group]
For the metric
\[
d_p(x,y) = \sum_{p \text{ prime} } |v_p(x) - v_p(y)|,
\]
where $v_{p}(x)$ is the $p$-adic valuation of $x \in \mathbb{N}$ -- i.e., the exponent of $p$ in the unique prime factorization of $x$, the wobbling group $W(\mathbb N,d_p)$ is not shift-similar. 
\end{example}

\begin{proof}
We first construct two sequences of primes of the form $p_n = 2^n k_n+1$ and $q_n = 2 \ell_n +1$ for $k_n$ and $\ell_n$ odd, where $n \ge 2$. First, it is not hard to see the prime $p_n$ can be written as $p_n = 2^n k_n+1$ for some odd number $k_n$ if and only if $p_n \equiv 2^n + 1 \pmod{2^{n+1}}$. Likewise, the prime $q_n$ is of the form $q_n = 2 \ell_n +1$ for $\ell_n$ odd if and only if $q_n \equiv 3 \pmod 4$. Since $n \ge 2$, we also have that $p_n \equiv 1 \pmod 4$, and hence $p_n$ and $q_m$ cannot be equal for any $n \ge 2$ and $m \in \mathbb{N}$.  

Now, let us actually construct these primes. For every $n \ge 2$, since $\gcd(2^{n+1},2^n+1) = 1$, Dirichlet's theorem tells us that $\{2^{n+1} + (2^n+1) : t \in \mathbb{N} \}$ has infinitely many primes. Hence, for every $n \ge 2$, there exists $t_n \in \mathbb{N}$ such that 
\[
p_n := 2^{n+1} t_n + (2^n + 1) = 2^n (2t_n +1) + 1 = 2^n k_n + 1
\]
where $k_n := 2t_n+1$ is odd, and we can pick these to be distinct by constructing these inductively and removing the previously selected primes $p_1, \dots , p_{n-1}$ to pick $p_n$. Likewise, since $\gcd(3,4) = 1$, Dirichlet's theorem tells us that there are infinitely many primes in $\{4t_n + 3 : t \in \mathbb{N} \}$. Hence, for every $n \ge 2$, there exists $t_n \in \mathbb{N}$ such that
\[
q_n := 4t_n + 3 = 2(2t_n+1) +1 = 2 \ell_n + 1
\]
is prime, where $\ell_n := 2t_n + 1$ is odd. Again, we can inductively pick the $t_n$ and consequently the $\ell_n$ so that all of the $q_n$ are distinct. As we noted, since $p_n \equiv 1 \pmod 4$ and $q_n \equiv 3 \pmod 4$, $p_n$ and $q_m$ cannot be equal for any $n \ge 2$ and $m \in \mathbb{N}$. 

Now, define $g$ to swap $p_n$ and $q_n$ and fix all other natural numbers. This is a well-defined permutation. Moreover, it is clear that $d_p(g(x),x) \le 2$, so $g \in W(\mathbb{N},d_p)$. However, 
\begin{align*}
d_p(\psi_1(g)(p_n-1) , p_n-1) &= d_p(g(p_n)-1, p_n-1) \\
&= d_p(q_n-1,p_n-1) \\
&= d_p(2^n k_n, 2 \ell_m) \\
& \ge |v_2(2^n k_n) - v_2(2 \ell_n)| \\
&= n-1,
\end{align*}
which proves $\psi_1(g)$ is not in $W(\mathbb{N},d_p)$. 

Note that this also clearly shows $\sigma^{-1}$ is not bornologous with respect to this prime metric. Indeed, for every $n \ge 2$, we constructed primes $p_n$ and $q_n$ with $d_p(p_n,q_n) = 2$ but $d_p(p_n-1,q_n-1) \ge n-1$ or $d_p(\sigma^{-1}(p_n),\sigma^{-1}(q_n)) \ge n-1$. Moreover, the permutation $g$ we constructed is entirely supported on the primes. 
\end{proof}

\begin{observation}
    An interesting observation from the previous example is that any bijection entirely supported on the primes is in $W(\mathbb{N},d_p)$ because 
\[
d_p(p,q) \le 2
\]
for any primes $p$ and $q$, meaning $\sym(P)$ embeds into $W(\mathbb{N},d_p)$. Hence, $W(\mathbb{N},d_p)$ certainly contains property (T) groups (every property (T) group in fact), but unfortunately it is not shift-similar, nor does it act amenably on $\mathbb{N}$.
\end{observation} 

For a geometric example illustrating the failure of the conditions of Theorem \ref{thm:shiftbornologous}, 
we can put the even numbers into a circle, the odd numbers on the positive $x$-axis, and then use the usual Euclidean metric $d_{\mathbb{R}^2}$ on $\mathbb{R}^2$. More precisely, we have the following.

\begin{example}
    Let $\phi : \mathbb{N} \to \mathbb{R}^2$ be defined by $\phi(2k) = \frac{1}{2} (\cos (k), \sin (k))$ and $\phi(2k-1) = (k^2,0)$. Define $d(x,y) = d_{\mathbb{R}^2}(\phi(x),\phi (y))$. Then $W(\mathbb{N},d)$ is not shift-similar and $\sigma^{-1}$ is not bornologous with respect to this metric.
\end{example}
\begin{proof}
    Showing $W(\mathbb{N},d)$ is not shift-similar and that $\sigma^{-1}$ is not bornologous with respect to this metric is similar to the above example with the prime metric. We have $\sym (2 \mathbb{N})$ embeds into  $W(\mathbb{N},d)$ because $d(x,y) \le 1$ for all even numbers $x,y$. And if $g \in \sym (2 \mathbb{N}) \setminus S_{\infty}(2 \mathbb{N})$, then $\psi_1(g) \notin W(\mathbb{N},d)$. 
\end{proof}

\subsection{Ends of a Graphs}

We now restrict our class of countable metric spaces further to those coming from the ends of a graph that is not necessarily the Cayley graph of a group. Before characterizing when the wobbling group of a locally finite connected graph with one end and with two ends are shift-similar (Section \ref{sec:OneEnd} and \ref{Sec:TwoEnds}) we set our notation and recall the relevant definitions. 

Given a coarse equivalence between graphs, it induces a bijection between the ends of the graphs. Let us recall how this is done.  Let $\mathcal{G} = (V,E)$ be a locally finite, connected graph. Given a set of vertices $V_0 \subseteq V$, write $\mathcal{G}[V_0]$ for the induced subgraph whose vertices are $V_0$ and whose edges are $\{u,v\} \in E$ for $u,v \in V_0$. For convenience, define $\mathcal{G} \setminus V_0 = \mathcal{G}[V \setminus V_0]$. It is the natural subgraph formed by removing the subgraph spanned by the vertices $V_0$.  

Because we are working with coarse geometric notions, we opt to use ``proper coarse rays" to define the ends rather than traditional proper rays, as it will simplify some ideas. A \emph{coarse ray} is a sequence of vertices $(v_i)_{i \ge 1}$ in $\mathcal{G}$ such that $\sup_{i \in \mathbb{N}} d(v_i,v_{i+1}) < \infty$, and it is said to be \emph{proper} if for every finite set $F$ of vertices, we have $v_i \notin F$ for all sufficiently large $i \in \mathbb{N}$. Note that the vertices $v_i$ and $v_{i+1}$ are not necessarily adjacent, but the lengths of the shortest paths from $v_i$ to $v_{i+1}$ are uniformly bounded across all $i \in \mathbb{N}$ since $\sup_{i 
\in \mathbb{N}} d(v_i,v_{i+1}) < \infty$. 

For $k \in \mathbb{N}$, the $k$-th tail of the coarse ray is the sequence $(v_k,v_{k+1}, \dots )$. For each $i \in \mathbb{N}$, choose the shortest path from $v_i$ to $v_{i+1}$. We regard a coarse ray and any tail as the connected subgraph (subspace) formed by the concatenation of these shortest paths whose lengths are uniformly bounded. If $F$ is a finite set of vertices, then we can find a large enough $k \in \mathbb{N}$ such that no vertex in the tail $(v_k, v_{k+1}, \dots )$ belong to $F$. Hence, the tail lives entirely in the subgraph $\mathcal{G} \setminus F$. Since $(v_k, v_{k+1}, \dots )$ corresponds to a connected subgraph (subspace), it lives in exactly one connected component of the subgraph $\mathcal{G} \setminus F$. Given a coarse ray $(v_1, v_2, \dots )$, and finite set $F$ of vertices, we say it eventually lies entirely in a connected component of $\mathcal{G} \setminus F$ if an entire tail of the coarse ray is a subset of that connected component.

Two proper coarse rays $r_1$ and $r_2$ are \emph{equivalent} (i.e., point in the same direction or define the same end) if for every finite set $F$ of vertices, there exist a sufficiently large $k \in \mathbb{N}$ such that the $k$-th tails of $r_1$ and $r_2$ live inside the same connected component of $\mathcal{G} \setminus F$. This gives us an equivalence relation on the set of all proper coarse rays in $\mathcal{G}$, and we call these equivalence classes the ends of the graph $\mathcal{G}$, denoting the set as $Ends(\mathcal{G})$. 

Note that every equivalence class has a genuine proper ray representing it: by concatenating the shortest paths between successive vertices, we obtain an infinite, connected, proper subgraph from which we can easily extract a bonafide proper ray by pruning cycles and any backtracking. Thus, using proper coarse rays to define ends is equivalent to the usual way.

Now let us show a coarse equivalence between two graphs induces a bijection between their ends.
\begin{lemma}
    A coarse equivalence $f$ between the graphs $\mathcal{G}_1 = (V_1, E_1)$ and $\mathcal{G}_2 = (V_2,E_2)$ induces a bijection between $Ends(\mathcal{G}_1)$ and $Ends(\mathcal{G}_2)$.
\end{lemma}
\begin{proof}
Assuming the setting of the lemma, let $(v_i)_{i \ge 1}$ be a proper coarse ray in $\mathcal{G}_1$. This means that there exists $R > 0$ such that $d(v_i,v_{i+1}) \le R$ for every $i \in \mathbb{N}$. Since $f$ is bornologous, there exist $S > 0$ such that $d(f(v_i),f(v_{i+1}) \le S$ for every $i \in \mathbb{N}$, which proves $(f(v_i))_{i \ge 1}$ is a coarse ray. Let $F$ be a finite set of vertices in $\mathcal{G}_2$. Since $f$ is proper, $f^{-1}(F)$ is bounded and hence finite, from which we get $v_i \notin f^{-1}(F)$ or $f(v_i) \notin F$ for all sufficiently large $i \in \mathbb{N}$. Hence, $(f(v_i))_{i \ge 1}$ is a proper coarse ray.

We define the induced map $f_* : Ends (\mathcal{G}_1) \to Ends (\mathcal{G}_2)$ as $f_*([(v_i)_{i \ge 1}]) = [(f(v_i))_{i \ge 1}]$. Let us check that this is well-defined. Let $(u_i)_{i \ge 1}$ and $(v_i)_{i \ge 1}$ be proper coarse rays in $\mathcal{G}_1$ which are equivalent. Since $f$ is bornologous, we can find $S > 0$ such that whenever $u$ and $v$ are adjacent in $\mathcal{G}_1$, we have $d_2(f(u),f(v)) \le S$. Let $F$ be a finite set of vertices in $\mathcal{G}_2$ and consider the $S$-neighborhood or $S$-thickening of $F$:
\[
N_S(F) := \{v \in V_2 : d(v,F) \le S \}
\]
Since the graph is locally finite, this set is finite and hence $f^{-1}(N_S(F))$ is finite. Therefore, there exists $k \in \mathbb{N}$ such that the tails $(u_i)_{i \ge k}$ and $(v_i)_{i \ge k}$ belong to the same component of $\mathcal{G}_1 \setminus f^{-1}(N_S(F))$. 

Let $n,m \ge k$. Since $u_n$ and $v_m$ lie in the same infinite component of $\mathcal{G}_1 \setminus f^{-1}(N_S(F))$, there exists a path $(p_1, \dots , p_\ell)$ connecting $u_n$ and $v_m$ that lies in that component of $\mathcal{G}_1 \setminus f^{-1}(N_S(F))$. In particular, for each $i=1, \dots , \ell$, we have $p_i \notin f^{-1}(N_S(F))$ or $f(p_i) \notin N_S(F)$ which means $d_2(f(p_i),F) > S$. Since $p_i$ and $p_{i+1}$ are adjacent, we have $d_2(f(p_i),f(p_{i+1})) \le S$. We can connect $f(p_i)$ and $f(p_{i+1})$ by a shortest path whose length is no more than $S$. Now, if any vertex on this path belonged to $F$, then we would have $d(f(p_i),F) \le S$, which is a contradiction. Hence, this path lies entirely in $\mathcal{G}_2 \setminus F$, and hence $f(p_i)$ and $f(p_{i+1})$ belong to the same connected component of $\mathcal{G}_2 \setminus F$. Concatenating these paths, we get a path from $f(p_1) = f(u_n)$ and $f(p_\ell) = f(v_m)$, implying that $f(u_n)$ and $f(v_m)$ belong to the same infinite component. These were arbitrary points in the tails $(f(u_i))_{i \ge k}$ and $(f(v_i))_{i \ge k}$, so the tails must lie in the same connected component of $\mathcal{G}_2 \setminus F$ and therefore define the same end.

Now, let $g$ be the coarse inverse of $f$. By the same argument as above, $g$ induces a well-defined map $g_* : Ends(\mathcal{G}_2) \to Ends(\mathcal{G}_1)$ on the ends. Because $g$ is a coarse inverse, there exists $C > 0$ such that $d(gf(v),v) \le C$ for every $v \in V(\mathcal{G}_1)$. Let $(v_i)_{i \ge 1}$ be a proper coarse ray in $\mathcal{G}_1$. Then there exists $R > 0$ such that $d(v_i,v_{i+1}) \le R$ for every $i \in \mathbb{N}$, and because $g \circ f$ is bornologous , there exists $S > 0$ such that $d(gf(v_i),gf(v_{i+1})) \le S$ for every $i \in \mathbb{N}$. 

Let $M := \max \{R,S,C\}$, and let $F \subset V(\mathcal{G}_1)$ be finite. Because $(v_i)_{i \ge 1}$ and $(gf(v_i))_{i \ge 1}$ are proper coarse rays, there exists $k \in \mathbb{N}$ such that  $v_i$ and $gf(v_i)$ do not belong to $N_M(F)$ for all $i \ge k$, and hence $d(v_i,F) > M$ and $d(gf(v_i),F) > M$. Let $i \ge k$. Pick the shortest path from $v_i$ to $v_{i+1}$. Its length is at most $R \le M$, and if the path met $F$, we would have $d(v_i,F) \le R \le M$, a contradiction. The path lies entirely in $\mathcal{G}_1 \setminus F$ and therefore the tail $(v_i)_{i \ge k}$ is connected inside. Likewise, a shortest path from $gf(v_i)$ to $gf(v_{i+1})$ has length at most $S \le M$, and if it met $F$, we would have $d(gf(v_i),F) \le S \le M$, a contradiction. Thus, the tail $(gf(v_i))_{i \ge k}$ is also connected in $\mathcal{G}_1 \setminus F$. Finally, a shortest path from $v_i$ to $gf(v_i)$ has length at most $C \le M$, so if it met $F$, we would have $d(v_i,F) \le C \le M$, a contradiction. 

Now, let $n,m \ge k$. Connect $v_n$ to $v_k$ in the connected component of $\mathcal{G}_1 \setminus F$ that $(v_i)_{i \ge k}$, connect $v_k$ to $gf(v_k)$ by a shortest path, and then connect $gf(v_k)$ t o $gf(v_m)$ inside the component of $\mathcal{G}_1 \setminus F$ that $(gf(v_i))_{i \ge k}$ lives. This proves that the two tails live inside the same component of $\mathcal{G}_1 \setminus F$, and consequently $[(gf(v_i))_{i \ge 1}] = [(v_i)_{i \ge 1}]$ or $g_* \circ f_* = id_{Ends(\mathcal{G}_1)}$. By an identical argument, we can show $f_* \circ g_* = id_{Ends(\mathcal{G}_2)}$.

\end{proof}

\subsection{Wobbling Groups of One-ended Graphs}\label{sec:OneEnd}

If a locally finite, connected graph $\mathcal{G}$ admits a bounded-step enumeration of its vertices, then this places an interesting restriction on the topological structure of the graph. Indeed, such a graph necessarily has one end. In fact, a locally finite, connected graph having one end is equivalent to having a bounded-step enumeration of its vertices (Theorem \ref{thm:OneEndedShiftSim}). 

Before we prove this, let us recall some important facts and terminology in graph theory. First, a Hamiltonian path in a graph is a path that visits each vertex exactly once, and a Hamiltonian cycle is a Hamiltonian path that ends at a neighbor of its starting vertex. Next, given a finite connected graph $\mathcal{G}$, the \emph{cube} of a graph, denoted $\mathcal{G}^3$, is the graph formed from $\mathcal{G}$ by adding (possibly new) edges in $\mathcal{G}$ between any two vertices in $\mathcal{G}$ if the distance between is at most $3$. In \cite{sekanina60}, Sekanina proved that the cube $\mathcal{G}^3$ of any finite connected graph $\mathcal{G}$ with at least three vertices has a Hamiltonian cycle starting at any vertex. With this in mind, we can prove the above mentioned characterization of locally finite one-ended graphs, and consequently we prove the wobbling group of any one-ended graph has a shift-similar embedding via a bounded-step enumeration of its vertices.

\begin{theorem}\label{thm:OneEndedShiftSim}
Let $\mathcal{G} = (V,E)$ be a locally finite, connected graph equipped with its canonical path metric $d_\mathcal{G}$. Then $\mathcal{G}$ has one-end if and only if it admits a bounded-step enumeration $\phi : \mathbb{N} \to V$ such that 
\[
\sup_{n \in \mathbb{N}} d_{\mathcal{G}}(\phi (n), \phi (n+1)) < 3.
\]
This implies $T$ is bornologous which is equivalent to  $\phi^{-1} W(\mathcal{G}) \phi$ in $\sym (\mathbb{N})$ being shift-similar. If in addition $\mathcal{G} = (V,E)$ is an amenable graph, then $\phi^{-1}W(\mathcal{G}) \phi$ together with $\sigma$ acts amenably on $\mathbb{N}$.  
\end{theorem}
\begin{proof}
By way of contradiction, suppose a bounded-step enumeration $\phi : \mathbb{N} \to V$ exists but $\mathcal{G}$ is a graph with more than one end. This means that there exists a finite set $F$ of vertices such that $\mathcal{G} \setminus F$ has at least two infinite components. Define 
\[
r := \sup_{n \in \mathbb{N}} d_{\mathcal{G}}(\phi (n),\phi (n+1)) \in (0 ,\infty ).
\]
Let $D$ be the maximum distance from any vertex in $F$ to any vertex in the finitely many finite components of $\mathcal{G} \setminus F$. We can thicken $F$ by considering the $\delta$-neighborhood $N_{\delta}(F)$, where $\delta > \max (r/2,D)$. This neighborhood will absorb all the finitely many finite components of $\mathcal{G} \setminus F$, and it allows us to remove a finite number of vertices even deeper within the infinite components created by removing $F$. 

Since $N_{\delta}(F)$ is finite, we can enumerate it in finitely many steps, meaning there exists $n_0 \in \mathbb{N}$ such that $\phi (n) \notin N_{\delta}(F)$ and in particular $\phi (n) \notin F$ for every $n \ge n_0$. Because removing $F$ from $\mathcal{G}$ creates at least two infinite components, we have that $\phi (n_0)$ must enumerate a vertex in one of the infinite components of $\mathcal{G} \setminus F$, call it $A$. But we have already enumerated $N_{\delta}(F)$, so we are specifically enumerating the vertices in $A \setminus N_{\delta}(F)$. Hence, $\phi(n_0) \in A \setminus N_{\delta}(F)$. 

At some point, $\phi$ must begin to enumerate the vertices in another part of the graph. Because $\phi$ has already enumerated the vertices in $N_{\delta}(F)$, the only other place to start enumerating is another infinite component of $\mathcal{G} \setminus F$. More precisely, let $B$ be another infinite component of $\mathcal{G} \setminus F$ distinct from $A$. Because $B$ is infinite, there exists a vertex $v \in B$ not in $\{\phi(1), \dots, \phi (n_0)\}$. By surjectivity, there exists $n_1 > n_0$ such that $\phi (n_1) = v$.  

Since $\phi(n_0) \in A \setminus N_{\delta}(F)$ and $\phi(n_1) \in B \setminus N_{\delta}(F)$, there is some largest $m \in \{n_0, \dots , n_1-1\}$ such that $\phi (m) \in A \setminus N_{\delta}(F)$. By maximality, $\phi (m+1) \notin A$. This means that $\phi (m+1)$ must be in an infinite component of $\mathcal{G} \setminus F$ different from $A$ and also not in $N_{\delta}(F)$ because $\phi (m+1) \notin \{\phi (1), \dots, \phi (n_0)\} \supset N_{\delta}(F)$. Because $\phi(m) \in A \setminus N_{\delta}(F)$ and $\phi (m+1) \notin N_{\delta}(F)$ are in distinct infinite components of $\mathcal{G} \setminus F$, any shortest path between $\phi(m)$ and $\phi (m+1)$ in the original graph $\mathcal{G}$ must pass through $F$. Let $v \in F$ be a vertex on a shortest path between them. Hence, 
\[
d(\phi (m),\phi (m+1)) = d(\phi (m),v) + d(v, \phi (m+1)) \ge d(\phi (m),F) + d(\phi (m+1),F) > 2 \delta > r. 
\]
which is a contradiction.

Now, assume $\mathcal{G}$ is a graph with one-end. Since the graph is locally finite and hence countable, a classical result of Halin in \cite{halin64} tells us that there exists an end-faithful spanning tree $\mathcal{T}$ of $\mathcal{G}$ with metric $d_{\mathcal{T}}$, that is, a spanning tree with one end. Because the spanning tree $\mathcal{T}$ is constructed from the same vertex set but possibly having some edges removed, a path between vertices $v$ and $w$ as viewed in $\mathcal{T}$ is possibly longer than the path in $\mathcal{G}$, which is to say $d_{\mathcal{G}}(v,w) \le d_{\mathcal{T}}(v,w)$ for all $v,w \in V$. Thus, it suffices to define a bounded-step enumeration $\phi$ with respect to $d_{\mathcal{T}}$. 

Since $\mathcal{T}$ is a one-ended tree, it contains a coarse proper ray $(v_i)_{i \ge 1}$ defining this end. Without loss of generality, we may assume that $(v_i)_{i \ge 1}$ is a bona fide proper ray because the equivalence class of any coarse proper ray has a proper ray representing it. Since $\mathcal{T}$ has exactly one end, every neighbor of $v_i$ other than $v_{i-1}$ and $v_{i+1}$ must lie in a finite connected component of $\mathcal{T} \setminus \{v_i\}$. If this were not so, some neighbor $u$ of $v_i$ (where $u \neq v_{i-1}, v_{i+1}$) would lie in an infinite component of $\mathcal{T} \setminus \{v_i\}$. Since this component is infinite and $\mathcal{T}$ is a tree, it contains a proper ray $(u_i)_{i \ge 1}$ starting at $u$ that never visits any vertex on the ray $(v_i)_{i \ge 1}$. That is because if some $u_i$ with $i > 1$ were equal to some $v_j$ for $j \in \mathbb{N}$, this would create a cycle in $\mathcal{T}$.

We claim these are non-equivalent proper rays. Let $k \in \mathbb{N}$ and suppose the $k$-th tails of $(u_i)_{i \ge 1}$ and $(v_i)_{i \ge 1}$ were in the same infinite component of $\mathcal{T} \setminus \{v_i\}$. In particular, $u_k$ and $v_k$ could be connected by a path. Since $v_i$ is removed, this path cannot pass through $v_i$. But, again, this would create a cycle in the original tree $\mathcal{T}$ with $v_i$ added back in it, a contradiction. Hence, the $k$-th tails of $(u_i)_{i \ge 1}$ and $(v_i)_{i \ge 1}$ lie in different infinite components for every $k \in \mathbb{N}$, which is strictly stronger than non-equivalence. This means these rays are not equivalent and therefore define distinct ends, contradicting the one-endness of $\mathcal{T}$. 

If $v_i$ has a neighbor that is neither $v_{i-1}$ nor $v_{i+1}$, then by the above we know it belongs to a finite component. This finite component is necessarily a tree we denote $\mathcal{T}_i$. Otherwise, $\mathcal{T}_i$ is just the single vertex $v_i$. Moreover, if $\mathcal{T}_i$ and $\mathcal{T}_j$ shared a vertex, where $i \neq j$, then this would create a cycle. Consequently, the finite trees $\mathcal{T}_i$ rooted at $v_i$ are disjoint.

We define the enumeration $\phi : \mathbb{N} \to V$ inductively on each finite subtree $\mathcal{T}_i$. First, define $\phi (1) = v_1$. Now, we begin enumerating the vertices in $\mathcal{T}_1$. If $\mathcal{T}_1$ consists only of $v_1$, there is nothing to do. If $\mathcal{T}_1$ has an additional vertex, then define $\phi(2)$ to be this additional vertex, so that $d_\mathcal{T}(\phi(1),\phi(2)) =1$, and then define $\phi(3) = v_2$, so that $d_\mathcal{T}(\phi(2),\phi(3)) = 2$. However, if $\mathcal{T}_1$ has at least three vertices, then we appeal to Sekanina's cube theorem. Indeed, by Sekanina's cube theorem, there exists a Hamiltonian cycle $(w_1, \dots ,w_k)$ in $\mathcal{T}_1^3$, where $w_1 = v_1$, $w_k$ is a neighbor of $w_1$ in $\mathcal{T}$, and every vertex of $\mathcal{T}_1$ appears exactly once in the cycle. We emphasize that in the cube $\mathcal{T}_1^3$ the distance between $w_i$ and $w_{i+1}$ is exactly one, but by Sekanina's theorem we know that $d_{T}(w_i,w_{i+1}) \le 3$ for all $i =1, \dots , k-1$ and $d_\mathcal{T}(v_1,w_k) = d(w_1,w_k) \le 3$, which is good enough for our purposes. Now, define $\phi (2) = w_2$, $\phi (3) = w_3$, etc. up to $\phi (k+1) = w_k$, and define $\phi (k+2) = v_2$. Then clearly we have $d_T(\phi(n),\phi(n+1)) \le 3$ for all $n=1, \dots , k-1$, and 
\[
d_\mathcal{T}(\phi(k+1), \phi (k+2)) = d_\mathcal{T}(w_k,v_2) \le d_\mathcal{T} (w_k,v_1) + d_T(v_1,v_2) = 2
\]
Hence, in all cases, we can enumerate the vertices in $\mathcal{T}_1$ so that $d_{T}(\phi(n),\phi (n+1)) \le 3$ for all $n =1, \dots , k+1$. Finally, we repeat the process at $v_2$ and enumerate the vertices of $\mathcal{T}_2$, and induction allows us to extend $\phi$ to enumerate all vertices in $\mathcal{G}$ with 
\[
\sup_{n \in \mathbb{N}} d_{\mathcal{G}}(\phi(n),\phi(n+1)) \le 3.
\]
Hence, $\phi$ is a bounded-step enumeration of the vertices. By Lemma \ref{lemma:boundeStepBij}, $\phi^{-1} W(\mathcal{G}) \phi$ is shift-similar.

Now, assume that $\mathcal{G}$ is amenable. By Lemma \ref{lemma:AmeanMetric/Wobbling}, this is equivalent to $W(\mathcal{G}, d_{\mathcal{G}})$ acting amenably on $\mathcal{G}$. Hence, there exists a collection of finite subsets $\{F_m\}_{m \in \mathbb{N}}$ of $V$ such that 
\[
\frac{|g(F_m) \triangle F_m|}{|F_m|} \to 0
\]
for every $g \in W(\mathcal{G}, d_{\mathcal{G}})$. We can remove $\phi (1)$ from every set $F_m$ without affecting the F\o lner condition. We claim that $E_m := \phi^{-1}(F_m)$ form a sequence of joint F\o lner sets for $\phi^{-1}W(X) \phi$ and $\sigma$ acting on $\mathbb{N}$, and we only need to check $\sigma$. Define $T : V\setminus \{\phi (1)\} \to V$ by 
\[
T(v) = \phi \sigma^{-1}\phi^{-1} (v) = \phi(\phi^{-1}(v)-1).
\]
It not hard to check that $T$ has bounded displacement. Indeed, if $v = \phi (n)$ with $n \ge 2$, then 
\[
d_{\mathcal{G}}(T(v),v) = d_{\mathcal{G}}(T(\phi (n)), \phi(n)) = d_{\mathcal{G}}(\phi (n-1),\phi (n)) \le 3
\]
So, if $T$ were a bijection of the vertices, then $T$ would belong to $W(\mathcal{G},d_{\mathcal{G}})$, and we could directly use it to show $\{E_m\}_{m \ge 1}$ are F\o lner sets for $\sigma$. Unfortunately, $T$ is not a bijection, but $T$ can be ``approximated" by two wobbling group elements in $W(\mathcal{G}, d_{\mathcal{G}})$. To see this, we define $g_1 : V \to V$ to swap $\phi (2k-1)$ and $\phi (2k)$ for $k \in \mathbb{N}$, and define $g_2 : V \to V$ by $g(\phi(1)) = \phi (1)$ and define it to swap $\phi (2k)$ and $\phi (2k+1)$ for $k \in \mathbb{N}$. 
It is not difficult to see that $g_1$ and $g_2$ belong to  $W(\mathcal{G},d_{\mathcal{G}})$, and that for any $v \in V$, $T(v)$ equals either $g_1(v)$ or $g_2(v)$. In particular, this implies $T(F_m) \subseteq g_1(F_m) \cup g_2(F_m)$ and therefore $T(F_m) \setminus F_m \subseteq (g_1(F_m) \setminus F_m) \cup (g_2(F_m) \setminus F_m)$. Since $T$ is injective, we have $|T(F_m)| = |F_m|$ and therefore $|F_m \setminus T(F_m)| = |T(F_m) \setminus F_m|$. Hence, 
\begin{align*}
\frac{|\sigma (E_m) \triangle E_m|}{|E_m|} &= \frac{|\sigma \phi^{-1}(F_m) \triangle \phi^{-1}(F_m)|}{|\phi^{-1}(F_m)|} \\
&= \frac{| \phi^{-1}(F_m) \triangle \sigma^{-1} \phi^{-1}F_m|}{|F_m|} \\
&= \frac{| F_m \triangle \phi\sigma^{-1} \phi^{-1}F_m|}{|F_m|} \\
&= \frac{|F_m \triangle T(F_m)|}{|F_m|} \\
&= \frac{2|T(F_m) \setminus F_m|}{|F_m|} \\
&\le \frac{2(|g_1(F_m) \setminus F_m) \cup (g_2(F_m) \setminus F_m)|}{|F_m|} \\
&\le \frac{2|g_1(F_m) \setminus F_m|}{|F_m|} + \frac{2|g_2(F_m) \setminus F_m|}{|F_m|} \\
&\le \frac{|g_1(F_m) \triangle F_m|}{|F_m|} + \frac{|g_2(F_m) \triangle F_m|}{|F_m|} 
\end{align*}
which converges to $0$ as $m \to \infty$ because $g_1, g_2 \in W(\mathcal{G}, d_{\mathcal{G}})$ and $W(\mathcal{G}, d_{\mathcal{G}})$ acts amenably on $\mathcal{G}$.
\end{proof}

A very basic example of a one-ended graph in the realm of groups is the Cayley graph $\mathbb{Z}^2$ with the word metric coming from its standard generating set $S = \{(\pm 1, 0), (0, \pm 1)\}$. The word metric is just the taxi-cab metric. The taxi-cab metric, Euclidean metric, and the Chebyshev metric defined as $d_{\infty}((x_1,y_1),(x_2,y_2)) = \max (|x_1-x_2|,|y_1-y_2|)$ all come from norms on a finite dimensional vector space and therefore are all bi-Lipschitz equivalent. Consequently, these metrics give the same wobbling group. For our purposes, we will consider the Chebyshev metric. 

The metric balls centered at the origin of radius $m \in \mathbb{N}$ are the squares 
\[
F_m := [-m,m] \times [-m,m].
\]
Hence, the most natural enumeration is the square spiral enumeration. That is, $\phi (1) = (0,0)$, $\phi (2) = (1,0)$, $\phi (3) = (1,1)$, $\phi (4) = (0,1)$, $\phi (5) = (-1,1)$, etc. Hence, if $\phi (n) = (i,j)$, given how we have enumerated $\mathbb{Z}^2$, the only possibilities for $\phi (n+1)$ are $(i+1,j)$, $(i,j+1)$, $(i-1,j)$, or $(i,j-1)$. In any case, it is clear that $d_{\infty}(\phi (n),\phi(n+1)) = 1$ for every $n \in \mathbb{N}$. Moreover, 
\[
\phi^{-1}(F_m) = \{1, \dots , (2m+1)^2\}.
\]
For $g \in W(\mathbb{Z}^2)$ with wobbling constant $C \in \mathbb{N}$, we have  $g(F_m) \setminus F_m \subset F_{m+C} \setminus F_m$, and hence 
\[
|g(F_m) \setminus F_m| \le |F_{m+C} \setminus F_m| = |F_{m+C}| - |F_m| = 8mC + 4C^2 + 4C.
\]
Hence, 
\[
\frac{|g(F_m) \triangle F_m|}{|F_m|} \le \frac{2(8mC + 4C^2 + 4C)}{(2m+1)^2} \to 0
\]
as $m \to \infty$, which shows $\{F_m\}_{m \ge 1}$ are F\o lner sets for $W(\mathbb{Z}^2)$ acting on $\mathbb{Z}^2$. The preimages $\phi^{-1}(F_m ) = \{1, \dots , (2m+1)^2\}$ are single, contiguous intervals which are F\o lner sets for $\sigma$ by our characterization in Section \ref{sec:ResidFiniteGrp}. Hence, $\phi^{-1} W(\mathbb{Z}^2) \phi$ with the spiral enumeration is shift-similar and acts amenably on $\mathbb{N}$ with $\sigma$.

For non-group example, take the non-negative ray $[0,\infty)$ and  glue a full binary tree $T_m$ of depth $m$ at each natural number $m$. This is certainly a one-ended graph. Define $\phi$ so that each tree and its root is enumerated by disjoint consecutive intervals. This will ensure $\phi$ is a bounded-step enumeration and hence $\phi^{-1} W(\mathcal{G})\phi$ is shift-similar. Let $F_m$ be the set of all vertices in $T_m$. Then $|F_m| = 1 + 2 \dots + 2^m = 2^{m+1}-1$, and $|\partial F_m| = 2$ for $m \ge 2$ because the only boundary vertices are $m-1$ and $m+1$. Hence, 
\[
\frac{|\partial F_m|}{|F_m|}  = \frac{2}{2^{m+1}-1} \to 0
\]
as $m \to \infty$, meaning the graph is amenable which is equivalent to the wobbling group acting amenably on it which is in turn equivalent to $\phi^{-1} W(\mathcal{G}, d) \phi$ and $\sigma$ acting amenably on $\mathbb{N}$. 

Not only does Theorem \ref{thm:OneEndedShiftSim} produce a plethora of shift-similar groups, but it provides us with an alternative construction of a group $G$ that is the shift-similar closure of a property (T) group such that $L(G)$ and $L(H_n(G))$ are both property Gamma factors. This is in the done in the following theorem by adapting ideas from \cite{juschenko15}. 

\begin{theorem}\label{thm:Wobbling(T)Amean}
    There exists a shift-similar wobbling group that acts jointly amenably on $\mathbb N$ and contains a property (T) subgroup. 
\end{theorem}

\begin{proof}

First, let $\Gamma$ be any solvable group with one end and containing a free subsemigroup. For example, we can take $\Gamma$ to be the Baumslag-Solitar group $BS(1,2) = \langle a,b | bab^{-1} = a^2 \rangle$ with $\{a,ab\}$ generating a free subsemigroup. Since $\Gamma$ contains a free subsemigroup, a binary tree coarsely embeds into $\Gamma$. Since $\text{SL}_3(\mathbb{Z})$ has finite asymptotic dimension, it coarsely embeds into a finite product of binary trees. Hence, $SL_3(\mathbb{Z})$ coarsely embeds into $\Gamma^k$ for some $k \in \mathbb{N}$. By \cite[Lemma 4.2]{juschenko15}, there exists some finite group $F$ such that $W(SL_3(\mathbb{Z}))$ embeds into $W(\Gamma^k \times F)$, which implies $W(\Gamma^k \times F)$ contains a copy of $SL_3(\mathbb{Z})$. 

Recall that the number of ends of a group is a quasi-isometric invariant, that the $k$-fold direct product of a one-ended group is one-ended, and that a direct product with a finite group produces a quasi-isometric group.  Putting this all together, $\Gamma^k \times F$ also has one end. Hence, by Theorem \ref{thm:OneEndedShiftSim} above, there exists a bijection $\phi : \mathbb{N} \to \Gamma^k \times F$ such that $\phi^{-1} W(\Gamma^k \times F)\phi$ in $\sym (\mathbb{N})$ is shift-similar. Since $\Gamma^k \times F$ is amenable, the wobbling group $W(\Gamma^k \times F)$ acts amenably on the Cayley graph and therefore $\phi^{-1}W(\Gamma^k \times F) \phi$ acts amenably on $\mathbb{N}$. Not only that, Theorem \ref{thm:OneEndedShiftSim} also tells us it acts jointly amenably with $\sigma$. 
\end{proof}

\begin{observation}
    In the setting of the previous theorem, define $G$ to be the shift-similar closure of the copy of $SL_3(\mathbb{Z})$ contained in $\phi^{-1} W(\Gamma^k \times F) \phi$. Then $G$ is shift-similar by construction and acts jointly amenably with $\sigma$, being a subgroup of $\phi^{-1} W(\Gamma^k \times F) \phi$. Hence, both $L(G)$ and $L(H_n(G))$ are property Gamma factors.
\end{observation}

\subsection{Wobbling Groups of Two-ended Graphs}\label{Sec:TwoEnds}

The condition that 
\[
\sup_{m \in \mathbb{N}} d_{S}(\phi (m), \phi (m+1)) < \infty
\]
imposes a strong topological restriction on the number of ends of the graph $\mathcal{G}$---namely, $\mathcal{G}$ must have one end. One might wonder whether the weaker condition of $T$ being bornologous puts any restrictions on the topological (end) structure of a group or graph. We will see that this is not the case, that there are examples of groups and graphs for which $T$ is bornologous and have two-ends.

One could consider graphs with more ends, such as a graph with infinitely many ends like the Cayley graph of a free group. But infinite ended groups, in particular, are always non-amenable and hence the wobbling group does not act amenably on the Cayley graph, which means our machinery will not necessarily produce property Gamma factors. For our purposes, we just consider one-ended and two-ended graphs, especially because the two-ended case nicely follows from the one-ended case, and leave graphs with more than two ends to future work.

First, let us consider $\mathbb{Z}$. With the standard generating set, the word metric is just the Euclidean metric. Let $\phi : \mathbb{N} \to \mathbb{Z}$ be defined by $\phi (n) = \frac{n}{2}$ if $n$ is even and $\phi (n) = \frac{1-n}{2}$ if $n$ is odd. First, note that $\mathbb{Z}$ has two ends, so $\phi$ cannot be a bounded-step enumeration. Indeed,
\[
|\phi (2n)-\phi (2n+1)| = 2n \to \infty
\]
as $n \to \infty$. However, $T$ is still bornologous, which is easy to verify by cases based on parity. Hence, using this bijection $\phi : \mathbb{N} \to \mathbb{Z}$, $T$ is bornologous which is equivalent to $\phi^{-1} W(\mathbb{Z},d_{euc}) \phi$ being shift-similar. 

The intervals $F_m = [-m,m]$ are F\o lner sets for the wobbling group $W(\mathbb{Z}, d_{euc})$ acting on $\mathbb{Z}$, and the calculations show these are similar to the ones for $\mathbb{Z}^2$. Note the inverse of $\phi$ is given by $\phi^{-1}(n) = 2n - 1$ if $n > 0$ and $\phi^{-1}(n) = 2(1-n)$ if $n \le 0$. It is not hard to verify that $\phi^{-1}(F_m) = [1, 2m+1]$, which is a single contiguous intervals implying $\{\phi^{-1}(F_n)\}_{m \ge 1}$ are F\o lner sets for $\phi^{-1} W(\mathbb{Z}, d_{euc}) \phi$ and $\sigma$ acting on $\mathbb{N}$.

Before we go further, we note that although $W(\mathbb{N}, d_{euc})$ embeds into $W(\mathbb{Z},d_{euc})$ in a canonical way, these two wobbling groups cannot be isomorphic. Indeed, $W(\mathbb{Z},d_{euc})$ admits the homomorphism onto the integers given by
\begin{align*}
f &\mapsto |\{x \ge 0 : f(x) < 0 \}| - |\{x < 0 : f(x) \ge 0\}| \\
&= \sum_{x \in \mathbb{Z}} (\chi_{\mathbb{Z}_{<0}} (f(x)) - \chi_{\mathbb{Z}_{<0}} (x)),
\end{align*}
which measures the net flow across $0$. Whereas $W(\mathbb{N},d_{euc})$ does not have a non-trivial homomorphism to $\mathbb{Z}$ because it is generated by its involutions. In fact, we have the following theorem. 

\begin{theorem}
Every permutation in $W(\mathbb{N},d_{euc})$ can be written as a product of at most three involutions. Consequently, $W(\mathbb{N},d_{euc})$ is not isomorphic to $W(\mathbb{Z},d_{euc})$. 
\end{theorem}
\begin{proof}
Let $g \in W(\mathbb{N},d_{euc})$. Then there exists $R \in \mathbb{N}$ such that $|g(x) - x| \le R$ for every $x \in \mathbb{N}$. The trick is to construct intervals $I_k$ of equal size partitioning $\mathbb{N}$ (even uniformly bounded is enough) and an involution $\tau \in W(\mathbb{N},d_{euc})$ such that $h := \tau g$ restricts to a permutation on each finite set $I_k$, and use the fact that a permutation of finite set can be written as a product of two involutions.

To this end, let $c_0 = 0$ and $c_k = 2Rk$ for $k \in \mathbb{N}$, and define $I_k = \{c_{k-1}+1, \dots , c_k\}$. Also, define $F_k  = \{1, \dots , c_k\}$ and 
\[
A_k = g(F_k) \setminus F_k \quad \text{ and } \quad B_k = F_k \setminus g(F_k) 
\]
Bounded displacement of $g$ implies 
\[
A_k \subseteq \{c_k+1, \dots , c_k + R \} \quad \text{ and } \quad B_k \subseteq \{c_k - R + 1, \dots , c_k\},
\]
and both are contained in 
\[
J_k = \{c_k-R+1, \dots , c_k+R\}
\]
Since $A_k$ and $B_k$ are disjoint with the same cardinality, we can pick any bijection from $A_k$ to $B_k$ that swaps pairs of points and extend it to a permutation of $\mathbb{N}$ by letting it act by identity outside of $A_k \cup B_k$. Call it $\tau_k$. Then $\{\tau_k\}_{k \ge 1}$ are involutions with disjoint support and therefore pairwise commute, meaning
\[
\tau = \tau_1 \tau_2 \dots 
\]
\looseness=-1 is a well-defined involution, and it is not difficult to show $|\tau (x)-x| \le 2R-1$ for every $x \in \mathbb{N}$.

Since $\tau_{\ell}(F_k) = F_k$ for $\ell \neq k$ and the factors of $\tau$ commute, we have
\[
\tau g(F_k) = \tau_k ((g(F_k) \cap F_k) \sqcup A_k) = \tau_k((g(F_k) \cap F_k)) \sqcup \tau_k(A_k) = g(F_k) \cap F_k) \sqcup B_k = F_k
\]
so
\[
\tau g(I_k) = \tau g(F_k \setminus F_{k-1}) = \tau g(F_k) \setminus \tau g(F_{k-1}) = F_k \setminus F_{k-1} = I_k.
\]
Hence, $h := \tau g$ restricts to a permutation on each $I_k$, implying $h|_{I_k} = a_k b_k$ for involutions $a_k,b_k \in \sym (I_k)$ which can be extended to $\mathbb{N}$ using the same notation. If $a = a_1 a_2 \dots $ and $b =b_1 b_2 \dots $, then $h = ab$, and it is not difficult to argue that $a,b \in W(\mathbb{N},d_{euc})$ with displacement constants $2R-1$ like $\tau$. Hence, $\tau g = ab$ or $g = \tau ab$, where $\tau a, b \in W(\mathbb{N},d_{euc})$ are involutions. 
\end{proof}

For two-ended groups, the integers represent the entire geometric story in terms of these wobbling groups. If $G$ is any two-ended group, Stallings' theorem tells us that $G$ is virtually cyclic. Pick an infinite cyclic normal subgroup $H=\langle z\rangle$ of index $m \in \mathbb{N}$, take right coset representatives $g_0,\dots,g_{m-1}$, use \(\{z^{\pm1},g_i^{\pm1}\}\) as a finite generating set. Then it is not hard to argue by cases that $mk+i\mapsto z^k g_i$ is a bi-Lipschitz map from $\mathbb{Z} \to G$, so they have the same wobbling group.

However, there are locally finite, connected two-ended graphs which are not bi-Lipschitz equivalent to the integers. The following is an example of such a two-ended graph.

\begin{example}[Integers-by-Finite Binary Tree]\label{ex:TwoEndsNotZ}
Let $\mathcal{G}$ be the graph constructed by taking the integers with an edge between consecutive integers, and gluing the root of a complete binary tree $T_n$ of depth $n$ to each natural number $n \in \mathbb{N}$. Note $|T_n| = 2^{n+1}-1$ and $\mathrm{diam}(T_n) = 2n$.
\end{example}

Although this example certainly contains an isometric copy of $\mathbb{Z}$, it is not bi-Lipschitz equivalent to $\mathbb{Z}$, and the wobbling group does not admit a shift-similar embedding because it lacks the ``right amount" of symmetry.  

\begin{proposition}\label{prop:TwoEndsNotBLEqivToZ}
    Let $\mathcal{G} = (V,E)$ be the integers-by-finite binary tree from Example \ref{ex:TwoEndsNotZ}. Then $\mathcal{G}$ is not bi-Lipschitz equivalent to $\mathbb Z$ and it's wobbling group does not admit a shift-similar embedding.
\end{proposition}

\begin{proof}
First, we prove $\mathcal{G} = (V,E)$ is not bi-Lipschitz equivalent to $\mathbb{Z}$. By way of contradiction, suppose there exists an injective Lipschitz map $f : V \to \mathbb{Z}$. Then there exists $S > 0$ such that whenever $u, v \in V$ are adjacent, we have $|f(u)-f(v)| \le S$. Let $u,v$ be arbitrary vertices in $T_n$. Then $d(u,v) \le \diam (T_n) = 2n$, so there exists a shortest path $u=v_0,v_1, \dots , v_k = v$ in $T_n$ whose length is at most $2n$. Hence,
\[
|f(u)-f(v)| \le \sum_{i=1}^{k} |f(v_i)-f(v_{i-1})| \le k S \le 2nS
\]
and therefore
\[
\diam (f(T_n)) = \max_{u,v \in T_n} |f(u)-f(v)| \le 2n S.
\]
Now, recall that for any finite set $A \subseteq \mathbb{Z}$, we have $\diam (A) = \max A - \min A \ge |A|-1$. Hence, on the other hand, we have
\[
\diam(f(T_n)) \ge |f(T_n)|-1 = |T_n| - 1 \ge 2^{n+1} -2. 
\]
Combing the inequalities gives us $2^{n+1}-2\le 2nS$ or
\[
S \ge \frac{2^{n+1}-2}{2n} \to \infty
\]
as $n \to \infty$, which is a contradiction. 

Now, we prove the second claim. By way of contradiction, suppose there is an enumeration $\phi$ of the vertices for which $T$ is bornologous. Then there exists $S \in \mathbb{N}$ such that $d(T(u),T(v)) \le S$ whenever $u$ and $v$ are adjacent. Pick $N \in \mathbb{N}$ large enough so that $\phi (1)$ belongs to the ``left half" of $\mathcal{G}$. That is, it belongs to
\[
L := \{\dots , N-1, N\} \cup \bigcup_{1 \le n \le N}  T_n, 
\]
the left half. Define the right half as
\[
R := \{N+1,N+2, \dots \} \cup \bigcup_{n > N } T_n. 
\]
First, let us prove the number of transitions or crossings from $R$ to $L$ via $T$ is infinite. To do this, recall that if $A \subset \mathbb{N}$ is infinite and co-infinite (i.e., $\mathbb{N} \setminus A$ is infinite), then there are infinitely many $n \in \mathbb{N}$ with $n \in A$ but $n+1 \notin A$, and also infinitely many $n$ with $n \notin A$ but $n+1 \in A$. 

Since $L$ is infinite and co-infinite in the graph, $\phi^{-1}(L)$ is an infinite and co-infinite subset of $\mathbb{N}$. Hence, there exist infinitely many $n \in \mathbb{N}$ with $n \in \phi^{-1}(L)$ and $n+1 \notin \phi^{-1}(L)$, and hence infinitely many $n \in \mathbb{N}$ with $\phi (n) \in L$ and $\phi (n+1) \in R$. If $n \in \mathbb{N}$ is anyone of these infinitely many natural numbers, set $v = \phi (n+1) \in R$. Hence,
\[
T(v) = T(\phi(n+1)) = \phi \sigma^{-1}\phi^{-1}(\phi(n+1)) = \phi (n) \in L
\]
Hence, there are infinitely many $v \in R$ such that $T(v) \in L$. 

Now, given $r > 0$, let $v$ be a vertex in a tree $T_n$ of depth $n \ge r$. We claim that 
\[
|B(v,r) | \ge 2^{\floor{r/2}+1}
- 1. 
\]
Let $d$ be the depth of $v$ in $T_n$. First, if $d \le \floor{r/2}$, then the subtree rooted at $v$ has depth at least
\[
n-d \ge r - \floor{r/2} \ge \floor{r/2},
\]
meaning $B(v,r)$ can accommodate a full binary tree of depth $\floor{r/2}$, of which there are $2^{\floor{r/2}+1}
- 1$ many vertices, rooted at $v$. If $d > \floor{r/2}$, then move $\floor{r/2}$ units from $v$ towards the root $n$ to an ancestor $u$ within the tree $T_n$. The subtree rooted at $u$ has depth 
\[
n - \left(d-\floor{\frac{r}{2}} \right) = n-d + \floor{\frac{r}{2}} \ge \floor{\frac{r}{2}}. 
\]
This shows that $u$ roots a full binary tree of depth at least $\floor{r/2}$. If $w$ is any vertex in that subtree rooted at $u$, then 
\[
d(v,w) \le d(v,u) + d(u,w) \le \floor{\frac{r}{2}} + \floor{\frac{r}{2}} \le r
\]
which implies $w \in B(v,r)$. Therefore, the whole subtree rooted at $u$ of depth $\floor{r/2}$ is contained in $B(v,r)$. In either case, we can conclude 
\[
|B(v,r)| \ge 2^{\floor{r/2}+1}-1. 
\]

Pick $r > 0$ large enough so that 
\[
2^{\floor{r/2}+1}-1 > 2Sr+1
\]
holds. Recall, there are infinitely many $v \in R$ for which $T(v) \in L$. Let $Y$ denote the set of all such vertices. Since $T$ is injective, $T(Y)$ is an infinite subset of $L$. Since the left hand side $L$ consists of only finitely many trees (or possibly none) with finitely many vertices, the set $T(Y)$ contains infinitely many vertices from the left integer spine $\{\dots , N-1, N \}$. Hence, $T(Y)$ is unbounded to the left, meaning that for every $M \in \mathbb{Z}$, there exists $v \in Y$ such that $T(v) < M$. In particular, we can pick $v \in Y$ such that 
\[
T(v) + Sr < 1. 
\]
This implies that the metric ball $B(T(v),Sr)$ lies entirely to the left of all the attached trees. Indeed, if $B(T(v),Sr)$ contained a vertex from some tree, path connectedness of $B(T(v),Sr)$ implies it contains $1$ since any path from $T(v)$ to a vertex in a tree passes through $1$. Hence,
\[
Sr \ge d(T(v), 1) = 1 - T(v) \implies T(v) + Sr \ge 1,
\]
which is a contradiction. Hence, $B(T(v),Sr)$  is just an interval on the integer spine and, therefore,
\[
|B(T(v),Sr)| = 2Sr+1.
\]
Now, on the right hand side $R$, only finitely many vertices lie in trees $T_n$ with $n < r$. Because $Y$ is infinite, we can discard those finitely many vertices, and pick $v \in Y$ with $v \in T_n$ for some $n \ge r$, and therefore the metric ball $B(v,r)$ contains a complete binary tree of depth $\floor{r/2}$. Moreover, among the infinitely many $v \in Y$ with $T(v)$ arbitrarily far to the left, only finitely many fail to lie in a tree $T_n$ with $n \ge r$. Therefore, we can pick $v \in Y$ satisfying both $v \in T_n$ with $n \ge r$ and $T(v) + Sr < 1$. 

Finally, since $T$ is an injective Lipschitz map, we have $T(B(v,r)) \subset B(T(v),Sr)$. Hence, 
\[
2Sr+1 = |B(T(v),Sr)| \ge |T(B(v,r))| = |B(v,r)| \ge 2^{\floor{r/2}+1}-1,
\]
which contradicts how we picked $r$.
\end{proof}

The lack of symmetry in the above example is precisely what prevents the wobbling group from having a shift-similar embedding. If we symmetrize the graph by adding binary trees at the negative integers of the same depth as their positive counterparts, we can construct an enumeration $\phi : \mathbb{N} \to V$ for which $T$ is bornologous. The key is that this graph possesses a symmetry---namely, an involution swapping the ends. Specifically, the vertices in this graph can be described as $(z,v)$ for $z \in \mathbb{Z}$ and $v \in T_z$, where $v$ is a finite word in $\{0,1\}^{<|z|}$. The map \((z,v) \mapsto (-z,v)\) then defines an involution swapping the ends. More generally, we only need the graph to be "coarsely" symmetric, meaning the two-ended graph admits a bijective coarse involution swapping its ends.

In fact, we can prove something considerably more general about metric spaces admitting a coarse involution,
which will have implications for two-ended graphs. A bornologous map $\iota : X \to X$ is a \emph{coarse involution} if 
\[
\sup_{x \in X} d(\iota^2(x),x) < \infty.
\]

\begin{definition}
    We say a countable metric space $(X,d)$ is a \emph{coarse symmetric space} if $X$ can be 
    decomposed as $X = X_- \sqcup X_+ \sqcup F$, where $F$ is some finite region, and $X$ admits a bijective, coarse involution $\iota : X \to X$ such that $r(X_+) = X_-$.
\end{definition}

We now prove the wobbling group of a coarse symmetric space admits a shift-similar embedding 

\begin{theorem}\label{thm:CoarseInvoBornolo}
Let $(X,d)$ be a coarse symmetric space. If $X_+$ admits a bounded-step enumeration, then there exists a bijection $\phi : \mathbb{N} \to X$ such that $T$ is bornologous with respect to $d$.

\end{theorem}
\begin{proof}
Let $N:= |F|$. Let $\phi$ enumerate the points in $F$ in any arbitrary manner using the natural numbers in $\{1, \dots , N\}$. We use the numbers $\{N+1, N+2, \dots \}$ to enumerate $X_-$ and $X_+$. Every element in this set is of the form $N+2n$ or $N+2n-1$ for some $n \in \mathbb{N}$. Let $\iota:X \to X$ be a coarse involution for $X$. Define $\phi : \mathbb{N} \to X$ by $\phi (N+2n) = a_n$ and $\phi (N+2n-1) = \iota (a_n)$ for $n \in \mathbb{N}$, which certainly defines a global enumeration of $X$. We will show $T$ and $\iota $ are close on the cofinite set $X \setminus F$. First, let $x \in X_+$. Then $x = a_n = \phi (N+2n)$ for some $n \in \mathbb{N}$. Hence,
\begin{align*}
d(\iota (x),T(x)) &= d(\iota (a_n), T(a_n)) \\
&= d(\iota (a_n), \phi (\phi^{-1}\phi (N+2n)-1)) \\
& = d(\iota (a_n),\phi (N+2n-1)) \\
&= d(\iota (a_n), \iota (a_n)) \\
&= 0.
\end{align*}
This actually shows that $T$ and $\iota$ agree on the right hand side of $X$. Now, let $x \in X_- = r(X_+)$. Then $x = r(a_n) = \phi (N+2n-1)$ for some $n \in \mathbb{N}$. Hence, 
\begin{align*}
d(\iota (x),T(x)) &= d(\iota^2(a_n), T(r(a_n))) \\
&= d(\iota^2(a_n), \phi (\phi^{-1} \phi (N+2n-1)-1)) \\
&= d(\iota^2(a_n),\phi (N+2(n-1)) \\
&= d(\iota^2(a_n),a_{n-1}) \\
& \le d(\iota^2(a_n),a_n)) + d(a_n, a_{n-1}) \\
& \le \sup_{x \in X} d(\iota^2(x),x) + \sup_{n \in \mathbb{N}} d(a_n,a_{n+1}).
\end{align*}
Hence, $T$ is close to $\iota$ on a cofinite set in $X$ and hence close to it on all of $X$. Since $\iota$ is bornologous, it follows $T$ must also be bornologous. Hence, $\phi W(X,d) \phi^{-1}$ is shift-similar. 
\end{proof}

Let $\iota : X \to X$ is a coarse involution of a metric space, and define 
\[
C := \sup_{x \in X} d(\iota^2(x),x ) < \infty.
\]
Because $d(\iota^2(x),x) = d(\iota (\iota (x)),x)$, it is clear that $\iota$ is its own coarse inverse. Let $B \subseteq X$ be a bounded subset. Because $\iota$ is bornologous, $\iota (B)$ is also bounded. Now, if $x \in \iota^{-1}(B)$, then $\iota^2(x) \in \iota (B)$, and hence
\[
d(x, \iota (B)) \le d(x,\iota^2(x)) \le C,
\]
which proves $\iota^{-1}(B) \subseteq N_C(\iota (B))$. Because the neighborhood of a bounded set is bounded, it follows $\iota^{-1}(B)$ is bounded. This proves that $\iota$ is coarsely proper. 

Hence, a coarse involution is necessarily a coarse equivalence. If $\iota$ is also a coarse involution of a graph, then we know it induces a bijection on the ends of the graph. Also, if $\iota$ is, moreover, a bijection, then for every $x \in X$, we have
\[
d(\iota^{-1}(y), \iota(x)) = d(\iota^{-1}(x), \iota^2(\iota^{-1}(x)) \le C,
\]
which shows $\iota^{-1}$ is bornologous because it is close to $\iota$ which is bornologous.

Our next theorem is about two-ended graphs which admit a bijective coarse involution.

\begin{theorem}\label{thm:TwoEndsShiftSim}
Let $\mathcal{G} = (V,E)$ be a locally finite, connected graph with two-ends. Suppose it admits a bijective, coarse involution $\iota : V \to V$ on the vertices which swaps the two ends of $\mathcal{G}$. Then $\mathcal{G}$ admits an enumeration $\phi : \mathbb{N} \to V$ such that $T$ is bornologous. If $\mathcal{G}$ is additionally amenable and the vertices have uniformly bounded degree, then $\phi^{-1} W(\mathcal{G}) \phi$ and $\sigma$ act amenably on $\mathbb{N}$.
\end{theorem}
\begin{proof}
Let $\iota$ be a bijective, coarse involution swapping the ends of $\mathcal{G}$. Since $\mathcal{G}$ has two ends, there exists a finite set $F_0 \subset V$ such that $\mathcal{G} \setminus F_0$ splits into two infinite components and possibly finitely many finite components. We can enlarge $F_0$ to absorb the vertices of the finitely many finite components of $\mathcal{G} \setminus F_0$. Doing so allows us to partition $V \setminus F_0$ into two sets of vertices $A$ and $B$ with the induced graphs $\mathcal{G}[A]$ and $\mathcal{G}[B]$ being the two infinite components of $\mathcal{G} \setminus F_0$ without any finite components. 

Because $\iota$ is bornologous, we can find $S > 0$ such that whenever $u$ and $v$ are adjacent vertices, we have $d(\iota (u), \iota (v)) \le S$. Define the finite set of vertices $E$ as 
\[
E = F_0 \cup \iota^{-1}(N_S(F_0)).
\]
Since the graph spanned by $A$ is infinite with one end and the graph spanned by $E$ is finite, the graph $\mathcal{G}[A] \setminus E$ has exactly one infinite component. The same is true of $\mathcal{G}[B] \setminus E$. Let $V_1$ and $V_2$ denote the vertices of this unique infinite components in $\mathcal{G}[A] \setminus E$ and $\mathcal{G}[B] \setminus E$, respectively. Note that both $A \setminus V_1$ and $B \setminus V_2$ are finite (possibly empty). Indeed, in $\mathcal{G}[A]\setminus E$, $V_1$ is the set of vertices of the unique infinite component and all other components are finite. Therefore,
\[
A\setminus V_1=(A\cap E)\cup F_A,
\]
where $F_A$ is the union of the vertices in the finitely many finite components of $\mathcal{G}[A]\setminus E$. Since $E$ is finite, $A\cap E$ is finite, so $A\setminus V_1$ is finite. The same argument applies to $B\setminus V_2$.

We claim that $\iota (V_1) \subseteq B$ and $\iota (V_2) \subseteq A$. By symmetry, we only prove the former inclusion. Let $v \in V_1$. Since $\mathcal{G}[V_1]$ is infinite and connected, there is a ray $(v_i)_{i \ge 1}$ in $\mathcal{G}[V_1]$ starting at $v$. For every $i \in \mathbb{N}$, $v_i \notin E$, so in particular $\iota (v_i) \notin N_S(F_0)$ and hence $d(\iota (v_i), F) > S$. Since $v_i$ and $v_{i+1}$ are adjacent, we have 
\[
d(\iota (v_i), \iota (v_{i+1})) \le S.
\]
\looseness=-1 Hence, any shortest path connecting $\iota (v_i)$ and $\iota (v_{i+1})$ must have length at most $S$, and it must avoid $F_0$ otherwise we would have both $d(\iota (v_i), F_0) \le S$ and $d(\iota (v_{i+1}), F_0) \le S$. This implies there is a path from $\iota (v_i)$ to $\iota (v_{i+1})$ avoiding $F_0$ for every $i \in \mathbb{N}$, meaning that all the vertices of the proper coarse ray $(\iota (v_i))_{i \ge 1}$ lie in the same necessarily infinite component of $\mathcal{G} \setminus F_0$.

Note that the ray $(v_i)_{i \ge 1}$ is chosen inside $\mathcal{G}[V_1]$, and since $V_1 \subseteq A$, the entire ray lies in $\mathcal{G}[A]$. Moreover, $(v_i)_{i \ge 1}$ represents an end of $\mathcal{G}$ living in the infinite component $\mathcal{G}[A]$. Since $\iota$ swaps the two ends, the coarse ray $(\iota(v_i))_{i \ge 1}$ must eventually lie in the component $\mathcal{G}[B]$; otherwise $\iota$ would not actually swap the ends. But all the vertices of $(\iota(v_i))_{i \ge 1}$ lie in the same infinite component of $\mathcal{G}\setminus F_0$. Since a tail of $(\iota(v_i))_{i \ge 1}$ lies in the infinite component $\mathcal{G}[B]$, the entire coarse ray $(\iota(v_i))_{i \ge 1}$ lies in $\mathcal{G}[B]$. In particular, $\iota(v) \in B$.

Next, we prove $\iota (V_1)$ is cofinite in $B$. To this end, let $v \in B \setminus \iota (V_1)$, and write $u = \iota^{-1}(v)$. Then clearly neither $u \in V_1$, nor can $u \in V_2$ because otherwise 
\[
v = \iota (u) \in \iota (V_2) \subseteq A
\]
contrary to $v \in B$. Because of our partition of our vertices, this leaves us with 
\[
u \in F_0 \cup (A \setminus V_1) \cup (B \setminus V_2)
\]
and therefore 
\[
B \setminus \iota (V_1) \subseteq \iota (F_0 \cup (A \setminus V_1) \cup (B \setminus V_2)),
\]
which is finite. Because $V_1 \subseteq A$ and $\iota (V_1) \subseteq B$, these sets are disjoint, and we have that
\[
F := V \setminus (V_1 \sqcup \iota (V_1))
\]
is finite. 

Now, define $V_+ := V_1$ and $V_- := \iota (V_1) = \iota (V_+)$. Since $\mathcal{G}[V_+]$ is a locally finite, connected graph with one-end, Theorem \ref{thm:OneEndedShiftSim} says that there exists a bounded-step enumeration $a : \mathbb{N} \to V_+$. Since we have a bijective coarse involution $\iota$, a decomposition $V  = V_+ \sqcup V_- \sqcup F$ with $\iota (V_+) = V_-$, and a bounded-step enumeration of $V_+$, Theorem \ref{thm:CoarseInvoBornolo} tells us there exists a bijection $\phi : \mathbb{N} \to V$ such that $T$ is bornologous.

Now, further assume $\mathcal{G}$ is amenable and that the vertices have uniformly bounded degree. Let $\{H_m\}_{m \ge 1}$ be F\o lner sets for $\mathcal{G}$. For convenience, we can remove the points $v_0 := \phi (1)$ and $\iota^{-1}(v_0)$ from $H_m$ and still have F\o lner sets because 
\[
|\partial (H_m \setminus \{v_0, \iota^{-1}(v_0)\})| \le |\partial H_m| + 2D,
\]
where $D$ is the maximum degree, and the new sets will have cardinality at least $|H_m|-2$, and the ratio of the new sets still tends to $0$ because it is bounded above by
\[
\frac{|\partial H_m| + 2D}{|H_m| - 2} = \frac{|\partial H_m|}{|H_m|} \cdot \frac{|H_m|}{|H_m| - 2} + \frac{2D}{|H_m| - 2}
\]
which converges to $0$. 

Now, consider the sets $F_m := H_m \cup \iota (H_m)$. We claim these are also F\o lner sets. First, 
\[
\partial (H_m \cup \iota (H_m)) \subseteq \partial H_m \cup \partial (\iota (H_m)) 
\]
Let $w \in \partial (\iota (H_m))$. Then there exists $v \in H_m$ such that $w = \iota (v)$ and there exists $z \notin \iota (H_m)$ adjacent to $w = \iota (v)$. Define $u := \iota^{-1}(z)$. Then $u \notin H_m$. Since $\iota^{-1}$ is bornologous, there exists $R > 0$ such that 
\[
R \ge d(\iota^{-1}(w), \iota^{-1}(z)) = d(v,u)
\]
because $w$ and $z$ are adjacent. Since $d(v,u) \le R$, there is a path from $v \in H_m$ to $u \notin H_m$ of length at most $R$. Let $b$ be the last vertex in $H_m$ that lies along this path. Then the next vertex on the path is outside $H_m$, so $b \in\partial H_m$. Hence,  $b \in \partial H_m$ with $d(v,b) \le R$. Now, since $\iota$ is bornologous, given $R > 0$, there exists $S > 0$ such that 
\[
S \ge d(\iota (v), \iota(b)) = d(w,\iota (b))
\]
because $d(v,b) \le R$, implying that $w = \iota (v) \in N_{S}(\iota (\partial H_m))$. Since the graph  $\mathcal{G}$ has bounded geometry, given $S > 0$, there exists a constant $M(S) >0$ such that 
\[
|N_{S}(\iota (\partial H_m))| \le M(S) |\iota (\partial H_m)| = M(S) |\partial H_m|.
\]
Moreover, since $H_m \subseteq F_m$, we have $|F_m|^{-1} \le |H_m|^{-1}$. Putting this all together, we have
\begin{align*}
\frac{|\partial F_m|}{|F_m|} & \le \frac{|\partial (H_m \cup \iota (H_m))|}{|H_m|} \\
&\le \frac{|\partial H_m| + |\partial (\iota (H_m))|}{|H_m|} \\
&\le \frac{|\partial H_m| + |N_{S}(\iota (\partial H_m))|}{|H_m|} \\
&\le \frac{|\partial H_m| + M(S) |\partial H_m|}{|H_m|} \\
&\le \frac{(1+M(S))|\partial H_m|}{|H_m|} \to 0
\end{align*}
as $m \to \infty$. Now, let us show $\{\phi^{-1}(F_m)\}_{m \ge 1}$ are F\o lner sets for $\sigma$. First, we show $\{F_m\}_{m \ge 1}$ are F\o lner sets for the bijective coarse involution $\iota$.

Given a finite set of vertices $A$, if $v \notin A$ but $d(v,A) \le R$, then there exists $a \in A$ such that $d(v,a) \le R$. Hence, there exists a path from $v \notin A$ to $a \in A$ whose length is at most $R$. Because $a \in A$, the path eventually enters $A$. Let $v_1$ be the first vertex on the path just before it enters $A$. Hence, $v_1 \in \partial A$ and $d(v,v_1) \le R$, which proves
\[
\{v \notin A : d(v,A) \le R \} \subseteq \bigcup_{a \in \partial A} B(a,R) 
\]
and hence 
\[
|\{v \notin A : d(v,A) \le R \}| \le \left| \bigcup_{a \in \partial A} B(a,R)  \right| \le \sum_{a \in \partial A} |B(a,R)| \le M(R) |\partial A|. \tag{1}
\]

Define 
\[
C := \sup_{v \in V} d(v,\iota^{-2} (v)) = \sup_{v \in V} d(\iota^{2} (v), v) < \infty
\]
If $v \in \iota (F_m) \setminus F_m$, then $v \in \iota^2(H_m) \setminus H_m$. Hence, there exists $h \in H_m$ such that $v = \iota^2 (h)$ but $v \notin H_m$, and therefore $d(v,h) = d(v,\iota^{-2}(v)) \le C$. This shows that if $v \in \iota (F_m) \setminus F_m$, then $v \notin H_m$ and $d(v,H_m) \le C$. Hence, by (1), we have 
\[
|\iota (F_m) \setminus F_m| \le |\{v \notin H_m : d(v,H_m) \le C \}| \le M(C) |\partial H_m|
\]
and therefore 
\[
\frac{|\iota (F_m) \triangle F_m|}{|F_m|} \le 2 M(C) \frac{|\partial H_m|}{|H_m|} \to 0
\]
as $m \to \infty$. 

Using the fact that $T$ and $\iota$ are close (in the proof of Theorem \ref{thm:CoarseInvoBornolo}) we prove $\{F_m\}_{m \ge 1}$ are F\o lner sets for $T$. Define 
\[
L := \sup_{v \in V \setminus \{v_0\}} d(T(v),\iota (v)) < \infty.
\]
Let $w \in T(F_m) \setminus \iota (F_m)$. Then $w = T(v)$ for some $v \in F_m$. Note that 
\[
d(w,\iota (v)) = d(T(v),\iota (v)) \le L
\]
so $w$ lies outside of $\iota (F_m)$ but is within $L$ units of it. Choose a shortest path from $w=T(v)$ to $\iota (v) \in \iota (F_m)$, and let $u$  be the first vertex on the path just before it enters $\iota (F_m)$. Hence, $u \in \partial (\iota (F_m))$ and $d(w,u) \le L$. Hence, 
\[
T(F_m) \setminus \iota (F_m) \subseteq \bigcup_{u \in \partial (\iota (F_m))} B(u,L)
\]
and consequently
\[
|T(F_m) \setminus \iota (F_m)| \le M(L) |\partial (\iota (F_m))|.
\]
Since $\iota^{-1}$ is bornologous, there exists $S > 0$ such that $d(\iota^{-1}(u), \iota^{-1}(v)) \le S$ whenever $u$ and $v$ are adjacent. Using the same boundary argument, we have 
\[
|\partial (\iota (F_m))| \le M(S) |\partial F_m|.
\]
Putting this all together, we obtain 
\begin{align*}
\frac{|T(F_m) \triangle F_m|}{|F_m|} &\le \frac{|T(F_m) \triangle \iota (F_m)| + |\iota (F_m) \triangle F_m|}{|F_m|} \\
&\le 2 M(L) M(S) \frac{|\partial F_m|}{|F_m|} + \frac{|\iota (F_m) \triangle F_m|}{|F_m|} \to 0 \\
\end{align*}
as $m \to \infty$. Finally, 
\begin{align*}
\frac{|\sigma (\phi^{-1}(F_m)) \triangle \phi^{-1}(F_m)|}{|\phi^{-1}(F_m)|} &= \frac{|\phi^{-1}(F_m) \triangle \sigma^{-1} \phi^{-1}(F_m)|}{|F_m|} \\ 
&= \frac{|F_m \triangle \phi \sigma^{-1} \phi^{-1}(F_m)|}{|F_m|} \\
&= \frac{|F_m \triangle T(F_m)|}{|F_m|} \to 0
\end{align*}
as $m \to \infty$. 
\end{proof}

Examples satisfying Theorem \ref{thm:TwoEndsShiftSim} include the symmetrized binary tree example. More generally, let $\{\mathcal{G}_z\}_{z \in \mathbb{Z}}$ be a family of finite connected graphs with $\mathcal{G}_{-z} = \mathcal{G}_z$ for every $z \in \mathbb{Z}$ and $|V(\mathcal{G}_z)| \to \infty$ as $|z| \to \infty$, and glue $\mathcal{G}_z$ at each $z \in \mathbb{Z}$. It is clear that this admits a coarse involution swapping the ends. We can choose $\mathcal{G}_z$ so that the vertices in the entire graph have uniformly bounded degree. Define $F_m = \{0\} \cup \bigcup_{i=-m}^{m} V(\mathcal{G}_m)$ for $m \in \mathbb{N}$. It is easy to see that all these examples are amenable. Indeed, $|\partial F_m | = 2$ because the only boundary vertices are $-(m+1)$ and $m+1$. Hence,
\[
\frac{|\partial F_m|}{|F_m|} \le \frac{2}{|V(\mathcal{G}_m)|} \to 0
\]
as $m \to \infty$. 

If a two-ended graph does not admit a bijective coarse involution, then its wobbling group can fail to have a shift-similar embedding. Indeed, this is the case for the integers-by-finite binary tree from Example \ref{ex:TwoEndsNotZ}. In Proposition \ref{prop:TwoEndsNotBLEqivToZ} we proved that it is not bi-Lipschitz equivalent to $\mathbb{Z}$ and does not admit a shift-similar embedding. We now show this example does not admit an injective bornologous map that swaps the ends, let alone a bijective coarse involution that does so. 

\begin{proposition}\label{prop:TwoEndsNoSwapEnds}
     Let $\mathcal G$ be the integers-by-finite binary tree graph from Example \ref{ex:TwoEndsNotZ}. Then there does not exist an injective bornologous map swapping the ends.
\end{proposition}
   
\begin{proof}
To seek a contradiction, suppose $f$ is an injective bornologous (Lipschitz) map on the vertices. Then there exists $S \in \mathbb{N}$ such that $d(f(u),f(v)) \le S$ whenever $u$ and $v$ are adjacent. By construction, removing the vertex $0$ splits the graph into two disjoint infinite components: the left component $\mathbb{Z}_{<0}$, and the right component 
\[
\mathbb{Z}_{\ge 1} \cup \bigcup_{n \ge 1} T_n.
\]
It is clear that $(-n)_{n \ge 1}$ and $(n)_{n \ge 1}$ are the proper rays representing the left and right ends, respectively. Since $f$ swaps the ends of the graph, we know that $(f(n)))_{n \ge 1}$ and $(-n)_{n \ge 1}$ represent the same end. Because $f$ is injective Lipschitz, $(f(n))_{n \ge 1}$ is a coarse ray which is proper. Hence, for any $m \in \mathbb{Z}_{<0}$, the proper coarse ray $(f(n))_{n \ge 1}$ must escape $\{m, m+1 , \dots , 0\}$, meaning $f(n) \notin \{m, \dots ,0\}$ for all sufficiently large $n \in \mathbb{N}$. Since $(f(n))_{n \ge 1}$ and $(-n)_{n \ge 1}$ define the same end, there exists $k \in \mathbb{N}$ such that their $k$-th tail must lie in the same infinite component of $\mathcal{G} \setminus \{0\}$. This means we must have $f(n) < m$ for all sufficiently large $n \in \mathbb{N}$. That is, for every $m \in \mathbb{Z}_{<0}$, we have that $f(n) < m$ for all sufficiently large $n \in \mathbb{N}$, which is to say $f(n) \to - \infty$ as $n \to \infty$. 

\looseness=-1 Now, for all sufficiently large $n \in \mathbb{N}$, define $R_n = -f(n)$, so $R_n \to \infty$ as $n \to \infty$. Also, define 
\[
d_n = \min (n, \floor{ \frac{R_n}{2S}})
\]
Then $d_n \to \infty$ as $n \to \infty$. Moreover, 
\[
Sd_n \le S \frac{R_n}{2S} = \frac{R_n}{2} = \frac{-f(n)}{2}
\]
and hence for all sufficiently large $n \in \mathbb{N}$ we have
\[
f(n) + Sd_n  \le f(i) -  \frac{f(n)}{2} = \frac{f(n)}{2} < 1
\]
because $f(n) \le 0$ for all sufficiently large $n \in \mathbb{N}$. Hence, the metric ball $B(f(n),Sd_n)$ does not contain any vertex in any tree and therefore lies entirely to the left of them. If this were not so, then by the argument we gave before, $B(f(n),Sd_n)$ would contain $1$ which would contradict the inequality $f(n) + Sd_n < 1$. Hence, $B(f(n),Sd_n)$ is just a metric ball on the integer spine $\mathbb{Z}_{\le 0}$, meaning 
\[
|B(f(n),S d_n)| = 2Sd_n+1
\]

Now, define $T_n^{d_n}$ to be the subtree of $T_n$ rooted at $n$ of depth $d_n$ (this is well-defined because $d_n \le n$). Given a vertex $v$ in $T_n^{d_n}$, we have $d(n,v) \le d_n$, so the Lipschitz property gives us
\[
d(f(n),f(v)) \le Sd_n
\]
and hence $f(V(T_n^{d_n})) \subseteq B(f(n),Sd_n)$. Finally, 
\[
2Sd_n + 1 = |B(f(n),Sd_n)| \ge |f(V(T_n^{d_n}))| = |V(T_n^{d_n})|= 2^{d_n+1}-1
\]
which fails for all sufficiently large $d_n$. Hence, an injective bornologous (Lipschitz) swapping of the ends cannot exist. 
\end{proof}

However, the integers-by-finite binary tree is an amenable graph. Indeed, the set $F_n$ defined to be all the vertices from $-n$ to $n$, including all the vertices in the complete finite binary trees at the integer vertices, form F\o lner sets for the graph. Combing the amenability of the action, Proposition \ref{prop:TwoEndsNoSwapEnds}, and Lemma \ref{lemma:AmeanMetric/Wobbling} we get the following observation.

\begin{observation}
    The integers-by-finite binary tree two ended graph $\mathcal G$ from Example \ref{ex:TwoEndsNotZ} does not satisfy Theorem \ref{thm:TwoEndsShiftSim}, its wobbling group acts amenably on $\mathcal G$, and is not shift-similar. 
\end{observation}

\bibliographystyle{amsalpha} 
\bibliography{stable}     

\end{document}